\documentclass{article}
\usepackage{latexsym,amsmath,amssymb,amscd}
\usepackage[dvips]{color}
\begin{document}

\title{Linearized Analysis of Barotropic Perturbations around Spherically Symmetric Gaseous Stars governed by the Euler-Poisson Equations }
\author{Juhi Jang \footnote{Department of Mathematics, University of Southern California, and Korea Institute for Advanced Study, Seoul, Republic of Korea, E-mail: juhijang@usc.edu}   and Tetu Makino \footnote{Professor Emeritus at Yamaguchi University, Japan, E-mail: makino@yamaguchi-u.ac.jp}}
\date{\today}
\maketitle

\newtheorem{Lemma}{Lemma}
\newtheorem{Proposition}{Proposition}
\newtheorem{Theorem}{Theorem}
\newtheorem{Definition}{Definition}
\newtheorem{Remark}{Remark}
\newtheorem{Corollary}{Corollary}
\newtheorem{Notation}{Notation}

\numberwithin{equation}{section}

\begin{abstract}
 The hydrodynamic evolution of self-gravitating gaseous stars is governed by the Euler-Poisson equations. We study the structure of the linear approximation of barotropic perturbations around spherically symmetric equilibria based on functional analytic tools. 
In contrast to folklore, we show that the spectrum of the linearized operator for general perturbations is not of the Sturm-Liouville type unless the perturbations are restricted. Among others, we prove that it is of the Sturm-Liouville type for irrotational perturbations. \\

Key Words and Phrases. Gaseous star. Linear Oscillation. Self-adjoint operator. Friedrichs extension. Spectrum of Sturm-Liouville type. Spherical Harmonics. Irrotational fluid flow.

2010 Mathematical Subject Classification Numbers. 35P05, 35L51, 35Q31, 35Q85, 46N20, 76N15.
\end{abstract}


\section{Introduction}

We consider mathematical models of gaseous stars governed by the system of Euler-Poisson equations in the Cartesian  co-ordinates $(t,\vec{x})=(t, x^1, x^2, x^3)$:
\begin{subequations}
\begin{align}
&\frac{\partial\rho}{\partial t}+\sum_{k=1}^3
\frac{\partial}{\partial x^k}(\rho v^k)=0, \\
&\rho\Big(\frac{\partial v^j}{\partial t}+
\sum_{k=1}^3v^k\frac{\partial v^j}{\partial x^k}\Big)
+\frac{\partial P}{\partial x^j}+\rho\frac{\partial\Phi}{\partial x^j}=0,\quad (j=1,2,3), \\
&\triangle \Phi=4\pi\mathsf{G}\rho. \label{Poisson}
\end{align}
\end{subequations}

Here $\rho$ is the density, $P$ the pressure, and  $\vec{v}=(v^1, v^2, v^3)^T$ is the velocity field, while 
 $\mathsf{G}$ is a positive constant and we denote
$$\triangle \Phi=\sum_{k=1}^3\Big(\frac{\partial}{\partial x^k}\Big)^2\Phi.
$$
Supposing that the support of $\rho(t,\cdot)$ is compact, we replace the Poisson equation
\eqref{Poisson} by the Newtonian potential
\begin{equation}
\Phi=-4\pi\mathsf{G}\mathcal{K}\rho(t,\cdot),
\end{equation}
where the integral operator $\mathcal{K}$ is defined as
\begin{equation}
\mathcal{K}f(\vec{x})=\frac{1}{4\pi}\int\frac{f(\vec{x}')}{\|\vec{x}-\vec{x}'\|}d\vec{x}'. \label{Def_K}
\end{equation}

We put the following assumption of a barotropic pressure law :\\

{\bf (A)}: {\it $P$ is a given function of $\rho >0$ such that $P>0, dP/d\rho >0$ for
$\rho >0$ and there is a smooth function $\Lambda$ on $\mathbb{R}$ such that
$\Lambda(0)=0$ and
\begin{equation}
P=\mathsf{A}\rho^{\gamma}(1+\Lambda(\rho^{\gamma-1}))
\end{equation}
for $\rho >0$, where $\mathsf{A}, \gamma $ are positive constants such that $1<\gamma<2$.}\\

Under this assumption we denote by $u$ the enthalpy variable : 
\begin{equation}
u=\int_0^{\rho}\frac{dP}{\rho},
\end{equation}
so that  $ \displaystyle \frac{du}{d\rho} = \frac{1}{\rho}\frac{dP}{d\rho}$.  

Note that 
\begin{equation}
u=\frac{\mathsf{A}\gamma}{\gamma-1}\rho^{\gamma-1}(1+\Lambda_u(\rho^{\gamma-1})),
\end{equation}
where
$$1+\Lambda_u(X)=\frac{1}{X}\int_0^X\Big[1+\Lambda(X')+
\frac{\gamma-1}{\gamma}X'D\Lambda(X')\Big]dX'.
$$

We shall define the function $\mathsf{f}$ on $\mathbb{R}$ by
\begin{equation}
 \mathsf{f}(u)=
\begin{cases}
\rho \qquad\mbox{if}\qquad u>0 \\
0 \qquad\mbox{if} \qquad u\leq 0
\end{cases}. 
\end{equation}
Then there is a smooth function $\Lambda_{\rho}$ on $\mathbb{R}$ such that
$\Lambda_{\rho}(0)=0$ and 
\begin{equation}
\mathsf{f}(u)=\Big(\frac{\gamma-1}{\mathsf{A}\gamma}\Big)^{\frac{1}{\gamma-1}}
\max(u, 0)^{\frac{1}{\gamma-1}}(1+\Lambda_{\rho}(u)). 
\end{equation}\\

Let us write
\begin{equation}
r:=\|\vec{x}\|=\sqrt{(x^1)^2+(x^2)^2+(x^3)^2}.
\end{equation}

Moreover we put the following assumption:\\

{\bf (S)}: {\it There is a spherically symmetric  equilibrium $\rho=\bar{\rho}(\vec{x}), \vec{v}=\vec{\bar{v}}=
(0,0,0)^T$,  
which enjoys the following properties: $\bar{\rho}(\vec{x})=
\mathsf{f}(\bar{u}^{\mathsf{r}}(\|\vec{x}\|))$ with
a function $\bar{u}^{\mathsf{r}} \in C^{2,\alpha}([0,+\infty[)\cap C^{\infty}([0,R[)$, where $ 0<\alpha < 1$ and $R$ is a finite positive number, such that

S.1) it holds that  there is a finite positive constant $R$ such that
$$ \bar{u}^{\mathsf{r}} (r) >0 \qquad \Leftrightarrow \qquad 0\leq r <R; $$ 

and

S.2) it holds that }

$$ \frac{d\bar{u}^{\mathsf{r}}}{dr} < 0 \qquad \mbox{for} \quad 0<\forall r <+\infty. $$\\

The assumption {\bf (S)} is realized either if   $P=\mathsf{A}\rho^{\gamma},
\displaystyle \frac{6}{5}<\gamma < 2$ (see \cite{Chandra}, \cite{JosephL}), or if  
$\gamma \geq 4/3$ with arbitrarily given central density $\rho_{\mathsf{O}}=\rho(O)>0$ (see \cite{TM1984}). \\

Under this assumption, we fix such a stationary solution $\bar{\rho}$
with $\bar{u}: \vec{x} \mapsto \bar{u}^{\mathsf{r}}(\|\vec{x}\|)$. 

We shall use the following 

\begin{Notation}
\begin{align}
&B_R:=\{ \vec{x}\in\mathbb{R}^3\  |\  r=\|\vec{x}\| <R \}=\{ \bar{\rho} >0 \}, \\
&\overline{B_R}:=\{ \vec{x}\in\mathbb{R}^3\  |\  r=\|\vec{x}\| \leq R \}.
\end{align}
\end{Notation}

\

Despite extensive studies of spherically symmetric equilibria, typically, the Lane-Emden models, mathematically rigorous theory on the structure of evolutions of perturbations around them has not yet been fully established. The first step is to understand the linear theory. 

We consider small perturbation from this fixed equilibrium by the Lagrangian co-ordinate system, which will be denoted by the same letters $(t, x^1,x^2,x^3)$ of the Eulerian co-ordinate system. The perturbations $\xi^j$ of $x^j$ is defined by
\begin{equation}
x^j+\xi^j(t,\vec{x})=\varphi^j(t, \vec{x}),
\end{equation}
where $t \mapsto \varphi^j(t, \vec{x})$ is the co-ordinate of the  stream line, that is,
the solution of the initial value problem of the ordinary differential equation
$$\frac{\partial}{\partial t}\varphi^j(t, \vec{x})=v^j(t, \mbox{\boldmath$\varphi$}(t,\vec{x})),
\qquad \varphi^j(0, \vec{x})=x^j. $$ It is known that the linearized approximation of the equations for the perturbations is
\begin{subequations}
\begin{align}
\frac{\partial^2\xi^1}{\partial t^2}&+
\frac{\partial}{\partial x^1}\Big(-\overline{\frac{1}{\rho}\frac{dP}{d\rho}}g +4\pi\mathsf{G}\mathcal{K}g\Big)
= 0, \label{KN.5.81a}\\
\frac{\partial^2\xi^2}{\partial t^2}&+
\frac{\partial}{\partial x^2}\Big(-\overline{\frac{1}{\rho}\frac{dP}{d\rho}}g +4\pi\mathsf{G}\mathcal{K}g\Big)
= 0, \label{KN.5.81b}\\
\frac{\partial^2\xi^3}{\partial t^2}&+
\frac{\partial}{\partial x^3}\Big(-\overline{\frac{1}{\rho}\frac{dP}{d\rho}}g +4\pi\mathsf{G}\mathcal{K}g\Big)
= 0, \label{KN.5.81c}
\end{align}
\end{subequations}
where
\begin{equation}
g:= \sum_k\frac{\partial}{\partial x^k}(\bar{\rho}\xi^k). \label{KN.5.82}
\end{equation}
Here we assume that the initial perturbation of the density
$\rho(0, \vec{x})-\bar{\rho}(\vec{x})$ is supposed to vanish. 
As for the derivation of the linearized equations, see \cite[Sect. 56]{LedouxW}, \cite[pp. 139-140.]{Batchelor},  \cite{LyndenBO}, \cite[p. 500, (1)]{Lebovitz}.
Note that 
$\xi^j(0,\vec{x})=0$ and 
$\displaystyle \frac{\partial\xi^j}{\partial t}(0,\vec{x})$ is the initial perturbed velocity $v^j(0,\vec{x})$, which is supposed to be sufficiently small. \\

The aim of our study is to perform a functional
analysis of the associated vector valued integro-differential operator $\mathbf{L}$ whose components appear in the second terms of the left-hand sides of the equations \eqref{KN.5.81a} \eqref{KN.5.81b} \eqref{KN.5.81c}.

In the astrophysical literature, it often appears that the variational principle has been taken for granted, that is, the infimum of the associated quadratic form, say,  $(\mathbf{L}\mbox{\boldmath$\xi$}|\mbox{\boldmath$\xi$})$
under the constraint $\|\mbox{\boldmath$\xi$}\|=1$ in a suitable Hilbert space, is attained at a minimizer, and they would give an eigenvalue and an eigenfunction of the considered integro-differential operator $\mathbf{L}$. See, e.g., \cite{Chandra64}, \cite[\S 5.3. p.438]{Gough}, \cite[\S 3.4.1. p.76]{Tassoul}. But we shall show that this is doubtful, and the spectrum of the integro-differential operator $\mathbf{L}$, which  actually can be considered as a self-adjoint operator bounded from below in an appropriate functional Hilbert space, turns out not to be of the Sturm-Liouville type, which we mean by  is the spectrum consists of eigenvalues with finite multiplicities accumulating to $+\infty$ (see Definition \ref{Def3}).

Actually, if we restrict ourselves to spherically symmetric perturbations, the spectrum of the operator considered, which is an ordinary differential operator, turns out to be of the Sturm-Liouville type. But this fact is far from obviousness due to the strong singularity of the equilibrium at the vacuum boundary. Mathematically rigorous proof of this fact was given independently by H. Beyer \cite{Beyer}
and S. S. Lin \cite{Lin} in 1995-1997, long after astrophysicists began the discussions
mentioning the name 'Sturm-Liouville'. However, even if we restrict ourselves to axisymmetric perturbations, the spectrum of $\mathbf{L}$ cannot be of the Sturm-Liouville type (see Remark \ref{Rem.as}). There should exist a suitable restriction of the perturbations which leads us to the spectrum of the Sturm-Liouville type between spherical symmetry and axial symmetry.  An answer to this requirement is the irrotational fields of 
perturbations.

 In 1995 H. R. Beyer and B. G. Schmidt \cite{BeyerSch} proposed the necessity of mathematically rigorous justifications of the discussions on the nature of the spectrum of radial and non-radial stellar oscillations in the astrophysical literature, saying `Strangely enough this astrophysical problem was never considered in the mathematical literature'. Our study is on the line of this issue. We note that in \cite{BeyerSch}, non-radial oscillations were considered beyond barotropic perturbations with possible appearance of continuous spectrum. In our case, the spectrum of $\mathbf{L}$, although it  is not of the Sturm-Liouville type, coincides with the point spectrum and the continuous spectrum is empty. \\

{\bf  Plan of the paper and summary of the results}\\

We discuss the plan of the article and provide a brief summary of the results of the rest of Sections. 


In Section 2, we shall show that the vector valued integro-differential operator $\mathbf{L}$ acting on the vector valued perturbations $\mbox{\boldmath$\xi$}$ under considerations turns out to be a self-adjoint operator in a suitable Hilbert space $\mathfrak{H}$  of perturbation functions $\mbox{\boldmath$\xi$}$, using the Friedrichs extension theory. Moreover we shall show that the usual variational principle does not work. Roughly speaking, the reason for this inconvenience is that the control of the norm of the divergence of the vector-valued functions  of perturbations is not enough to deduce the usual compactness argument in the basic norm. 

In Section 3, we introduce a scalar valued integro-differential operator $\mathcal{N}$ which acts on  scalar functions $g$ derived by the vector perturbations
$\mbox{\boldmath$\xi$}$ as $g=\mathrm{div}(\bar{\rho}\mbox{\boldmath$\xi$})$ associated with the original vector operator $\mathbf{L}$. For this auxiliary operator $\mathcal{N}$, we shall show that it can be considered as a self-adjoint operator with spectrum of the Sturm-Liouville type in a suitable Hilbert space $\mathfrak{G}$ of scalar functions.  In order to do it, we rely on the theory of  
 weighted Sobolev spaces which permits us to treat the strong singularity of the equilibrium at the vacuum boundary.  

In Section 4, we apply the result about the scalar operator $\mathcal{N}$ to the study of the structure of the spectrum of the original vector operator $\mathbf{L}$. As a result, the structure of the spectrum of $\mathbf{L}$ is clarified immediately. We shall show that the spectrum $\sigma(\mathbf{L})$ coincides with $\sigma(\mathcal N) \cup \{0\}$ and that the kernel of $\mathbf{L}$ is infinite-dimensional and the multiplicity of nonzero eigenvalues is finite. 

In the remaining Sections, we shall discuss about more concrete information of the non-zero eigenvalues and associated eigenfunctions of $\mathbf{L}$ by using the spherical harmonics.  

In Section 5 we shall derive the information on multiplicities of the non-zero eigenvalues whose eigenfunctions are given by a specific form coming from the spherical harmonics ($Y_{lm}$) expansion of the Helmholtz decomposition of $\mbox{\boldmath$\xi$}$. It turns out that for each $l\in \mathbb N$,  we can find eigenvalues with
multiplicity 
 at least $2l+1$ for $\mathcal N$ and $\mathbf L$.  
 
In Section 6 we shall discuss on the regularity of $g=\mathrm{div}(\bar{\rho}\mbox{\boldmath$\xi$})$ and its expansion with respect to spherical harmonics, which will be foundational for the results in subsequent Sections. We shall show that $g$ is continuous in the interior of the domain and its spherical harmonics expansion is justified. A key is to establish the unique solvability and regularity of a degenerate elliptic problem $\text{div} ( \bar\rho\, \text{grad}\, U ) = f$ in suitable weighted Sobolev spaces. Such an elliptic problem is responsible for irrotational part of the perturbation.

In Section 7 we shall derive a comparison between eigenvalues of $\mathbf{L}$ and the least eigenvalue of the spherically symmetric perturbation problem. This will be done by writing an eigenvalue of $\mathbf{L}$ into the series expansion using the spherical harmonics and showing the non-negativity of each term in the series for $l \geq 1$. In the case of $P=A\rho^\gamma$, $\gamma\geq 4/3$, it turns out that no negative eigenvalues exist for $\mathbf{L}$, and the bottom of the essential spectrum is 0. 

In Section 8, we shall formulate a concept of irrotational (curl-free) vector fields, and by using the theory of weighted Sobolev spaces and the elliptic regularity result, we shall prove that when the operator $\mathbf{L}$ is restricted to these irrotational perturbations, its spectrum turns out to be of the Sturm-Liouville type. Moreover, we shall argue that the dimension of the kernel for irrotational perturbations is is the same as the one for spherically symmetric perturbations.

\section{Functional analysis of the integro-differential operator $\mathbf{L}$}

Let us consider the linearized system of equations for the perturbation 
$$\mbox{\boldmath$\xi$}=\sum_k\xi^k\frac{\partial}{\partial x^k}.$$\\

The  equations \eqref{KN.5.81a}, \eqref{KN.5.81b}, \eqref{KN.5.81c} read
\begin{equation}
\frac{\partial^2\mbox{\boldmath$\xi$}}{\partial t^2}+
\mathrm{grad}G=0, \label{1}
\end{equation}
where
\begin{align}
&G=-\frac{1}{\rho}\frac{dP}{d\rho}g+4\pi\mathsf{G}\mathcal{K}g, \\
&g=\mathrm{div}(\rho\mbox{\boldmath$\xi$}).
\end{align}

Here and hereafter we write
$ \displaystyle \rho, u, \frac{1}{\rho}\frac{dP}{d\rho}, $ etc
for 
$\displaystyle \bar{\rho}, \bar{u}, \overline{\frac{1}{\rho}\frac{dP}{d\rho}},  $ etc
for the sake of abbreviation. \\

Here recall
$$\mathcal{K}g(\vec{x})=\frac{1}{4\pi}\int\frac{g(\vec{x}')}{\|\vec{x}-\vec{x}'\|}d\vec{x}'.
\eqno\eqref{Def_K}
$$
Therefore, if the support of $g$ is included in $\overline{B_R}$ and $g \in C^{\alpha}(\overline{B_R})$ for $0<\exists \alpha <1$, then $\mathcal{K}g \in C^{2,\alpha}(\overline{B_R})$ and
$$-\triangle\mathcal{K}g=g \quad \mbox{in}\quad B_R.
$$\\

We are going to perform a functional analysis of the integro-differential operator in the equation \eqref{1}.

\subsection{Self-adjoint operator by the Friedrichs extension}

We consider the integro-differential operator $\mathbf{L}$ acting on the field $\mbox{\boldmath$\xi$}$:
\begin{equation}
\mathbf{L}\mbox{\boldmath$\xi$}=\mathrm{grad} G
\end{equation}
in the Hilbert space $\mathfrak{H}$ defined by the following

\begin{Definition}
$\mathfrak{H}$ is the Hilbert space of vector valued functions on $B_R$ endowed with the inner product 
\begin{equation}
(\mathbf{f}_1|\mathbf{f}_2)_{\mathfrak{H}}=
\int_{B_R}(\mathbf{f}_1(\vec{x})|\mathbf{f}_2(\vec{x}))\rho(\vec{x})d\vec{x}.
\end{equation}
\end{Definition}

Let us start with  the operator $\mathfrak{A}$ defined on $C_0^{\infty}(B_R; \mathbb{C}^3)$ which maps
$\mbox{\boldmath$\xi$}$ to $\mathbf{L}\mbox{\boldmath$\xi$}$.

\begin{Proposition}\label{Prop1}
$\mathfrak{A}$ is symmetric, that is,
\begin{equation}
(\mathbf{L}\mbox{\boldmath$\xi$}_1|\mbox{\boldmath$\xi$}_2)_{\mathfrak{H}}=
(\mbox{\boldmath$\xi$}_1|\mathbf{L}\mbox{\boldmath$\xi$}_2)_{\mathfrak{H}}
\end{equation}
for any $\mbox{\boldmath$\xi$}_1,\mbox{\boldmath$\xi$}_2 \in C_0^{\infty}(B_R)$.
\end{Proposition}


Proof. By integration by parts, we see
\begin{equation}
({\mathbf{L}}\mbox{\boldmath$\xi$}_1|\mbox{\boldmath$\xi$}_2)_{\mathfrak{H}}=
\int_{B_R}
\Big(\frac{1}{\rho}\frac{dP}{d\rho}g_1-
4\pi\mathsf{G}\mathcal{K}g_1\Big)g_2^*
d\vec{x}
\end{equation}

But
we see 
$$\int (\mathcal{K}g_1)(\vec{x})g_2(\vec{x})^* d\vec{x}= 
\frac{1}{4\pi}\int\int
\frac{1}{\|x-x'\|}g_1(\vec{x}')g_2(\vec{x})^* d\vec{x}'d\vec{x}
$$
is symmetric.
$\square$\\

Here and hereafter we use the following\\

\begin{Notation}
 $z^*$ stands for the complex conjugate $x-\sqrt{-1}y$ 
to the complex number $z=x+\sqrt{-1}y, x, y \in \mathbb{R}$.
\end{Notation}


\begin{Proposition}\label{Prop2}
The operator $\mathfrak{A}$ 
    is bounded from below,
that is, there exists a constant $C$ such that
\begin{equation}
(\mathbf{L}\mbox{\boldmath$\xi$}|\mbox{\boldmath$\xi$})_{\mathfrak{H}}
\geq -C\|\mbox{\boldmath$\xi$}\|_{\mathfrak{H}}^2
\quad \mbox{for}\quad \forall \mbox{\boldmath$\xi$}\in C_0^{\infty}(B_R).
\end{equation}
\end{Proposition}

Proof. We are considering
$$(\mathbf{L}\mbox{\boldmath$\xi$}|\mbox{\boldmath$\xi$})_{\mathfrak{H}}=
\int\frac{1}{\rho}\frac{dP}{d\rho}|g|^2d\vec{x}
-4\pi\mathsf{G}\int
(\mathcal{K}g)g^*{d}\vec{x}.$$

We are going to show the estimate
\begin{equation}
\clubsuit : =-\int (\mathcal{K}g)g^* d\vec{x} \geq -\rho_{\mathsf{O}}
\|\mbox{\boldmath$\xi$}\|_{\mathfrak{H}}^2. \label{4.11}
\end{equation}\\
Here $\rho_{\mathsf{O}}=\bar{\rho}(O)$.\\

We put
\begin{equation}
\Psi:=-\mathcal{K}g.
\end{equation}
Recall that $\triangle\Psi=g$ says 
\begin{equation}
\mathrm{div}(\mathrm{grad}\Psi)=\mathrm{div}(\rho\mbox{\boldmath$\xi$}).
\end{equation}
This means that, if we define $\mathbf{C}$ by
\begin{equation}
\mathbf{C}:=\rho\mbox{\boldmath$\xi$}-\mathrm{grad}\Psi,
\end{equation}
then we have 
\begin{equation}
\mathrm{div}\mathbf{C}=0.
\end{equation}\\

By integration by parts,  we have 
$$
\clubsuit=\int \Psi\cdot \mathrm{div}(\rho\mbox{\boldmath$\xi$})^* d\vec{x} =-\int(\mathrm{grad}\Psi|\rho\mbox{\boldmath$\xi$})d\vec{x}.
$$

Recall that 
$$ \rho\mbox{\boldmath$\xi$}=\mathbf{C}+\mathrm{grad}\Psi.$$
Since the support of $g=\mathrm{div}(\rho\mbox{\boldmath$\xi$})$ is a compact subset of $B_R$, we have 
$$\|\mathrm{grad}\Psi \|=O\Big(\frac{1}{r^2}\Big)\quad\mbox{as}\quad
r\rightarrow +\infty$$
so that
$$\mathbf{C}=O\Big(\frac{1}{r^2}\Big)\quad\mbox{as}\quad
r\rightarrow +\infty,
$$
too. Hence we can write
$$\clubsuit=-\int\|\mathrm{grad}\Psi\|^2d\vec{x}
-\int
(\mathrm{grad}\Psi|\mathbf{C})d\vec{x}.
$$\\

Since $\displaystyle \Psi=O\Big(\frac{1}{r}\Big), \mathbf{C}=O\Big(\frac{1}{r^2}\Big)$ as $r \rightarrow +\infty$, we get
\begin{equation}
\int
(\mathrm{grad}\Psi|\mathbf{C})d\vec{x}
=-\int\Psi(\mathrm{div}\mathbf{C})^* {d}\vec{x}=0.\label{zero}
\end{equation}

Summing up, we have 
\begin{equation}
\clubsuit=-\int
\|\mathrm{grad}\Psi\|^2d\vec{x}.
\end{equation}
But, using \eqref{zero}, we see
\[
\begin{split}
\int\|\mathrm{grad}\Psi\|^2d\vec{x} &\leq 
\int\|\mathrm{grad}\Psi\|^2d\vec{x} + \int
\|\mathbf{C} \|^2d\vec{x}\\
&= \int\|\mathrm{grad}\Psi+\mathbf{C}\|^2d\vec{x} \\
&=  \int\|\rho \mbox{\boldmath$\xi$}\|^2d\vec{x} \\
&=\int\|\rho \mbox{\boldmath$\xi$}\|^2d\vec{x} \leq \rho_{\mathsf{O}}\|\mbox{\boldmath$\xi$}\|_{\mathfrak{H}}^2
\end{split}
\]

This completes the proof. $\square$\\

Proposition \ref{Prop1},
Proposition \ref{Prop2} imply the following conclusion:\\

The operator $\mathfrak{A}$, defined by
$\mathsf{D}(\mathfrak{A})=C_0^{\infty}(B_R),
\mathfrak{A}\mbox{\boldmath$\xi$} =\mathbf{L}\mbox{\boldmath$\xi$}$, admits the Friedrichs extension 
$\mathfrak{T}$ which is a self-adjoint operator in the Hilbert space $\mathfrak{H}$. 
See e.g. \cite[Chapter VI, \S 2.3]{Kato}. The domain $\mathsf{D}(\mathfrak{T})$ of the operator $\mathfrak{T}$ is 
\begin{align}
\mathsf{D}(\mathfrak{T})=&\mathfrak{H}_1\cap \mathsf{D}(\mathfrak{A}^*) \nonumber \\
=&
\{ \mbox{\boldmath$\xi$}\in\mathfrak{H}_1\  |\  \mathbf{L}\mbox{\boldmath$\xi$}\in\mathfrak{H} \quad\mbox{in distribution sense} \}.
\end{align}

Here $\mathfrak{H}_1$ is defined by the following

\begin{Definition}\label{Def.H1}

1) Let $L^2(du/d\rho)$ denote the Hilbert space $\displaystyle L^2\Big(B_R, \frac{1}{\rho}\frac{dP}{d\rho}d\vec{x}\Big)=L^2\Big(B_R, \frac{du}{d\rho}d\vec{x}\Big)$ endowed with the norm
$\|\cdot\|_{L^2(du/d\rho)}$ defined by
$$\|g\|_{L^2(du/d\rho)}^2:=
\int_{B_R}|g(\vec{x})|^2\frac{1}{\rho}\frac{dP}{d\rho}(\vec{x})d\vec{x}. $$

2)
$\mathfrak{H}_1$ is the set of all $\mbox{\boldmath$\xi$}\in\mathfrak{H}$
with $g=\mathrm{div}(\rho\mathrm{grad}\mbox{\boldmath$\xi$}) \in L^2(du/d\rho)$ 
in distribution sense such that there exists a sequence $\mbox{\boldmath$\varphi$}_n \in C_0^{\infty}(B_R)$ such that
$$\mbox{\boldmath$\varphi$}_n\rightarrow \mbox{\boldmath$\xi$} \quad\mbox{in}\quad \mathfrak{H}\quad
\mbox{as}\quad n\rightarrow \infty,$$
and
$\mathrm{div}(\rho\mathrm{grad}\mbox{\boldmath$\varphi$}_n) \rightarrow g$ in $L^2(du/d\rho)$. Here 
`$g=\mathrm{div}(\rho\mbox{\boldmath$\xi$})$ in distribution sense' means
$$\int_{B_R} g\varphi^* =\int_{B_R}(\mbox{\boldmath$\xi$}|\mathrm{grad}\varphi)\rho
$$
for $\forall \varphi \in C_0^{\infty}(B_R)$.
\end{Definition}

Here we have used the following
\begin{Proposition}
If $\mbox{\boldmath$\varphi$}_n \in C_0^{\infty}(B_R)$ satisfies
$\mbox{\boldmath$\varphi$}_n \rightarrow \mbox{\boldmath$\xi$}$ in $\mathfrak{H}$ and
$\mathrm{div}(\rho\mbox{\boldmath$\varphi$}_n) \rightarrow g=\mathrm{div}(\rho\mbox{\boldmath$\xi$})$ in 
$L^2(du/d\rho)$, then
$Q[\mbox{\boldmath$\varphi$}_m-\mbox{\boldmath$\varphi$}_n]
\rightarrow 0$, where
$$Q[\mbox{\boldmath$\varphi$}]:=
\int|\mathrm{div}(\rho\mbox{\boldmath$\varphi$})|^2\frac{1}{\rho}\frac{dP}{d\rho}
-4\pi\mathsf{G}\int
(\mathcal{K}\mathrm{div}(\rho\mbox{\boldmath$\varphi$}))(\mathrm{div}(\rho\mbox{\boldmath$\varphi$}))^*.
$$
\end{Proposition}

Proof. Write
$$Q[\mbox{\boldmath$\varphi$}]=\check{Q}_{00}[\mathrm{div}(\rho\mbox{\boldmath$\varphi$})]+
\check{Q}_{01}[\mathrm{div}(\rho\mbox{\boldmath$\varphi$})],$$
where
\begin{align*}
&\check{Q}_{00}[g]=\int |g|^2\frac{1}{\rho}\frac{dP}{d\rho}, \\
&\check{Q}_{01}[g]=-4\pi\mathsf{G}\int
(\mathcal{K}g)g^*.
\end{align*}
If $g_n=\mathrm{div}(\rho\mbox{\boldmath$\varphi$}_n)$ converges in $L^2(du/d\rho)$, then clearly
$\check{Q}_{00}[g_m-g_n]\rightarrow 0$. Thus, it is sufficient to show that
$\check{Q}_{01}[g_m-g_n]\rightarrow 0$. We put
$$\Psi=-\mathcal{K}g,
\qquad
\Psi_n=-\mathcal{K}g_n=-\mathcal{K}\mathrm{div}(\rho\mbox{\boldmath$\varphi$}_n).
$$
Then we have
\begin{align*}
|\Psi(\vec{x})|&=\Big|
\frac{1}{4\pi}\int\frac{g(\vec{x}')}{\|\vec{x}-\vec{x}'\|}d\vec{x}'\Big| \\
&\leq \frac{1}{4\pi}\sqrt{ \int_{\|\vec{x}'\|\leq R}\frac{d\vec{x}'}{\|\vec{x}-\vec{x}'\|^2} }
\|g\|_{L^2(B_R)} \\
&\leq \sqrt{\|\vec{x}\|+R}\|g\|_{L^2(B_R)},
\end{align*}
therefore
$$\|\Psi\|_{L^2(B_R)} \lesssim
\|g\|_{L^2(B_R)}
\lesssim
\|g\|_{L^2(du/d\rho)}.
$$
Thus
\begin{align*}
&\Big|\int\mathcal{K}(g_m-g_n)g_m^*\Big| = \Big|\int(\Psi_m-\Psi_n)g_m^* \Big| \\
&\leq \|\Psi_m-\Psi_n\|_{L^2(B_R)}\|g_m\|_{L^2(B_R)} \\
&\lesssim \|g_m-g_n\|_{L^2(du/d\rho)}\|g_m\|_{L^2(du/d\rho)} \rightarrow 0
\end{align*}
and
\begin{align*}
&\Big|\int\mathcal{K}g_n(g_m-g_n)^*\Big|
\leq \|\Psi_n\|_{L^2(B_R)}\|g_m-g_n\|_{L^2(B_R)} \\
&\lesssim \|g_n\|_{L^2(du/d\rho)}\|g_m-g_n\|_{L^2(du/d\rho)}
\rightarrow 0.
\end{align*}
Therefore $\check{Q}_{01}[g_m-g_n]\rightarrow 0$.
 $\square$
\\

And `$\mathbf{L}\mbox{\boldmath$\xi$}\in \mathfrak{H}$ in distribution sense' means that there exists $\mathbf{f}\in\mathfrak{H}$ such that
$$(\mbox{\boldmath$\xi$}|\mathbf{L}\mbox{\boldmath$\varphi$})_{\mathfrak{H}}
=(\mathbf{f}|\mbox{\boldmath$\varphi$})_{\mathfrak{H}}$$
for any $\mbox{\boldmath$\varphi$}\in C_0^{\infty}(B_R)$.\\

Hereafter we shall denote by the same letter $\mathbf{L}$ the Friedrichs extension $\mathfrak{T}$. Thus we can claim the following

\begin{Theorem}\label{Th.1}
The operator $\mathbf{L}$ is a self-adjoint operator bounded from below in the Hilbert space $\mathfrak{H}$.
\end{Theorem}

\begin{Remark}\label{R3}
A characterization of the domain
$\mathsf{D}(\mathbf{L})$ is not transparent. But 
$\displaystyle \mbox{\boldmath$\xi$}=ry(r)\frac{\partial}{\partial r} $ belongs to
$ \mathsf{D}(\mathbf{L})$ if $y \in \mathsf{D}( \mathcal{L}^{\mathsf{ss}})$, where 
$\mathcal{L}^{\mathsf{ss}}$ is the Friedrichs extension of the ordinary differential operator
$$\mathcal{L}^{\mathsf{ss}}y=-\frac{1}{\rho r^4}
\frac{d}{dr}\Big(\gamma r^4P\frac{dy}{dr}\Big)
-(3\gamma-4)\frac{1}{r}\frac{du}{dr}y$$
in $L^2([0,R]; \rho r^4dr)$,   which will be denoted by $\mathfrak{W}$ later by Definition \ref{Def.X}. 
Here we have supposed $P=\mathsf{A}\rho^{\gamma}, \frac{6}{5} < \gamma <2$ for the simplicity.
 Clearly
$$ \mathbf{L}\mbox{\boldmath$\xi$}=r\Big(\mathcal{L}^{\mathsf{ss}}y\Big)\frac{\partial}{\partial r} 
\qquad\mbox{for}\qquad \displaystyle \mbox{\boldmath$\xi$}=ry\frac{\partial}{\partial r}.$$ 
A characterization of the domain $\mathsf{D}(\mathcal{L}^{\mathsf{ss}})$ is known. See 
\cite[p. 554]{OJM}. In particular we have  $C^2([0,R]) \subset \mathsf{D}(\mathcal{L}^{\mathsf{ss}})$. Note that it is not necessary that $y(r) \in \mathsf{D}(\mathcal{L}^{\mathsf{ss}})$ vanishes at $r=R-0$,
while $r^2(\gamma P\rho)^{1/4}y$ vanishes at $r=R-0$ if $y=O(1)$.
\end{Remark}

\subsection{Does the variational principle work?}

Since $Q[\mbox{\boldmath$\xi$}]$ is bounded from below, we have the finite infimum
\begin{equation}
\underline{\lambda} :=\inf \{ Q[\mbox{\boldmath$\xi$}]\   |\   \mbox{\boldmath$\xi$}\in\mathfrak{H}_1, \|\mbox{\boldmath$\xi$}\|_{\mathfrak{H}}=1\}.
\end{equation}

At the moment, we have no information on the signature of $\underline{\lambda}$. So,
let us put 
\begin{equation}
\mathcal{Q}(\mbox{\boldmath$\xi$}_1,\mbox{\boldmath$\xi$}_2):=Q(\mbox{\boldmath$\xi$}_1,\mbox{\boldmath$\xi$}_2)+(1-\underline{\lambda})(\mbox{\boldmath$\xi$}_1|\mbox{\boldmath$\xi$}_2)_{\mathfrak{H}}
\end{equation}
so that 
$$\inf \{ \mathcal{Q}[\mbox{\boldmath$\xi$}]\  |\  \mbox{\boldmath$\xi$}\in \mathfrak{H}_1, 
\|\mbox{\boldmath$\xi$}\|_{\mathfrak{H}}=1\} =1.
$$

We consider that $\mathfrak{H}_1$ is a Hilbert space endowed with this inner product
$\mathcal{Q}(\cdot,\cdot)$ and the norm
$\sqrt{\mathcal{Q}[\mbox{\boldmath$\xi$}]}=\sqrt{\mathcal{Q}(\mbox{\boldmath$\xi$},\mbox{\boldmath$\xi$})}$.\\

We might want to show that $\underline{\lambda}$ is the least eigenvalue of the operator 
$\mathbf{L}$. The usual argument to show it is as following: \\

Let us consider a minimizing sequence $(\mbox{\boldmath$\xi$}_n)_n$ of 
$\mathcal{Q}$, that is, $\mbox{\boldmath$\xi$}_n\in\mathfrak{H}_1, 
\|\mbox{\boldmath$\xi$}_n\|_{\mathfrak{H}}=1$ and
$\mathcal{Q}[\mbox{\boldmath$\xi$}_n]\rightarrow 1$
or $Q[\mbox{\boldmath$\xi$}_n]\rightarrow \underline{\lambda}$. So, we want to prove that
there exists a subsequence which converges to a limit $\mbox{\boldmath$\xi$}_{\infty}\in \mathsf{D}(\mathbf{L})$ and satisfies
$\mathbf{L}\mbox{\boldmath$\xi$}_{\infty}=\underline{\lambda}\mbox{\boldmath$\xi$}_{\infty}$, that is,
$\underline{\lambda}$ is an eigenvalue of $\mathbf{L}$ and $\mbox{\boldmath$\xi$}_{\infty}$ is an associated eigenfunction.
In order to do this, the usual argument adopts the theory developed in
\cite[Kapitel VII]{CH} assuming the following `Rellichscher Auswahlsatz' (\cite[p. 489, \S VII.3. Satz 2]{CH}):\\

($\heartsuit $): {\it If $\mbox{\boldmath$\Xi$}_n\in\mathfrak{H}_1$ satisfies
$\mathcal{Q}[\mbox{\boldmath$\Xi$}_n]\leq \exists M$, then there is a subsequence
$(\mbox{\boldmath$\Xi$}_{n_k})_k$ which converges in $\mathfrak{H}$. In other words, the imbedding
$\mathfrak{H}_1 \hookrightarrow \mathfrak{H}$ is compact.}\\

However the above argument doesn't work. In order to show this, let us introduce the following definition:

\begin{Definition}\label{Def3}
The spectrum $\sigma(T)$ of a self-adjoint operator $T$ in an infinitely dimensional Hilbert space
$\mathsf{X}$ is said to be of the Sturm-Liouville type if $\sigma(T)$ consists of
eigenvalues with finite multiplicities which accumulate to $+\infty$.
\end{Definition}
\begin{Remark}
To be of the Sturm-Liouville type for the spectrum in sense of the above definition is said to be `discrete' in many literatures on the spectral theory of differential operators. (E.g., \cite[p. 187]{Kato}, \cite[p. 73]{Davies}, \cite[p.132]{Helffer} and so on.) But we avoid this terminology, since the spectrum can be a discrete subset of the complex plane $\mathbb{C}$ but does contain a point which is not an eigenvalue of finite multiplicity. Actually, later (see Theorem \ref{Th.3}), we shall show that it is the case for $\mathbf{L}$, say, the spectrum $\sigma(\mathbf{L})$ consists of eigenvalues of finite multiplicity $\lambda_j$,
$\lambda_j\not=0$, and $\{0\}$ for which the dimension of $\mathrm{Ker}(0-\mathbf{L})$ is infinity, say, $\{0\}$ is an essential spectrum, although it is  isolated, when we consider $\mathbf{L}$ in the Hilbert space
$\mathfrak{F}=\{ \mbox{\boldmath$\xi$} \in \mathfrak{H}\  |\  \mathrm{div}(\rho\mbox{\boldmath$\xi$}) \in \mathfrak{G} \}$.
\end{Remark}

Then the Riesz-Schauder's theorem (see \cite[Chapter III, Theorem 6.29]{Kato}) reads \\

\noindent {\bf The Riesz-Schauder's theorem: } {\it If a resolvent $(\lambda-T)^{-1}$ of the self-adjoint operator $T$ bounded from below is a compact operator, then the spectrum of the operator $T$ is of the Sturm-Liouville type.}\\

\begin{Remark} Recall that if $\lambda, \mu$ belong to the resolvent set of $T$, then 
the resolvent $(\lambda-T)^{-1}$ is a compact operator if and only if $(\mu-T)^{-1}$ is so, for the resolvent equation
$$(\lambda-T)^{-1}-(\mu-T)^{-1}=
(\mu-\lambda)(\lambda-T)^{-1}(\mu-T)^{-1}$$ holds.
\end{Remark}

Now the operator $\mathbf{L}+1-\underline{\lambda}$ has the inverse 
$(\mathbf{L}+1-\underline{\lambda})^{-1}$ which is a bounded linear operator 
on $\mathfrak{H}$ into itself with operator norm $\leq 1$. It is a compact operator if $(\heartsuit$) holds.
Therefore if ($\heartsuit$) hold, then by the Riesz-Schauder's theorem would imply that the spectrum of $\mathbf{L}$ is of the Sturm-Liuouville type.\\

However we see that $(\heartsuit)$ does not hold unfortunately.

Let us consider the functional space
\begin{equation}
\mathfrak{N}=\{ \mbox{\boldmath$\xi$} \in\mathfrak{H}\  |\   \mathrm{div}(\rho \mbox{\boldmath$\xi$})=0 \quad\mbox{in distribution sense} \}.
\end{equation}
Here `$\mathrm{div}(\rho \mbox{\boldmath$\xi$})=0$ in distribution sense' means 
$$ 
\int (\mbox{\boldmath$\xi$}|\mathrm{grad}\varphi)\rho d\vec{x} =0
\quad\forall \varphi \in C_0^{\infty}(B_R).
$$

Clearly $\mathfrak{N}\subset \mbox{Ker}\mathbf{L}$, that is, any $\mbox{\boldmath$\xi$}\in\mathfrak{N}, \not=0$ is an eigenfunction of
$\mathbf{L}$ associated with the eigenvalue 0.
But the dimension of $\mathfrak{N}$ is infinite.  In fact, for any vector function $
\mathbf{A} \in C_0^{\infty}(B_R)$,
the vector function
$$\mbox{\boldmath$\xi$}=\frac{1}{\rho}\mbox{curl}\mathbf{A}
$$
belongs to $\mathfrak{N}$.  This is a contradiction to the assertion that the spectrum of $\mathbf{L}$ was of the Sturm-Liouville type. \\

Therefore we cannot expect the `Rellichscher Auswahlsatz' ($\heartsuit$), if we do not
limit  ourselves to spherically symmetric perturbations or curl-free perturbations. Then it is doubtful that the so called `variational principle'   that the minimum $\underline{\lambda}$ should be attained at an eigenfunction $\mbox{\boldmath$\xi$}_{\infty}$. \\

In fact, as proved later, the situation is as follows:

One consider the so-called `Max-Min principle' by putting
\begin{align*}
\mu_n:=&\sup_{\psi_1,\cdots,\psi_{n-1}}
\inf\{(\mathbf{L}\phi|\phi)_{\mathfrak{F}}\  | \\
& |\qquad   \phi \in [\mathrm{span}(\psi_1,\cdots,\psi_{n-1})]^{\perp}, \phi\in \mathsf{D}(\mathbf{L}), \|\phi\|_{\mathfrak{F}=1} \}. 
\end{align*}
See \cite[Chapter 11]{Helffer}. Now, suppose that $P=\mathsf{A}\rho^{\gamma}, \frac{4}{3}<\gamma <2$, for simplicity. Then, as shown later, (see Theorem \ref{Th.3}, Theorem \ref{Th.6}, Corollary \ref{Cor.6.1}) the spectrum of the self-adjoint operator $\mathbf{L}$ considered in the Hilbert space
$\mathfrak{F}$ is $\{ 0\}\cup\{\lambda_n | n=1,2,\cdots \}$, where $0<\lambda_1\leq \lambda_2 \leq \cdots$ are positive eigenvalues of finite multiplicities but $0$ is an eigenvalue of infinite multiplicity, so that $\{ 0\}$ is an essential spectrum. Therefore, according to \cite[Theorem 11.7]{Helffer}, $\mu_n$ is blocked by this essential spectrum $0$ so that $\mu_n=0$ for $\forall n$ and the smallest positive eigenvalue $\lambda_1$ cannot be reached by the variational principle. 

Of course, if $P=\mathsf{A}\rho^{\gamma}, \frac{6}{5}< \gamma <\frac{4}{3}$, the least eigenvalue $\lambda_1$ is negative, and negative eigenvalues must be reached by the variational principle,
or Max-Min principle, 
before the essential spectrum $0$. 
But this is a judgment of a matter from the results of the discussions which will be developed later, after the next Section. The `Rellichscher Auswahlsatz' $(\heartsuit)$ cannot be applied here and now. \\

Thus we seek another formulation of the linearized analysis. 

 \section{Another formulation of the linearized problem}
 
 
Let us introduce the variables
 \begin{equation}
\hat{\mbox{\boldmath$\xi$}}
=\rho\mbox{\boldmath$\xi$}.
 \end{equation}
Then we have
 \begin{equation}
 g=\mathrm{div}\hat{\mbox{\boldmath$\xi$}}
 \end{equation}
 and
 the linearized equations are reduced to
 \begin{equation}
 \frac{\partial^2\hat{\mbox{\boldmath$\xi$}}}{\partial t^2}+
\hat{\mathbf{M}}g=0,\label{5.3}
 \end{equation}
 where
 
\begin{equation}
\hat{\mathbf{M}}g=\rho\mathrm{grad}\Big(
-\frac{1}{\rho}\frac{dP}{d\rho}g+4\pi\mathsf{G}\mathcal{K}g \Big).
\end{equation}

 Taking the divergence of \eqref{5.3}, we get
 \begin{equation}
 \frac{\partial^2g}{\partial t^2}+\mathcal{N} g=0,
 \end{equation}
 where
 
\begin{align}
\mathcal{N}g
=&\mathrm{div}\hat{\mathbf{M}}g  \nonumber \\
=&\mathrm{div}\Big(\rho\mathrm{grad}
\Big(-\frac{1}{\rho}\frac{dP}{d\rho}g+4\pi\mathsf{G}\mathcal{K}g\Big)\Big) \nonumber \\
=&-\rho\frac{d\rho}{dP}
\mathrm{div}
\Big(
\frac{1}{\rho}\Big(\frac{dP}{d\rho}\Big)^2
\mathrm{grad}g\Big)-
\Big[\triangle\Big(\frac{dP}{d\rho}-u\Big)\Big]g \nonumber \\
&+4\pi\mathsf{G}\mathrm{div}\Big(\rho\mathrm{grad}
(\mathcal{K}g)\Big).
\end{align}\\

  Thus, under the initial conditions
\begin{subequations}
\begin{align}
&\hat{\mbox{\boldmath$\xi$}}|_{t=0}=\vec{0},\quad g|_{t=0}=0, \\
&\frac{\partial\hat{\mbox{\boldmath$\xi$}}}{\partial t}\Big|_{t=0}=
\overset{\circ}{\hat{\mathbf{v}}}
(\not=\vec{0}),
\quad
\frac{\partial g}{\partial t}\Big|_{t=0}=
\overset{\circ}{g}=\mathrm{div}\overset{\circ}{\hat{\mathbf{v}}},
\end{align}
\end{subequations}
we consider the system
\begin{subequations}
\begin{align}
&\frac{\partial^2\hat{\mbox{\boldmath$\xi$}}}{\partial t^2}+\hat{\mathbf{M}}g=0, \label{5.10a}\\
&\frac{\partial^2g}{\partial t^2}+
\mathcal{N} g=0,\label{5.10b}
\end{align}
\end{subequations}
under the constraint
\begin{equation}
g=\mathrm{div}\hat{\mbox{\boldmath$\xi$}},
\end{equation}
which is conserved during the evolution.

Note that \eqref{5.10a} is a system of ordinary differential equations
re $t$ with parameter $x$ of $\hat{\mbox{\boldmath$\xi$}}$ provided that
$g$ is given, and 
\eqref{5.10b} is a single second order wave equation of $g$
provided that $\hat{\mbox{\boldmath$\xi$}}$ is given. 

Thus we are going to perform a functional analysis of the operator $\mathcal{N}$ only.\\

Let us consider the operator $\mathcal{N}$ in a Hilbert space $\mathfrak{G}$ defined below.

First let us consider the Hilbert space 
$L^2(du/d\rho)= L^2\Big(B_R; \frac{1}{\rho}\frac{dP}{d\rho}d\vec{x}\Big)
=L^2\Big(B_R; \frac{du}{d\rho}d\vec{x}\Big) $ endowed the norm 
$\|\cdot\|_{L^2(du/d\rho)}$ defined by Definition \ref{Def.H1}.
If $ g \in L^2\Big(B_R; \frac{1}{\rho}\frac{dP}{d\rho}d\vec{x}\Big) $, then $g \in L^1(B_R)$ and
$$\int_{B_R}|g|d\vec{x} \leq \|g\|_{L^2(du/d\rho)}
\sqrt{\int_{B_R}\frac{d\rho}{du}d\vec{x}}. $$
Here we note $\frac{d\rho}{du} \lesssim (R-\|\vec{x}\|)^{\frac{2-\gamma}{\gamma-1}} =O(1)$. Therefore $g \mapsto \int_{B_R}g(\vec{x})d\vec{x}$ is a continuous linear functional on
$L^2\Big(B_R; \frac{1}{\rho}\frac{dP}{d\rho}d\vec{x}\Big)$ and $\int g=0$ gives a closed subspace of $L^2\Big(B_R; \frac{1}{\rho}\frac{dP}{d\rho}d\vec{x}\Big)$. So we put the following

\begin{Definition}
We put
\begin{equation}
\mathfrak{G} :=\{ g \in L^2\Big(B_R; \frac{1}{\rho}\frac{dP}{d\rho}d\vec{x}\Big)\quad | \quad
\int_{B_R}g(\vec{x})d\vec{x} =0 \},
\end{equation}
which is a Hilbert space endowed with the  norm $\|\cdot\|_{\mathfrak{G}}=
\|\cdot\|_{L^2(du/d\rho)}$.
\end{Definition}

Let us analyze the operator $\mathcal{N}$, which can be decomposed as

\begin{subequations}
\begin{align}
\mathcal{N}g&=\mathcal{N}_{00}g+\mathcal{N}_{01}g, \\
\mathcal{N}_{00}g&=
\mathrm{div}\Big(\rho \mathrm{grad}\Big(-\frac{1}{\rho}\frac{dP}{d\rho}g\Big)\Big) \nonumber \\ 
&=-\rho\frac{d\rho}{dP}\mathrm{div}
\Big(\frac{1}{\rho}\Big(\frac{dP}{d\rho}\Big)^2
\mathrm{grad}g\Big)-\Big[\triangle\Big(\frac{dP}{d\rho}-u\Big)\Big] g, \\
\mathcal{N}_{01}g&=4\pi\mathsf{G}
\mathrm{div}(\rho\mathrm{grad}(\mathcal{K}g)).
\end{align}
\end{subequations}\\

If $g_1, g_2 \in C_0^{\infty}(B_R)$, then 
$$(\mathcal{N}_{00}g_1|g_2)_{\mathfrak{G}}=Q(g_1,g_2), $$
where
\begin{align}
Q(g_1,g_2)&=
\int
\frac{1}{\rho}\Big(\frac{dP}{d\rho}\Big)^2
(\mathrm{grad}g_1|\mathrm{grad}g_2)d\vec{x} + \nonumber \\
&-\int\Big[\triangle\Big(\frac{dP}{d\rho}-u\Big)\Big]
g_1g_2^*\frac{1}{\rho}\frac{dP}{d\rho}
d\vec{x}.
\end{align}
Let $\mathcal{A}_{00}: g\mapsto \mathcal{N}_{00}g$ be the operator on the domain
$\mathsf{D}(\mathcal{A}_{00})=C_0^{\infty}(B_R)$.

Note that it is easy to see that $\mathcal{N}_{00}g, \mathcal{N}_{01}g \in \mathfrak{G}$, for which
$\int_{B_R} =0$, provided that $g \in C_0^{\infty}(B_R)$.

Then $\mathcal{A}_{00}$ is symmetric and bounded from below, since  
$$
\Big|\triangle\Big(\frac{dP}{d\rho}-u\Big)\Big|\leq C.
$$

Taking $\kappa$ sufficiently large, we have that the bilinear form
$$\mathcal{Q}(g_1,g_2)=Q(g_1,g_2)+\kappa (g_1|g_2)_{\mathfrak{G}}$$
associated with  $\mathcal{A}_{00}+\kappa$ satisfies
$$\mathcal{Q}[g]\geq \|g\|_{\mathfrak{G}}^2.$$
Thus we get the Friedrichs extension $\mathcal{T}_{00}$ of
$\mathcal{A}_{00}$ which is self-adjoint. Then
we have the bounded inverse $(\mathcal{T}_{00}+\kappa)^{-1}$
with the operator norm $\leq 1$. We want to show the inverse is compact (completely continuous), that is, the closure in
$\mathfrak{G}$ of
$$\mathcal{B}=\{ g\ |\ \mathcal{Q}[g]\leq 1\}$$
is compact.\\

Recall that
\begin{align*}
\mathcal{Q}[g]&=
\int
\frac{1}{\rho}\Big(\frac{dP}{d\rho}\Big)^2\|\mathrm{grad}g\|^2d\vec{x} +\\
&+\int
\Big(-\Big[\triangle\Big(\frac{dP}{d\rho}-u\Big)\Big] +\kappa\Big)|g|^2\frac{1}{\rho}\frac{dP}{d\rho}
d\vec{x} \\
&=
\int
\frac{1}{\rho}\Big(\frac{dP}{d\rho}\Big)^2
\|\mathrm{grad}g\|
d\vec{x} +\\
&+\int
\Big(-\Big[\triangle\Big(\frac{dP}{d\rho}-u\Big)\Big] +\kappa\Big)|g|^2\frac{1}{\rho}\frac{dP}{d\rho}
d\vec{x} 
\end{align*}
can control $\|\mathrm{grad}g\|$. 
We are assuming that
$$-\Big[\triangle\Big(\frac{dP}{d\rho}-u\Big)\Big] +\kappa \geq 1.$$\\

Precisely let us introduce the following :

\begin{Definition}
We  put 
\begin{align}
\mathfrak{G}_1&=\{ g  \in \mathfrak{G}\  |\  \exists \varphi_n \in C_0^{\infty}(B_R) \quad\mbox{such that} \quad 
\|\varphi_n-g\|_{\mathfrak{G}}\rightarrow 0\quad (n\rightarrow \infty), \nonumber \\
&\mbox{and}\quad \mathcal{Q}[\varphi_m-\varphi_n]\rightarrow 0 \quad(m,n\rightarrow\infty) \}
\end{align} 
and
we consider $\mathfrak{G}_1$ as a Hilbert space endowed with the inner product
\begin{equation}
(g_1|g_2)_{\mathfrak{G}_1}=
\int\Big[
\frac{1}{\rho}\Big(\frac{dP}{d\rho}\Big)^2
(\mathrm{grad}g_1|\mathrm{grad}g_2)
+\frac{1}{\rho}\frac{dP}{d\rho}g_1g_2^*\Big]d\vec{x}.
\end{equation}
\end{Definition}

Of course
$$\frac{1}{C}\|g\|_{\mathfrak{G}_1}^2
\leq \mathcal{Q}[g]\leq C\|g\|_{\mathfrak{G}_1}^2.
$$

We claim
\begin{Proposition}\label{Prop3}
The closure of the unit ball of $\mathfrak{G}_1$ in $\mathfrak{G}$ is compact, and therefore the resolvent $(\mathcal{T}_{00}+\kappa)^{-1}$ is a compact operator in $\mathfrak{G}$.
\end{Proposition}

Before giving a proof of this Proposition, we introduce the notations to denote various weighted spaces according to   \cite{GurkaO}:\\

\begin{Notation} 
1)  If $s, s_0, s_1$ are positive continuous functions on $B_R$, we consider the Hilbert spaces 
$L^2(B_R, s), W^{1,2}(B_R, s_0, s_1)$ endowed with the norms
$\|\cdot\|_{L^2(s)}, \|\cdot\|_{W^{1,2}(s_0,s_1)}$ defined by
\begin{equation}
\|f\|_{L^2(s)}^2=\int_{B_R}|f(\vec{x})|^2s(\vec{x})d\vec{x},
\end{equation}
and
\begin{equation}
\|f\|_{W^{1,2}(s_0,s_1)}^2=
\|f\|_{L^2(s_0)}^2+\|\nabla f\|_{L^2(s_1)}^2.
\end{equation}

2) We shall use the function $\mathsf{d}$ defined by
\begin{equation}
\mathsf{d}(\vec{x})=\mbox{dist.}(\vec{x}, \partial B_R)=R-\|\vec{x}\| .
\end{equation}

Recall that 
$$ \frac{1}{C}u(\vec{x})\leq \mathsf{d}(\vec{x})\leq C u(\vec{x}) \quad\mbox{on}\quad B_R, 
$$ 
thanks to the physical vacuum boundary condition. Therefore $L^2(B_R, \mathsf{d}^{\alpha})=L^2(B_R, u^{\alpha})$ and so on. 

3) We denote
\begin{equation}
\nu:=\frac{1}{\gamma-1}.
\end{equation}
\end{Notation}

\

Since we are supposing $1<\gamma <2$, we have $1<\nu <+\infty$. 

We observe that $\mathfrak{H}=L^2(B_R,\rho )=L^2(B_R, \mathsf{d}^{\nu})$.

On the other hand $$\mathfrak{G}=
 L^2(B_R, \frac{1}{\rho}\frac{dP}{d\rho}) \cap\{ g | \int g=0\}
=L^2(B_R, \mathsf{d}^{1-\nu}) \cap\{ g |  \int g=0\}.$$
 And $\mathfrak{G}_1$ is the closure of $C_0^{\infty}(B_R)$ in
$$\displaystyle W^{1,2}\Big(B_R, \frac{1}{\rho}\frac{dP}{d\rho}, \frac{1}{\rho}
\Big(\frac{dP}{d\rho}\Big)^2\Big)\cap\{ \int =0\}=
W^{1,2}(B_R, \mathsf{d}^{1-\nu}, \mathsf{d}^{2-\nu})\cap\{\int =0\}.$$.


Proof of Proposition \ref{Prop3}.  We have to show that the imbedding of 

\noindent  $W_0^{1,2}(B_R, \mathsf{d}^{-\nu+1},\mathsf{d}^{-\nu+2})$ into
$L^2(B_R, \mathsf{d}^{-\nu+1})$ is compact, 

\noindent since $\mathfrak{G}_1 \subset
W_0^{1,2}(B_R, \mathsf{d}^{-\nu+1},\mathsf{d}^{-\nu+2})$, which stands for the closure of
$C_0^{\infty}(B_R)$ in
$W^{1,2}(B_R, \mathsf{d}^{-\nu+1},\mathsf{d}^{-\nu+2})$.

First, by the Hardy imbedding inequality,  \cite[Theorem 8.4]{Kufner}
or \cite[Theorem 8.7]{GurkaO.2}, we can imbed 
$W_0^{1,2}(B_R, \mathsf{d}^{-\nu+2}, \mathsf{d}^{-\nu+2})$ 
 continuously into $L^2(B_R, \mathsf{d}^{-\nu})$. 

But it is clear that
$W_0^{1,2}(B_R, \mathsf{d}^{-\nu+1}, \mathsf{d}^{-\nu+2})$
is continuously imbedded into
$W_0^{1,2}(B_R, \mathsf{d}^{-\nu+2},\mathsf{d}^{-\nu+2})$.
Thus 
$W_0^{1,2}(B_R, \mathsf{d}^{-\nu+1}, \mathsf{d}^{-\nu+2})$
is continuously imbedded into
$W^{1,2}(B_R, \mathsf{d}^{-\nu}, \mathsf{d}^{-\nu+2})$.

On the other hand, by \cite[Theorem 2.4]{GurkaO} the imbedding of
$W^{1,2}(B_R, \mathsf{d}^{-\nu}, \mathsf{d}^{-\nu+2})$ into $L^2(B_R, \mathsf{d}^{-\nu+1})$ is compact.
This completes the proof.
$\square$ \\


Now let us deal with the perturbation $\mathcal{N}_{01}$.\\

Let $\mathcal{A}: g \mapsto \mathcal{N}g=\mathcal{N}_{00}g+\mathcal{N}_{01}g$ be the operator on the domain $\mathsf{D}(\mathcal{A})=C_0^{\infty}(B_R)$. \\

Here let us note that $\mathcal{N}_{01}g \in \mathfrak{G}$ if $g \in C_0^{\infty}(B_R)$.
In fact, keeping in mind that $-\triangle\mathcal{K}g=g$, we have 
$$ \mathcal{N}_{01}g=4\pi\mathsf{G}\Big[
-\rho g +(\mathrm{grad}\rho|\mathrm{grad}\mathcal{K}g) \Big],
$$
where $\rho g \in C_0^{\infty}(B_R)$,  $\mathrm{grad}\mathcal{K}g
\in C^1(\overline{B_R})$ and $\mathrm{grad}\rho \in L^2(B_R, \mathsf{d}^{-\nu+1})$ since $\nu >1$. \\

If $g_1, g_2 \in C_0^{\infty}(B_R)$, then
\begin{equation}
(\mathcal{N} g_1|g_2)_{\mathfrak{G}}=
Q(g_1,g_2)+Q_{01}(g_1,g_2),
\end{equation}
where
\begin{align}
Q_{01}(g_1,g_2)&:=(\mathcal{N}_{01}g_1|g_2)_{\mathfrak{G}}
=4\pi\mathsf{G}(\mathrm{div}(\rho\mathrm{grad}
\mathcal{K}g_1)|g_2)_{\mathfrak{G}} \\
&=4\pi\mathsf{G}((\mathrm{grad}\rho|\mathrm{grad}\mathcal{K}g_1)|
g_2)_{\mathfrak{G}}
-4\pi\mathsf{G}(\rho g_1|g_2)_{\mathfrak{G}},
\end{align}
since we have the identity
\begin{align*}
\mathcal{N}_{01}g&=4\pi\mathsf{G}
\mathrm{div}(\rho\mathrm{grad}\mathcal{K}g) \\
&=4\pi\mathsf{G}(\mathrm{grad}\rho|\mathrm{grad}\mathcal{K}g)+
\rho\triangle 4\pi\mathsf{G}\mathcal{K}g \\
&=4\pi\mathsf{G}(\mathrm{grad}\rho|\mathrm{grad}\mathcal{K}g)
-4\pi\mathsf{G}\rho g.
\end{align*}

We claim that
\begin{equation}
|Q_{01}[g]|=|Q_{01}(g,g)|\leq C\|g\|_{\mathfrak{G}}^2. \label{5.35}
\end{equation}

By putting
\begin{equation}
 \Psi=-\mathcal{K}g,
\end{equation}
we have
\begin{align*}
|\Psi(\vec{x})|&=\Big|
\frac{1}{4\pi}\int
\frac{g(\vec{x}')}{\|\vec{x}-\vec{x}'\|}d\vec{x}'\Big| \\
&\leq \frac{1}{4\pi}
\sqrt{\int_{\|\vec{x}'\|\leq R}\frac{d\vec{x}'}{\|\vec{x}-\vec{x}'\|^2}}
\cdot \|g\|_{L^2(B_R)} \\
&\lesssim \sqrt{\|\vec{x}\|+R}\|g\|_{L^2(B_R)},
\end{align*}
therefore
$$\|\Psi\|_{L^2(B_R)}
\lesssim \sqrt{R}\|g\|_{L^2(B_R)}.
$$
It follows that
\begin{align*}
\|\mathrm{grad} \Psi\|_{L^2(B_R)}^2
&\leq \|\mathrm{grad} \Psi\|_{L^2(\mathbb{R}^3)}^2 \\
&=\int_{\mathbb{R}^3}
(\mathrm{grad} \Psi|\mathrm{grad} \Psi)d\vec{x} =
-\int_{\mathbb{R}^3}(\triangle \Psi|\Psi)d\vec{x} \\
&=-\int_{B_R}(g|\Psi)d\vec{x}
\lesssim \sqrt{R} \|g\|_{L^2(B_R)}^2.
\end{align*}\\

Now we see
\begin{align*}
(*)&:=\Big|((\mathrm{grad}\rho | \mathrm{grad}(\mathcal{K}g_1))|g_2)_{\mathfrak{G}}\Big| \\
&=
\Big|\int
(\mathrm{grad}\rho|\mathrm{grad}\Psi_1)g_2^*
\frac{1}{\rho}\frac{dP}{d\rho}d\vec{x}\Big| \\
&=\Big|\int
(\mathrm{grad}\rho|\mathrm{grad}\Psi_1)g_2^*
\frac{du}{d\rho}d\vec{x}\Big| \\
&=
\Big|\int(\mathrm{grad}u|\mathrm{grad}\Psi_1)g_2^*{d}\vec{x} \Big| \\
&\lesssim \|\mathrm{grad} \Psi_1\|_{L^2(B_R)}
\|g_2\|_{L^2(B_R)} \\
&\lesssim
\|g_1\|_{L^2(B_R)} 
\|g_2\|_{L^2(B_R)},
\end{align*}
since $\mathrm{grad} u $ is bounded and
$$(\mathrm{grad}u|\mathrm{grad}\Psi_1)=
(\nabla u|\nabla \Psi_1)_{\mathbb{R}^3}.$$

However, since $\displaystyle \frac{1}{\rho}\frac{dP}{d\rho}\geq \frac{1}{C}$, we have
$$\|g\|_{L^2(B_R)}\lesssim \|g\|_{\mathfrak{G}}.$$

Thus
$$(*)\lesssim \|g_1\|_{\mathfrak{G}}\|g_2\|_{\mathfrak{G}}.$$
This shows \eqref{5.35}.

Therefore, taking $\kappa$ sufficiently large, we have the bilinear form
\begin{equation}
\mathcal{Q}_1(g_1,g_2)=Q(g_1,g_2)+
Q_{01}(g_1,g_2)+\kappa (g_1|g_2)_{\mathfrak{G}}
\end{equation}
associated with the operator $\mathcal{A}+\kappa $ satisfies
$$\mathcal{Q}_1[g]\geq \|g\|_{\mathfrak{G}}^2$$
and $\{ g\in \mathfrak{G}_1 | \mathcal{Q}_1[g]\leq 1\}$ is a compact subset of
$\mathfrak{G}$. Summing up, we have: \\

The operator $\mathcal{A}$ has the Friedrichs extension
$\mathcal{T}$, which is a self-adjoint operator in $\mathfrak{G}$, and
the resolvent $(\mathcal{T}+\kappa)^{-1}$ is a 
compact operator in $\mathfrak{G}$. The spectrum of $\mathcal{T}$ is of the Sturm-Liouville type. \\
 
 Hereafter let us use the letter $\mathcal{N}$ for the Friedrichs extension $\mathcal{T}$. Thus we have the following

\begin{Theorem}\label{Th.2}
The operator $\mathcal{N}$ is a self-adjoint operator bounded from below in the Hilbert space $\mathfrak{G}$ and its spectrum $\sigma(\mathcal{N})$ is of the Sturm-Liouville type.
\end{Theorem} 

\section{Spectral analysis of $\mathbf{L}$}

We go back to the spectral analysis of $\mathbf{L}$.

\subsection{From $\sigma(\mathcal{N})$ to $\sigma(\mathbf{L})$}

Since we deal with the variable $\hat{\mbox{\boldmath$\xi$}}=\rho\mbox{\boldmath$\xi$}$ in this section, we shall use the following notations: 

\begin{Definition}

1) We define the Hilbert space $\hat{\mathfrak{H}}$ of functions on $B_R$
endowed with the inner product
\begin{equation}
(\hat{\mbox{\boldmath$\xi$}}_1|\hat{\mbox{\boldmath$\xi$}}_2)_{\hat{\mathfrak{H}}}=
\int
(\hat{\mbox{\boldmath$\xi$}}_1|\hat{\mbox{\boldmath$\xi$}}_2)\frac{1}{\rho}d\vec{x}
\end{equation}
so that $\mbox{\boldmath$\xi$}\in\mathfrak{H} \Leftrightarrow \rho \mbox{\boldmath$\xi$}\in \hat{\mathfrak{H}}$ and
$$(\rho \mbox{\boldmath$\xi$}_1|\rho\mbox{\boldmath$\xi$}_2)_{\hat{\mathfrak{H}}}=
(\mbox{\boldmath$\xi$}_1|\mbox{\boldmath$\xi$}_2)_{\mathfrak{H}}.$$

2) We put
\begin{equation}
\hat{\mathfrak{F}}=\{ \hat{\mbox{\boldmath$\xi$}} \in \hat{\mathfrak{H}} \  |\  
\mathrm{div}\hat{\mbox{\boldmath$\xi$}} \in \mathfrak{G} \quad
\mbox{in the distribution sense} \},
\end{equation}
defining the inner product
\begin{equation}
(\hat{\mbox{\boldmath$\xi$}}_1|\hat{\mbox{\boldmath$\xi$}}_2)_{\hat{\mathfrak{F}}}
=(\hat{\mbox{\boldmath$\xi$}}_1|\hat{\mbox{\boldmath$\xi$}}_2)_{\hat{\mathfrak{H}}}+
(\mathrm{div}\hat{\mbox{\boldmath$\xi$}}_1|\mathrm{div}\hat{\mbox{\boldmath$\xi$}}_2)_{\mathfrak{G}}.
\end{equation}
\end{Definition}
Here `$g=\mathrm{div}\hat{\mbox{\boldmath$\xi$}}$ in the distribution sense' means 
\begin{equation}
(g|\varphi)_{\mathfrak{G}}=-\int
\Big(\hat{\mbox{\boldmath$\xi$}}\Big|\mathrm{grad}\Big(\frac{1}{\rho}\frac{dP}{d\rho}\varphi\Big)\Big)d\vec{x}
\quad\forall \varphi
\in {C}_0^{\infty}(B_R). \label{5.23bis}
\end{equation}\\

Let us put
\begin{equation}
\hat{\mathbf{L}}\hat{\mbox{\boldmath$\xi$}}=
\hat{\mathbf{M}}(\mathrm{div}\hat{\mbox{\boldmath$\xi$}}),
\end{equation}
for 
\begin{equation}
\hat{\mbox{\boldmath$\xi$}}\in \mathsf{D}(\hat{\mathbf{L}})
=\{ \hat{\mbox{\boldmath$\xi$}}\in\hat{\mathfrak{F}}\  |\  
\mathrm{div}\hat{\mbox{\boldmath$\xi$}}\in \mathsf{D}(\mathcal{N})\}.
\end{equation}\\

Note that
\begin{equation}
\rho{\mathbf{L}}\mbox{\boldmath$\xi$}=\hat{\mathbf{L}}(\rho \mbox{\boldmath$\xi$})\quad\mbox{or}\quad
\rho{\mathbf{L}}\Big(\frac{1}{\rho}\hat{\mbox{\boldmath$\xi$}}\Big)
=\hat{\mathbf{L}}(\hat{\mbox{\boldmath$\xi$}}).
\end{equation}
and, 
as for the operator $\hat{\mathbf{M}}$, we note following

\begin{Proposition}\label{Prop.L0}
The operator $\hat{\mathbf{M}}$ maps $\mathfrak{G}_1$ into $\hat{\mathfrak{H}}$.
\end{Proposition}

Proof. Consider
$$\hat{\mathbf{M}}g=L_1+4\pi\mathsf{G}L_2, $$
where
$$
 L_1=\rho\mathrm{grad}\Big(-\frac{1}{\rho}\frac{dP}{d\rho}g\Big), \quad
 L_2=\rho\mathrm{grad}\mathcal{K}g.
$$

Since
$$\mathrm{grad}\Big(-\frac{1}{\rho}\frac{dP}{d\rho}g\Big)=
-\frac{1}{\rho}\frac{dP}{d\rho}
\mathrm{grad}g+
\mathrm{grad}\Big(-\frac{1}{\rho}\frac{dP}{d\rho}\Big)
\cdot g,
$$
and $$\mathrm{grad}\Big(-\frac{1}{\rho}\frac{dP}{d\rho}\Big) \lesssim\frac{1}{\rho},$$ 
we see
\begin{align*}
\|L_1\|_{\hat{\mathfrak{H}}}^2& \lesssim 
\int \frac{1}{\rho}\Big(\frac{dP}{d\rho}\Big)^2\|\mathrm{grad}g\|^2d\vec{x} 
+\int \frac{1}{\rho}|g|^2d\vec{x} \\
&\lesssim \|g\|_{W^{1,2}(\mathsf{d}^{-\nu}, \mathsf{d}^{2-\nu})}^2 \\
&\lesssim \|g\|_{W^{1,2}(\mathsf{d}^{1-\nu}, \mathsf{d}^{2-\nu})}^2 \\
& \lesssim \|g\|_{\mathfrak{G}_1}^2
\end{align*}

Therefore $$\|L_1\|_{\hat{\mathfrak{H}}}\lesssim \|g\|_{\mathfrak{G}_1}.$$

As for $L_2$, putting
$$\Psi=-\mathcal{K}g, $$
we see
\begin{align*}
\|L_2\|_{\hat{\mathfrak{H}}}^2&=
\int \rho \|\mathrm{grad}\Psi\|^2d\vec{x} \\
&\lesssim \rho_{\mathsf{O}}\|\mathrm{grad}\Psi\|_{L^2(B_R)}^2
\lesssim \rho_{\mathsf{O}}\sqrt{R}
\|g\|_{L^2(B_R)}^2 \\
&\lesssim \|g\|_{\mathfrak{G}}^2.
\end{align*}
Summing up, we have
$$\|\hat{\mathbf{M}}g\|_{\hat{\mathfrak{H}}}\lesssim
\|g\|_{\mathfrak{G}_1}.
$$
This completes the proof. $\square$\\

Take an arbitrary $\lambda \not=0$
from the resolvent set $\varrho(\mathcal{N})$ of the operator $\mathcal{N}$ and consider the resolvent
$(\lambda - \mathcal{N})^{-1}$, which is a bounded linear operator
from $\mathfrak{G}$ into $\mathfrak{G}_1$.
Then the equation
\begin{equation}
(\lambda -\hat{\mathbf{L}})\hat{\mbox{\boldmath$\xi$}}=\hat{\mathbf{f}}\in \hat{\mathfrak{F}}
\label{*1}
\end{equation}
can be solved as
\begin{equation}
\hat{\mbox{\boldmath$\xi$}}=\frac{1}{\lambda}
(\hat{\mathbf{f}}+
\hat{\mathbf{M}}(\lambda-\mathcal{N})^{-1}\mathrm{div}
\hat{\mathbf{f}}). \label{*2}
\end{equation}
In fact, \eqref{*2} implies
$$\lambda\hat{\mbox{\boldmath$\xi$}}=\hat{\mathbf{f}}+
\hat{\mathbf{M}}(\lambda-\mathcal{N})^{-1}
\mathrm{div}\hat{\mathbf{f}}; $$
putting $f=(\lambda-\mathcal{N})^{-1}\mathrm{div}\hat{\mathbf{f}}$,
we have
$$\mathrm{div}\hat{\mathbf{f}}=(\lambda-\mathcal{N})f,\quad
\lambda\hat{\mbox{\boldmath$\xi$}}=\hat{\mathbf{f}}+
\hat{\mathbf{M}}f $$
so that
\begin{align*}
\lambda\mathrm{div}\hat{\mbox{\boldmath$\xi$}}&=\mathrm{div}\hat{\mathbf{f}}+
\mathrm{div}\hat{\mathbf{M}}f \\
&=(\lambda-\mathcal{N})f+
\mathcal{N}f =\lambda f;
\end{align*}
therefore $\mathrm{div}\hat{\mbox{\boldmath$\xi$}}=f$, since $\lambda\not=0$; then we have
\begin{align*}
\lambda\hat{\mbox{\boldmath$\xi$}}&=\hat{\mathbf{f}}+
\hat{\mathbf{L}}(
\lambda-\mathcal{N})^{-1}\mathrm{div}\hat{\mathbf{f}} \\
&=\hat{\mathbf{f}}+\hat{\mathbf{M}}f \\
&=\hat{\mathbf{f}}+\hat{\mathbf{M}}\mathrm{div}\hat{\mbox{\boldmath$\xi$}} \\
&=\hat{\mathbf{f}}+\hat{\mathbf{L}}\hat{\mbox{\boldmath$\xi$}},
\end{align*}
that is \eqref{*1}.

Since we know by Proposition \ref{Prop.L0} that the operator $\hat{\mathbf{M}}(\lambda-\mathcal{N})^{-1}$
is a bounded linear operator from $\mathfrak{G}$ into
$\hat{\mathfrak{H}}$
 and since 
\begin{align*}
\|\mathrm{div}\hat{\mathbf{M}}(\lambda-\mathcal{N})^{-1}f\|_{\mathfrak{G}}&
=\|\mathcal{N}(\lambda-\mathcal{N})^{-1}f\|_{\mathfrak{G}} \\
&\leq \||\lambda(\lambda-\mathcal{N})^{-1}-I|\|\|f\|_{\mathfrak{G}},
\end{align*}
the operator $\hat{\mathbf{M}}(\lambda-\mathcal{N})^{-1}$ is a bounded operator from $\mathfrak{G}$ into 
$\hat{\mathfrak{F}}$. Therefore 
we can claim the following

\begin{Proposition}\label{Prop.O1}
If $\lambda \in \varrho(\mathcal{N}) \setminus\{0\}$, then the operator
$(\lambda-\hat{\mathbf{L}})^{-1}$ is a bounded linear operator from
$\hat{\mathfrak{F}}$ into $\hat{\mathfrak{F}}$, that is, $\lambda$ belongs to $\varrho(\hat{\mathbf{L}})$, the resolvent set of $\hat{\mathbf{L}}$.
\end{Proposition}

We can claim

\begin{Proposition}
The operator $\hat{\mathbf{L}}$ considered in the Hilbert space $\hat{\mathfrak{F}}$ is self-adjoint.
\end{Proposition}

Proof. First we note that
$$\mathcal{N}g=\mathrm{div}(\hat{\mathrm{M}}g) $$
for $g \in \mathsf{D}(\mathcal{N})$.Let $\mbox{\boldmath$\xi$}_{\mu} \in \mathsf{D}(\hat{\mathbf{L}})$ with
$g_{\mu} =\mathrm{div}\hat{\mbox{\boldmath$\xi$}}_{\mu} \in \mathsf{D}(\mathcal{N})$ , $\mu=1,2$.  Then
\begin{align*}
(\hat{\mathbf{L}}\mbox{\boldmath$\xi$}_1|\hat{\mbox{\boldmath$\xi$}}_2)_{\hat{\mathfrak{F}}}&=
(\hat{\mathrm{M}}g_1|\hat{\mbox{\boldmath$\xi$}}_2)_{\hat{\mathfrak{H}}}+(\mathcal{N}g_1|g_2)_{\mathfrak{G}} \\
&=(\rho\mathrm{grad}G_1|\hat{ \mbox{\boldmath$\xi$} }_2)_{\hat{\mathfrak{H}}}+(\mathcal{N}g_1|g_2)_{\mathfrak{G}} \\
&=\int (\mathrm{grad}G_1|\hat{\mbox{\boldmath$\xi$}}_2) +(\mathcal{N}g_1|g_2)_{\mathfrak{G}} \\
&=-\int G_1 g_2^* +(\mathcal{N}g_1|g_2)_{\mathfrak{G}} \\
&=\int \Big[\frac{1}{\rho}\frac{dP}{d\rho}g_1g_2^*-4\pi\mathsf{G}(\mathcal{K}g_1)g_2^*\Big]
+(\mathcal{N}g_1|g_2)_{\mathfrak{G}} \\
&=(\hat{\mbox{\boldmath$\xi$}}_1|\hat{\mathbf{L}}\hat{\mbox{\boldmath$\xi$}}_2)_{\hat{\mathfrak{F}}},
\end{align*}
since $\mathcal{K}$ is symmetric and $\mathcal{N}$ is self-adjoint in $\mathfrak{G}$.
Thus  $\hat{\mathbf{L}}$ is symmetric. 

$\hat{\mathbf{L}}$ is closable. In fact, let $\hat{\mbox{\boldmath$\xi$}}_n \in \mathsf{D}(\hat{\mathbf{L}}) \rightarrow 0$ in $\hat{\mathfrak{F}}$ and 
$\hat{\mathbf{L}}\hat{\mbox{\boldmath$\xi$}}_n \rightarrow \hat{\mathbf{f}}$ in
$\hat{\mathfrak{F}}$ as $ n\rightarrow \infty$. Then $g_n:=\mathrm{div}\hat{\mbox{\boldmath$\xi$}}_n\rightarrow 0$ in $\mathfrak{G}$, $g_n \in \mathsf{D}(\mathcal{N})$, and
\begin{align*}
\|\mathcal{N}g_n -\mathrm{div}\hat{\mathbf{f}}\|_{\mathfrak{G}}&=\|\mathrm{div}\hat{\mathbf{M}}g_n-\mathrm{div}\hat{\mathbf{f}}\|_{\mathfrak{G}} \\
&=\|\mathrm{div}(\hat{\mathbf{L}}\hat{\mbox{\boldmath$\xi$}}_n-\hat{\mathbf{f}})\|_{\mathfrak{G}} \\
&\leq \|\hat{\mathbf{L}}\hat{\mbox{\boldmath$\xi$}}_n-\hat{\mathbf{f}}\|_{\hat{\mathfrak{F}}} \rightarrow 0.
\end{align*}
Since $\mathcal{N}$ is closed, we can claim $\mathrm{div}\hat{\mathbf{f}}=0$, that is, $\mathcal{N}g_n \rightarrow 0$ in $\mathfrak{G}$.
Thus 
$$\|g_n\|_{\mathfrak{G}_1}^2 \leq ((\mathcal{N}+\kappa)g_n|g_n)_{\mathfrak{G}} \rightarrow 0, $$
where $\kappa$ is a sufficiently large constant. Therefore Proposition \ref{Prop.L0} 
says that $\hat{\mathbf{L}}\hat{\mbox{\boldmath$\xi$}}_n=\hat{\mathbf{M}}g_n \rightarrow 0$ in $\hat{\mathfrak{H}}$. Of course $\mathrm{div}\hat{\mathbf{L}}\mbox{\boldmath$\xi$}_n=\mathcal{N}g_n $, which $
\rightarrow 0$ in $\mathfrak{G}$. Thus $\hat{\mathbf{L}}\mbox{\boldmath$\xi$}_n \rightarrow 0$ in $\hat{\mathfrak{F}}$, that is, $\hat{\mathbf{f}}=0$. This means that $\hat{\mathbf{L}}$ is closable.

But, since $\varrho(\hat{\mathbf{L}})\not=\emptyset$ due to Proposition
\ref{Prop.O1}, \cite[Chapter III, Theorem 3.16]{Kato} guarantees that $\hat{\mathbf{L}}$ 
is self-adjoint.
$\square$.\\


We also claim

\begin{Proposition}\label{Prop.O2}
A non-zero eigenvalue of $\mathcal{N}$ is an eigenvalue of $\hat{\mathbf{L}}$, and
$\sigma(\mathcal{N})\setminus\{0\} \subset \sigma(\hat{\mathbf{L}})$. 
\end{Proposition}

Proof. Let $\lambda \not=0$ be an eigenvalue of $\mathcal{N}$ and let $g$ be an associated eigenfunction 
of $\mathcal{N}$ such that $\|g\|_{\mathfrak{G}}=1$. Put
$$\hat{\mbox{\boldmath$\xi$}}=\frac{1}{\lambda}\hat{\mathbf{M}}g. $$
Then $\hat{\mbox{\boldmath$\xi$}}\in \hat{\mathfrak{F}}$ and 
$\lambda\hat{\mbox{\boldmath$\xi$}}=\hat{\mathbf{L}}g$. Thus
$$\lambda\mathrm{div}\hat{\mbox{\boldmath$\xi$}}=\mathcal{N}g=\lambda g $$
so that $\mathrm{div}\hat{\mbox{\boldmath$\xi$}}=g$, since $\lambda\not=0$. Thus
$$\lambda \hat{\mbox{\boldmath$\xi$}}=\hat{\mathbf{M}}\mathrm{div}\hat{\mbox{\boldmath$\xi$}}=\hat{\mathbf{L}}\hat{\mbox{\boldmath$\xi$}},$$
that is, $\lambda$ is an eigenvalue of
$\hat{\mathbf{L}}$. $\square$\\

We should note the following 

\begin{Proposition}\label{Prop.NC}
The kernel of the operator $\hat{\mathbf{L}}$ is infinitely dimensional, and therefore
a resolvent $(\lambda -\hat{\mathbf{L}})^{-1}$ of the operator $\hat{\mathbf{L}}$ cannot be a compact operator and the spectrum $\sigma(\hat{\mathbf{L}})$ is not of the Sturm-Liouville type.
\end{Proposition}

Proof. The space
$$\hat{\mathfrak{N}}:=\{ \hat{\mbox{\boldmath$\xi$}}\in \hat{\mathfrak{H}}\   |\   \mathrm{div}\hat{\mbox{\boldmath$\xi$}}=0
\quad\mbox{in the distribution sense} \}$$
is infinitely dimensional, since
$\hat{\mbox{\boldmath$\xi$}}=
\mbox{curl}\mathbf{A}
$
belongs to $\hat{\mathfrak{N}}$ for arbitrary $\mathbf{A} \in {C}_0^{\infty}(B_R)$. Clearly
$\hat{\mathfrak{N}} \subset \mbox{Ker}\hat{\mathbf{L}}$, that is, any
function $\hat{\mbox{\boldmath$\xi$}}\in \hat{\mathfrak{N}} $ with
$\|\hat{\mbox{\boldmath$\xi$}}\|_{\hat{\mathfrak{H}}}\not=0$ is an eigenfunction of
$\hat{\mathbf{L}}$ associated with the eigenvalue 0. $\square$.\\

\begin{Remark}\label{Rem.as}
Even if we restrict ourselves to axisymmetric perturbations, the situation is the same, that is,
$$\hat{\mathfrak{N}}^{\mathsf{as}}:=
\{ \hat{\mbox{\boldmath$\xi$} }\in \hat{\mathfrak{H}} | \ \mbox{\boldmath$\xi$}\quad\mbox{is axisymmetric},\quad \mathrm{div}\hat{\mbox{\boldmath$\xi$}}=0 \}$$
is infinite dimensional.

 In fact, for arbitrary $A: (r,\zeta) \mapsto A(r,\zeta)  \in C_0^{\infty}([0,R]\times ]-1,1[)$
$$\hat{\mbox{\boldmath$\xi$}}=\Big(r\frac{\partial A}{\partial\zeta}\Big)\frac{\partial}{\partial r}+
\Big(-\frac{1}{r^2}\frac{\partial}{\partial r}(r^3A)\Big)\frac{\partial}{\partial\zeta}$$
belongs to $\hat{\mathfrak{N}}^{\mathsf{as}}$. Here $\zeta=x^3/r$. 
\end{Remark}

However we can claim that $\sigma(\hat{\mathbf{L}})\setminus \{0\}$ consists of eigenvalues of finite multiplicity, that is, we have

\begin{Proposition}\label{A0FE}
If $\lambda \in \sigma(\hat{\mathbf{L}})$ and if $\lambda \not=0$, then
$\lambda$ is an eigenvalue of the operator
$\hat{\mathbf{L}}$ of finite multiplicity.
\end{Proposition}

Poof. Let $\lambda \in \sigma(\hat{\mathbf{L}}) \setminus \{0\}$.
 By Proposition \ref{Prop.O1}, we know that $\sigma(\hat{\mathbf{L}})\subset
\sigma(\mathcal{N}) \cup \{0\}$. Therefore we have $\lambda \in
\sigma(\mathcal{N})$ so that $\lambda$ is an eigenvalue of $\mathcal{N}$ with finite multiplicity.
On the other hand,
$$\hat{\mbox{\boldmath$\xi$}}\in \mbox{Ker}(\lambda-\hat{\mathbf{L}})$$
if and only if
$$
\lambda \hat{\mbox{\boldmath$\xi$}}=\hat{\mathbf{L}}\hat{\mbox{\boldmath$\xi$}}=
\hat{\mathbf{M}}\mathrm{div}\hat{\mbox{\boldmath$\xi$}}.$$

If we put $g=\mathrm{div}\hat{\mbox{\boldmath$\xi$}}$, this implies
$\lambda g=\mathcal{N}g$, that is, $g \in
\mbox{Ker}(\lambda-\mathcal{N})$ and
$$\hat{\mbox{\boldmath$\xi$}}=\frac{1}{\lambda}\hat{\mathbf{M}}g. $$
Hence we can 
claim
$$\hat{\mbox{\boldmath$\xi$}}\in \mbox{Ker}(\lambda-\hat{\mathbf{L}})$$
if and only if
$$\hat{\mbox{\boldmath$\xi$}}=\frac{1}{\lambda}\hat{\mathbf{M}}g\quad\mbox{with}\quad
g \in \mbox{Ker}(\lambda-\mathcal{N}). $$
Therefore
$$\mbox{dim.Ker}(\lambda-\hat{\mathbf{L}})
\leq \mbox{dim.Ker}(\lambda -\mathcal{N}).$$
$\square$\\

Summing up, we have \\

{\it 
The operator $\hat{\mathbf{L}}$ is a self-adjoint operator in 
$\hat{\mathfrak{F}}$. Its spectrum $\sigma(\hat{\mathbf{L}})$ coincides with
$\sigma(\mathcal{N})\cup\{0\}$, while $\mathrm{dim}.\mathrm{Ker}(\hat{\mathbf{L}})=\infty$ and $\lambda \in \sigma(\hat{\mathbf{L}})
\setminus \{0\}$ is an eigenvalue of finite multiplicity. } \\

Translating this statement to that in terms in $\mathbf{L}$, we can claim the following

\begin{Theorem} \label{Th.3}
The operator $\mathbf{L}$ is a self-adjoint operator in 
$\mathfrak{F}$. Its spectrum $\sigma({\mathbf{L}})$ coincides with
$\sigma(\mathcal{N})\cup\{0\}$, while $\mathrm{dim}.\mathrm{Ker}({\mathbf{L}})=\infty$ and $\lambda \in \sigma({\mathbf{L}})
\setminus \{0\}$ is an eigenvalue of finite multiplicity.
\end{Theorem}

Here the Hilbert space $\mathfrak{F}$, which is nothing but $\frac{1}{\rho}\hat{\mathfrak{F}}$,  is defined by the following

\begin{Definition}\label{Def_F}
 We put
\begin{equation}
{\mathfrak{F}}=\{ {\mbox{\boldmath$\xi$}} \in {\mathfrak{H}} \  |\  
\mathrm{div}(\rho{\mbox{\boldmath$\xi$}}) \in \mathfrak{G} \quad
\mbox{in the distribution sense} \},
\end{equation}
defining the inner product
\begin{equation}
({\mbox{\boldmath$\xi$}}_1|{\mbox{\boldmath$\xi$}}_2)_{{\mathfrak{F}}}
=({\mbox{\boldmath$\xi$}}_1|{\mbox{\boldmath$\xi$}}_2)_{{\mathfrak{H}}}+
(\mathrm{div}(\rho{\mbox{\boldmath$\xi$}}_1)|\mathrm{div}(\rho{\mbox{\boldmath$\xi$}}_2))_{\mathfrak{G}}.
\end{equation}
\end{Definition}

Therefore $\mathbf{L}$ here stands for the Friedrichs extension of $\mathfrak{A}$ with domain
$C_0^{\infty}(B_R)$ in $\mathfrak{F}$, and is different from the Friedrichs extension in $\mathfrak{H}$ considered in Section 2. But this diverting of symbol may not lead confusions. 

Moreover, we note that $\lambda \in \varrho(\mathbf{L})$ if and only if 
$\lambda \in \varrho(\hat{\mathbf{L}})$, while 
$$ (\lambda-\mathbf{L})^{-1}\mathbf{f}=
\frac{1}{\rho}(\lambda-\hat{\mathbf{L}})^{-1}(\rho\mathbf{f}) $$
for $\lambda \in \varrho(\mathbf{L})$ and $\mathbf{f} \in \mathfrak{F}$.\\


\subsection{Non-time-periodic special solutions}

Now the system \eqref{5.10a}\eqref{5.10b} is split to 
\begin{equation}
\frac{\partial^2\hat{\mbox{\boldmath$\xi$}}}{\partial t^2}+\hat{\mathbf{M}}g=0 \label{5.38}
\end{equation}
and
\begin{equation}
\frac{\partial^2g}{\partial t^2}+\mathcal{N}g=0. \label{5.40}
\end{equation}

The equation \eqref{5.38} for $\hat{\mbox{\boldmath$\xi$}}$ is integrated as
\begin{equation}
\hat{\mbox{\boldmath$\xi$}}=t\overset{\circ}{\hat{\mathbf{v}}}-\int_0^t(t-s)\hat{\mathbf{M}}g(s)ds,
\end{equation}
where $\overset{\circ}{\hat{\mathbf{v}}}$ is the initial data $\partial\hat{\mbox{\boldmath$\xi$}}/\partial t|_{t=0}$, 
provided that $g \in C([0,+\infty[; \mathfrak{G}_1) $ is given. Therefore it is sufficient to consider the single equation \eqref{5.40} for $g$.

But the functional analysis of the partial differential operator $\mathcal{N}$ is already done. 

Let $\lambda$ be a positive eigenvalue of the operator
$\mathcal{N}$ and $\varphi(\vec{x})$ be an eigenfunction of $\mathcal{N}$ associated with this eigenvalue $\lambda$ such that $\|\varphi\|_{\mathfrak{G}}=1$. Then, for an arbitrary constant $E$, 
\begin{equation}
g=E\sin(\sqrt{\lambda}t)\varphi(\vec{x})
\end{equation}
is the time periodic solution of \eqref{5.40} with the initial conditions
\begin{equation}
g|_{t=0}=0,\qquad
\frac{\partial g}{\partial t}\Big|_{t=0}=E\sqrt{\lambda}\varphi(\vec{x}).
\end{equation}
Let the initial data $\overset{\circ}{\hat{\mathbf{v}}}$ for $\partial\hat{\mbox{\boldmath$\xi$}}/\partial t$ be such that
\begin{equation}
\mathrm{div}\overset{\circ}{\hat{\mathbf{v}}} =E\sqrt{\lambda}\varphi(\vec{x}). \label{5.44}
\end{equation}
Then the solution $\hat{\mbox{\boldmath$\xi$}}$ of \eqref{5.38} to the initial conditions is
\begin{align}
\hat{\mbox{\boldmath$\xi$}}&=t\overset{\circ}{\hat{\mathbf{v}}}-
\Big[\int_0^t(t-s)\sin (\sqrt{\lambda}t)ds\Big]\hat{\mathbf{M}}\varphi \nonumber \\
&=
t\overset{\circ}{\hat{\mathbf{v}}}+E\Big(\frac{1}{\lambda}\sin(\sqrt{\lambda}t)-\frac{t}{\sqrt{\lambda}}\Big)
\hat{\mathbf{M}}\varphi. \label{5.45}
\end{align}

Note that it seems to be possible that $\hat{\mbox{\boldmath$\xi$}}$ given by \eqref{5.45} is not time periodic. 
In fact, \eqref{5.45} can be written as 
\begin{equation}
\hat{\mbox{\boldmath$\xi$}}=\mathbf{B}t +\frac{E}{\lambda}\sin(\sqrt{\lambda}t)
\hat{\mathbf{M}}\varphi, \label{5.46}
\end{equation}
where 
\begin{equation}
\mathbf{B}=\overset{\circ}{\hat{\mathbf{v}}}-\frac{E}{\sqrt{\lambda}}\hat{\mathbf{M}}
\varphi.
\end{equation}\\

Let us note that, since $\mathrm{div}\hat{\mathbf{M}}=\mathcal{N}$ and
$\mathcal{N} \varphi=\lambda \varphi$, \eqref{5.44} implies 
\begin{equation}
\mathrm{div}\mathbf{B}=0.
\end{equation}

This means that there is a vector function $\mathbf{A}$ such that
\begin{equation}
\mathbf{B}=\mbox{curl}\mathbf{A}. \label{5.49}
\end{equation}

Let us consider the case in which the initial perturbation is spherically symmetric, that is, $\overset{\circ}{v}=\overset{\circ}{v}(r), \overset{\circ}{w}=0$, and  the eigenfunction is spherically symmetric, that is, $\varphi=\varphi(r)$, then $\mathbf{B}$ is of the form
$$
\mathbf{B}=
rb(r) \frac{\partial}{\partial r},
$$
and $\mathrm{div}\mathbf{B}=\displaystyle\frac{1}{r^2}\frac{d}{dr}(r^3 b)=0 $
implies $r^3b(r)=\mbox{Const.}$ on $0\leq r<R$. 
Since $\displaystyle \Big(\frac{1}{r^3},0\Big)^T \not\in \hat{\mathfrak{H}}$, 
we have
$b(r)=0$ and $\mathbf{B}=0$. In this case the solution $\hat{\mbox{\boldmath$\xi$}}$ given by \eqref{5.46} is time periodic.\\

However $\mathbf{B}\not=0$ in general. Actually we can take the initial data
$\overset{\circ}{\hat{\mathbf{v}}}$ so that
$$
\overset{\circ}{\hat{\mathbf{v}}} =\frac{E}{\sqrt{\lambda}}\hat{\mathbf{M}}\varphi
+\mbox{curl}\mathbf{A}
$$
with an arbitrary vector field $\mathbf{A} \in C_0^{\infty}(B_R)$ which
does not identically vanish. Then $\mathbf{B}\not=\vec{0}$ and
the solution $\hat{\mbox{\boldmath$\xi$}}$ given by \eqref{5.46} is not time periodic
but grows linearly in time $t$.\\

Let us suppose that $P=\mathsf{A}\rho^{\gamma}, \frac{6}{5} < \gamma <2$ for the simplicity.\\

If $\gamma <4/3$, there is a negative eigenvalue $\lambda$ for the operator 
$\mathcal{N}$ restricted to the space
of spherically symmetric functions
and the solution
$$\hat{\mbox{\boldmath$\xi$}}=\mathbf{B}t-
\frac{E}{\lambda}e^{-\lambda t}\hat{\mathbf{M}}\varphi
$$
exponentially grows in $t$ even if $\mathbf{B}=0$. ( See the discussion below.) Here $E\not=0, \|\varphi\|_{\mathfrak{G}}=1, \mathcal{N}\varphi=\lambda\varphi$.
For this case the non-linear instability in the sense of Ljapunov of the spherically
symmetric equilibria was rigorously proven in \cite{JJ1}.

But if we consider perturbations which are not spherically symmetric, we have a particular solution
$$\hat{\mbox{\boldmath$\xi$}}=t\overset{\circ}{\hat{\mathbf{v}}}$$
with $\mathrm{div}\overset{\circ}{\hat{\mathbf{v}}}=0$ but
$\overset{\circ}{\hat{\mathbf{v}}}\not=\vec{0}$, e.g.,
$$
\overset{\circ}{\hat{\mathbf{v}}}=
\mbox{curl}\mathbf{A}
$$
with an arbitrary vector field $\mathbf{A} \in C_0^{\infty}(B_R)$ which does not identically vanish. This solution is growing linearly in $t$,
and this phenomenon appears for all $6/5 <\gamma <2$, say, even for $4/3 \leq \gamma$.
This may imply nonlinear instability of the spherically symmetric equilibria with respect to
thoroughly general perturbations. \\

Let us recall the well-known result for spherically symmetric perturbations,
and note that it can be applied to the present problem.\\

Now let us consider
\begin{equation}
\hat{\mathfrak{H}}^{\mathsf{ss}}=\{ \hat{\mbox{\boldmath$\xi$}}\in\hat{\mathfrak{H}} \  |\  
\hat{\mbox{\boldmath$\xi$}}=r\hat{y}(r)\mathbf{e}_r\},
\end{equation}
where $\displaystyle \mathbf{e}_r=\frac{\partial}{\partial r}$.

Then we have
\begin{equation}
\hat{\mathbf{M}}g=
\rho r\Big(
\mathcal{L}^{\mathsf{ss}}y \Big)\mathbf{e}_r
\end{equation}
where
$$\mathcal{L}^{\mathsf{ss}}y=
-\frac{1}{\rho r^4}\frac{d}{dr}
\Big(\gamma r^4P\frac{dy}{dr}\Big)
-\frac{(3\gamma-4)}{r}\frac{du}{dr}y, $$
provided that
$$\hat{\mbox{\boldmath$\xi$}}=
\rho ry(r) \mathbf{e}_r \in \hat{\mathfrak{H}}^{\mathsf{ss}}. $$

Therefore an eigenvalue of $\mathcal{L}^{\mathsf{ss}}$, which is always positive if $\gamma >4/3$ , is an eigenvalue of $\hat{\mathbf{M}}$ and of $\mathcal{N}$. \\

Moreover, using results proven later,  we shall be able to claim the following assertion:\\

{\it If $P=\mathsf{A}\rho^{\gamma}$ with $\gamma > 4/3$, then the least eigenvalue 
$\lambda_1\not=0$ of $\hat{\mathbf{M}}$ is positive.}\\

Proof can be done thanks to Theorem \ref{Th.3} and Theorem \ref{Th.6} (its Corollary) . Actually $\lambda_1^{\mathsf{ss}} > 0$,
$\lambda_1^{\mathsf{ss}}$ being the least eigenvalue of the spherically symmetric problem, for $P=\mathsf{A}\rho^{\gamma}$ with $\gamma > 4/3$.\\



\section{Eigenfunctions given by spherical harmonics}

We are going to find an eigenvalue 
$\lambda$ of the operator $\mathbf{L}$ with an eigenfunction $\mbox{\boldmath$\xi$}$ 
of the form
\begin{equation}
\mbox{\boldmath$\xi$}=r\psi(r)Y_{lm}\mathbf{e}_r+r\chi(r)\nabla_sY_{lm}+
r\kappa(r)\nabla_s^{\perp}Y_{lm} \label{xilm}
\end{equation}
for which
$g=\mathrm{div}(\rho\mbox{\boldmath$\xi$})$ has the form
\begin{equation}
g=\check{g}(r)Y_{lm}(\vartheta, \phi)
\end{equation}
with 
\begin{equation}
\check{g}=\frac{1}{r^2}\frac{\partial}{\partial r}
(r^3\rho\psi)-l(l+1)\rho\chi. \label{Letglm}
\end{equation}\\

Here $l,m \in \mathbb{Z} $ are fixed integers such that $l \geq 0, |m|\leq l$ and the spherical harmonics $Y_{lm}$ for $l,m \in \mathbb{Z}, 0\leq l, |m|\leq l$ are defined by
$$Y_{lm}(\vartheta,\phi)=
\sqrt{\frac{2l+1}{4\pi}\cdot\frac{(l-m)!}{(l+m)!}}
P_l^m(\cos\vartheta)e^{\sqrt{-1}m\phi} $$
for $l=0,1,\cdots, m=0,\cdots, l$ and
$$Y_{l,-m}=(-1)^mY_{lm}^*. $$
See \cite{Jackson}. We use

\begin{Notation} We denote
\begin{align*}
&\nabla_sf:=\frac{\partial f}{\partial\vartheta}\mathbf{e}_{\vartheta}+
\frac{1}{\sin\vartheta}\frac{\partial f}{\partial\phi}\mathbf{e}_{\phi}, \\
&\nabla_s^{\perp}f:= \frac{1}{\sin\vartheta}\frac{\partial f}{\partial\phi}
\mathbf{e}_{\vartheta}
-\frac{\partial f}{\partial\vartheta}\mathbf{e}_{\phi}, \\
&\triangle_sf:=\frac{1}{\sin\vartheta}\frac{\partial}{\partial\vartheta}
\Big(\sin\vartheta\frac{\partial f}{\partial\vartheta}\Big)+
\frac{1}{\sin^2\vartheta}\frac{\partial^2f}{\partial\phi^2},
\end{align*}
while $(r, \vartheta, \phi)$ is the spherical co-ordinate system :
$$ x^1=r\sin\vartheta\cos\phi,\quad
x^2=r\sin\vartheta\sin\phi,\quad x^3=r\cos\vartheta, $$
and $\mathbf{e}_r, \mathbf{e}_{\vartheta}, \mathbf{e}_{\phi}$ are the unit vectors for the 
co-ordinate system, that is,
$$\mathbf{e}_r=\frac{\partial}{\partial r}, \quad
\mathbf{e}_{\vartheta}=\frac{1}{r}\frac{\partial}{\partial\vartheta}, \quad
\mathbf{e}_{\phi}=\frac{1}{r\sin\vartheta}\frac{\partial}{\partial \phi}. $$
\end{Notation}

 Note that
$$-\triangle_sY_{lm}=l(l+1)Y_{lm}. $$\\


Next 
we have
\begin{equation}
\mathcal{K}g=\mathcal{K}(\check{g}Y_{lm})=H(r)Y_{lm}(\vartheta, \phi),
\end{equation}
where
\begin{equation}
H(r)=\mathcal{H}_l\check{g}(r)
\end{equation}
with
\begin{equation}
\mathcal{H}_l(\check{g})(r):=\frac{1}{2l+1}
\Big[\int_r^{+\infty}\check{g}(r')
\Big(\frac{r}{r'}\Big)^lr'dr'+
\int_0^r\check{g}_(r')
\Big(\frac{r}{r'}\Big)^{-l-1}r'dr'\Big]. \label{LetHlm}
\end{equation}
Here we assume that $\check{g}(r)=0$ for $r\geq R$. In fact $H=\mathcal{H}_l(\check{g})$ is determined by solving the ordinary differential equation
\begin{equation}
-\frac{1}{r^2}\frac{d}{dr}\Big(r^2\frac{dH}{dr}\Big)+\frac{l(l+1)}{r^2}H=
\check{g},
\end{equation}
for fixed $t$ and given $\check{g}$ such that $\check{g}(r)=0$ for $r \geq R$. Note that a fundamental system of the homogeneous equation is $r^l,r^{-l-1}$.

Then we have
\begin{equation}
\mathrm{grad}\mathcal{K}(\check{g}Y_{lm})=
\frac{d H}{dr}Y_{lm}\mathbf{e}_r+
\frac{H}{r}\nabla_sY_{lm},
\end{equation}
and the eigenvalue problem $\mathbf{L}\mbox{\boldmath$\xi$}=\lambda \mbox{\boldmath$\xi$}$ reduces to
\begin{subequations}
\begin{align}
 \frac{1}{r} \frac{d\check{G}}{dr} &= \lambda \psi \\
\frac{1}{r^2}\check{G}&=\lambda \chi \\
 0&=\lambda \kappa,
\end{align}
\end{subequations}
where
\begin{equation}
\check{G}=\check{G}(r)=-\frac{1}{\rho}\frac{dP}{d\rho}\check{g}+4\pi\mathsf{G}\mathcal{H}_{l}(\check{g}).
\end{equation}

Note that
\begin{align}
\frac{d}{dr}\mathcal{H}_l(\check{g})=&
\frac{1}{2l+1}\Big[l\int_r^{+\infty}\check{g}(r')\Big(\frac{r}{r'}\Big)^{l-1}dr' \nonumber \\
&-(l+1)\int_0^r\check{g}(r')\Big(\frac{r}{r'}\Big)^{-l-2}dr'\Big]. \label{Let18}
\end{align}

We see that $\check{g}, H,  \check{G}$ are independent of $\kappa$ and
$\kappa$ should vanish when $\lambda \not=0$. \\



Now, for the space of perturbations $\mbox{\boldmath$\xi$}$, we are considering  the Hilbert space $\mathfrak{H}$ of vector valued functions on 
$B_R$ endowed with the inner product
$$
(\mathbf{f}_1|\mathbf{f}_2)_{\mathfrak{H}}:=\int_{B_R}(\mathbf{f}_1(\vec{x})|
\mathbf{f}_2({x}))\rho(\vec{x})d{x}.
$$
Thus we see that for $\mbox{\boldmath$\xi$}$ of the form \eqref{xilm} we have
\begin{equation}
\|\mbox{\boldmath$\xi$}\|_{\mathfrak{H}}^2=
\|\psi\|_{\mathfrak{W}}^2+
l(l+1)\Big(
\|\chi\|_{\mathfrak{W}}^2+\|\kappa\|_{\mathfrak{W}}^2\Big),
\end{equation}
where we denote
\begin{equation}
\|f\|_{\mathfrak{W}}^2:=\int_0^R|f(r)|^2\rho(r)r^4dr \label{NormX}
\end{equation}
for a function $f : [0,R) \rightarrow \mathbb{C}$.
Here we have used the fact that $Y_{lm}$ is normalized so that 
$$\int_0^{2\pi}\int_0^{\pi}|Y_{lm}(\vartheta,\phi)|^2\sin\vartheta d\vartheta d\phi =1
$$
and, since $-\triangle_sY_{lm}=l(l+1)Y_{lm}$, this implies
\begin{align*}
&\int_0^{2\pi}\int_0^{\pi}\|\nabla_sY_{lm}(\vartheta,\phi)\|^2\sin\vartheta d\vartheta d\phi = \\
&=
\int_0^{2\pi}\int_0^{\pi}\|\nabla_s^{\perp}Y_{lm}(\vartheta,\phi)\|^2\sin\vartheta d\vartheta d\phi 
=l(l+1).
\end{align*}
See \cite{Jackson}. \\

Here let us introduce the following Hilbert spaces of scalar functions:

\begin{Definition}\label{Def.X}
1)
$\mathfrak{W}$ is the Hilbert space of scalar functions on $[0,R[$ endowed with the inner product
\begin{equation}
(f_1|f_2)_{\mathfrak{W}}=\int_0^R
f_1(r)f_2(r)^*\rho(r)r^4dr
\end{equation}
and the norm $\|\cdot\|_{\mathfrak{W}}$ defined by \eqref{NormX}.

2) For $\beta \in \mathbb{R}$, we put
\begin{equation}
\mathfrak{X}_{\beta}:=L^2([0,R], (R-r)^{\beta}r^2dr)
\end{equation}
endowed with the inner product
\begin{equation}
(f_1|f_2)_{\mathfrak{X}_{\beta}}=\int_0^Rf_1f_2^*(R-r)^{\beta}r^2dr.
\end{equation}
\end{Definition}

Of course we have
$$\|(rf)^{\flat}\|_{\mathfrak{H}}=\sqrt{4\pi}\|f\|_{\mathfrak{W}},
$$
where $(rf)^{\flat}$ stands for the function on $B_R$ defined by
$$
(rf)^{\flat}(\vec{x})=rf(r)\quad\mbox{with}\quad r=\|\vec{x}\|.$$
And 
$$ f \in \mathfrak{W} \quad \Leftrightarrow rf \in \mathfrak{X}_{\nu}.$$

\subsection{Case of $l=0$ with $m=0$}

Now let us consider $l=m=0$.

Then $Y_{lm}=Y_{00}=\frac{1}{\sqrt{4\pi}}$ and
\begin{equation}
\mbox{\boldmath$\xi$}=\frac{r\psi}{\sqrt{4\pi}}\mathbf{e}_r,\qquad
\check{g}=\frac{1}{r^2}\frac{d}{dr}(r^3\rho\psi).
\end{equation}
Now 
$\psi$ is a solution of 
\begin{equation}
\frac{1}{r}\frac{\partial}{\partial r}
\Big(-\frac{1}{\rho}\frac{dP}{d\rho}
\frac{1}{r^2}\frac{\partial}{\partial r}(r^3\rho \psi)
\Big)+\frac{4\pi\mathsf{G}}{r}\frac{d H_{0}(\check{g})}{d r}=\lambda \psi. \label{Eq00pre}
\end{equation}
But we see
$$\frac{1}{r}\frac{\partial H_{0}(\check{g})}{d r}=-\rho\psi.$$
Therefore \eqref{Eq00pre} reads
$$
\frac{1}{r}\frac{\partial}{\partial r}
\Big(-\frac{1}{\rho}\frac{dP}{d\rho}
\frac{1}{r^2}\frac{\partial}{\partial r}(r^3\rho \psi)
\Big)
-4\pi\mathsf{G}\rho\psi=\lambda \psi. 
$$
Moreover, keeping in mind that $-\triangle u=4\pi\mathsf{G}\rho$, we can rewrite this equation as
\begin{equation}
-\frac{1}{\rho r^4}\frac{\partial}{\partial r}
\Big(r^4\Gamma P\frac{\partial\psi}{\partial r}\Big)
-\Big(3\frac{\rho}{\Gamma}\frac{d\Gamma}{d\rho}+3\Gamma-4\Big)
\frac{1}{r}\frac{du}{dr}\psi=\lambda\psi \label{Eq00}
\end{equation}
through a tedious calculation. Here
\begin{equation}
\Gamma:=\frac{\rho}{P}\frac{dP}{d\rho}.
\end{equation}

Suppose $P=\mathsf{A}\rho^{\gamma}$ exactly. Then, since $\Gamma=\gamma$,
\eqref{Eq00} turns out to be
$$
-
\frac{1}{\rho r^4}\frac{\partial}{\partial r}\Big(
\gamma r^4 P\frac{\partial\psi}{\partial r}\Big)-(3\gamma-4)\frac{1}{r}\frac{du}{dr}\psi=
\lambda \psi.
$$
This is nothing but the equation for spherically symmetric perturbations. \\

 It is known that the differential operator $\mathcal{L}^{\mathsf{ss}}$:
\begin{equation}
\mathcal{L}^{\mathsf{ss}}\psi=-\frac{1}{\rho r^4}
\frac{d}{dr}\Big(r^4\Gamma P\frac{d\psi}{dr}\Big)
-\Big(3\frac{\rho}{\Gamma}\frac{d\Gamma}{d\rho}+3\Gamma-4\Big)\frac{1}{r}
\frac{du}{dr}\psi
\end{equation}
can be considered as a self-adjoint operator bounded from below with compact resolvents in
$\mathfrak{W}$. Here we note that
$$\frac{1}{C}\leq \Gamma \leq C,\quad
\Big|\frac{\rho}{\Gamma}\frac{d\gamma}{d\rho}\Big|\leq C
$$
on $0 < \rho \leq \rho_{\mathsf{O}}$. It is well known that if
$P=\mathsf{A}\rho^{\gamma}$ exactly and $\gamma >4/3$, then the least eigenvalue of $\mathcal{L}^{\mathsf{ss}}$ is positive. It is the case generally if 
$$3\frac{\rho}{\Gamma}\frac{d\Gamma}{d\rho}+3\Gamma-4\geq \frac{1}{C}>0$$
on $0<\rho\leq\rho_{\mathsf{O}}$.

\subsection{Case of $l\geq 1$ with $|m|\leq l $}

Next let us consider $l \geq 1$ with $|m|\leq l$. 

Introduce the variable
\begin{equation}
\vec{U}=
\begin{bmatrix}
\psi \\
\\
\sqrt{l(l+1)}\chi
\end{bmatrix}
\in \mathfrak{W}^2=\mathfrak{W}\times\mathfrak{W}.
\end{equation}
Then the eigenvalue problem reads
\begin{equation}
\vec{L}_l\vec{U}=\lambda \vec{U}, 
\end{equation}
where
\begin{equation}
\vec{L}_{l}\vec{U}:=
\begin{bmatrix}
\displaystyle\frac{1}{r}\frac{d}{dr}\check{G} \\
\\
\displaystyle\sqrt{l(l+1)}\frac{1}{r^2}\check{G}
\end{bmatrix}.
\end{equation}
Then, provided that $\vec{U},\vec{U}' \in C_0^{\infty}([0,R[)^2$, it can be verified by integration by parts that
\begin{align}
(\vec{L}_l\vec{U}|\vec{U}')_{\mathfrak{W}^2}&=
\int_0^R\frac{1}{\rho}\frac{dP}{d\rho}\check{g}(\check{g}')^*r^2dr-4\pi\mathsf{G}
\int_0^RH\cdot(\check{g}')^*r^2dr \nonumber \\
&=
\int_0^R\frac{1}{\rho}\frac{dP}{d\rho}\hat{g}(\check{g}')^*r^2dr+ \nonumber \\
&-{4\pi\mathsf{G}}
\int_0^{+\infty}\Big[r^2
\Big(\frac{dH}{dr}\Big)\Big(\frac{dH'}{dr}\Big)^*+
l(l+1)H\cdot (H')^*\Big]dr,
\end{align}
where
$$\check{g} 
=\frac{1}{r^2}\frac{\partial}{\partial r}
(r^3\rho\psi)-l(l+1)\rho\chi. \eqno\eqref{Letglm} ,
$$
and
$$
H=\mathcal{H}_l\check{g}(r)
$$
with
$$
\mathcal{H}_l\check{g}(r)=\frac{1}{2l+1}
\Big[\int_r^{+\infty}\check{g}(r')
\Big(\frac{r}{r'}\Big)^lr'dr'+
\int_0^r\check{g}(r')
\Big(\frac{r}{r'}\Big)^{-l-1}r'dr'\Big]. \eqno\eqref{LetHlm}
$$
Note that $H(r)=\mathcal{H}_l\check{g}(r)$ is defined for $\forall r \geq 0$ so that
$$H=O(r^{-l-1}),
\quad \frac{dH}{dr}=O(r^{-l-2})
$$
as $r \rightarrow +\infty$ provided that $\check{g} (r)=0$ for $r\geq R$.
Therefore we see that $\vec{L}_{l}$ as an operator with domain $C_0^{\infty}([0,R[)^2$ in the Hilbert space $\mathfrak{W}^2$ 
is symmetric. The boundedness from below of  $\vec{L}_{l}$ can be verified as follows:

Consider 
$$I=\int_0^RH\check{g}^*r^2dr.$$

Then we have
\begin{align*}
I&=\int_0^R\Big[H\frac{d}{dr}(r^3\rho\psi^*)-l(l+1)H\rho\chi^*r^2\Big]dr \\
&=-\Big[
\int_0^R\frac{dH}{dr}r^3\rho\psi^*dr+l(l+1)
\int_0^RH\rho\chi^*r^2dr\Big] \\
&\leq\frac{\epsilon}{2}\int\frac{1}{r^2}\Big|\frac{dH}{dr}\Big|^2\rho r^4dr  +
\frac{1}{2\epsilon}\int|\psi|^2\rho r^4dr + \\
&
+\frac{\epsilon}{2}\int l(l+1)\frac{|H|^2}{r^4}\rho r^4dr
+\frac{1}{2\epsilon}\int l(l+1)|\chi|^2\rho r^4dr 
\end{align*}
\begin{align*}
&=\frac{\epsilon}{2}\int\Big(r^2\Big|\frac{dH}{dr}\Big|^2+l(l+1)|H|^2\Big)\rho dr +
\frac{1}{2\epsilon}\|\vec{U}\|_{\mathfrak{W}^2}^2 \\
&\leq\frac{\epsilon\rho_{\mathsf{O}}}{2}
\int\Big(r^2\Big|\frac{dH}{dr}\Big|^2+l(l+1)|H|^2\Big) dr+
\frac{1}{2\epsilon}\|\vec{U}\|_{\mathfrak{W}^2}^2 \\
&=\frac{\epsilon\rho_{\mathsf{O}}}{2}\int H\Big(
-\frac{d}{dr}r^2\frac{dH^*}{dr}+l(l+1)H^*\Big)dr +
\frac{1}{2\epsilon}\|\vec{U}\|_{\mathfrak{W}^2}^2 \\
&=\frac{\epsilon\rho_{\mathsf{O}}}{2}\int H\check{g}^* r^2dr +
\frac{1}{2\epsilon}\|\vec{U}\|_{\mathfrak{W}^2}^2 \\
&=\frac{\epsilon\rho_{\mathsf{O}}}{2}I+
\frac{1}{2\epsilon}\|\vec{U}\|_{\mathfrak{W}^2}^2 .
\end{align*}
Thus
$$ I\leq \Big( 1-\frac{ \epsilon\rho_{\mathsf{O}} }{2} \Big)^{-1}\frac{1}{2\epsilon}\|\vec{U}\|_{\mathfrak{W}^2}^2
$$
provided that $\displaystyle \frac{\epsilon\rho_{\mathsf{O}}}{2} <1$. Taking, e.g., $\epsilon=1/\rho_{\mathsf{O}}$, we have
$$I\leq \rho_{\mathsf{O}}\|\vec{U}\|_{\mathfrak{W}^2}^2.$$

Therefore $\vec{L}_{l}$
can be considered as a self-adjoint operator in $\mathfrak{W}^2$ bounded from below. 

However the resolvents of $\vec{L}_{l}$ cannot be compact, and the spectrum of 
$\vec{L}_{l}$ is not of the Sturm-Liouville type, since
$$ \vec{\mathfrak{N}}_{l}:=\{ \vec{U}\  |\  \check{g}=0 \}
\subset \mbox{Ker}\vec{L}_{l} $$
and $\mbox{dim}\vec{\mathfrak{N}}_{l}=\infty$.

\subsection{Analysis of the operator $\mathcal{N}_l$ acting on
$\check{g}$}

Let us derive the equation for $\check{g}$.

Since
$$\check{g}=\frac{1}{r^2}\frac{d}{d r}(r^3\rho U^1)-\sqrt{l(l+1)}\rho U^2 $$
for
$$\vec{U}=
\begin{bmatrix}
U^1 \\
\\
U^2
\end{bmatrix}
=
\begin{bmatrix}
\psi \\
\\
\sqrt{l(l+1)}\chi
\end{bmatrix},
$$
putting
\begin{equation}
\mathcal{N}_l\check{g}=\frac{1}{r^2}\frac{\partial}{\partial r}(r^3\rho M^1)-
\sqrt{l(l+1)}\rho M^2
\end{equation}
for
$$\vec{M}_l\check{g}=
\begin{bmatrix}
M^1 \\
\\
M^2
\end{bmatrix}
=
\begin{bmatrix}
\displaystyle \frac{1}{r}\frac{d}{dr} \\
\\
\displaystyle \frac{\sqrt{(l(l+1)}}{r^2}
\end{bmatrix}
\Big(-\frac{1}{\rho}\frac{dP}{d\rho}\check{g}+4\pi\mathsf{G}\mathcal{H}_l\check{g}\Big),
$$
we get the equation 
\begin{equation}
\mathcal{N}_l\check{g}=\lambda \check{g}.
\end{equation}

We see
\begin{align}
&\mathcal{N}_l\check{g}=-\frac{\rho}{r^2}\frac{d\rho}{dP}\frac{d}{d r}\Big(\frac{r^2}{\rho}\Big(\frac{dP}{d\rho}\Big)^2\frac{d \check{g}}{d r}\Big)+q_l\check{g}+
4\pi\mathsf{G}\frac{d\rho}{dr}\frac{d}{d r} \mathcal{H}_l\check{g}, \\
&\mbox{with} \nonumber \\
&q_l:=-\frac{\rho}{r^2}\frac{d}{dr}\Big(\frac{r^2}{\rho}\frac{d}{dr}\Big(\frac{dP}{d\rho}\Big)\Big)+\Big[\frac{1}{r^2}\frac{d}{dr}\Big(\frac{r^2}{\rho}\frac{d\rho}{dr}\Big)
+\frac{l(l+1)}{r^2}\Big]\frac{dP}{d\rho}-4\pi\mathsf{G}\rho.
\end{align}

We can assume that
$$q_l(r) \geq -C $$
for $0<r<R$, although, for $l \geq 1$, 
$$q_l \sim \frac{l(l+1)}{r^2}K \quad\mbox{as}\quad r\rightarrow +0
$$
with $\displaystyle K=\frac{dP}{d\rho}\Big|_{r=+0} >0$. Therefore
$\mathcal{N}_l$ defined on 
$ C_0^{\infty}([0,R[)$ is symmetric and bounded from below in the Hilbert space
$$
\mathfrak{Y}=\{ \check{g} | \|\check{g}\|_{\mathfrak{Y}}=\sqrt{(\check{g}|\check{g})_{\mathfrak{Y}}}<\infty,
\int_0^R\check{g}(r)r^2dr =0 \}
$$
endowed with the inner product
$$(\check{g}_1|\check{g}_2)_{\mathfrak{Y}}=
\int_0^R\check{g}_1(\check{g}_2)^*\frac{1}{\rho}\frac{dP}{d\rho}r^2dr,
$$
and its Friedrichs extension is self-adjoint. 
Hereafter the same letter $\mathcal{N}_l$ will denote the Freidrichs extension.\\

Here we are introducing  the Hilbert spaces  $\mathfrak{Y}$ and $\mathfrak{Y}_1$ as following

\begin{Definition}
The Hilbert spaces  $\mathfrak{Y}$ and $\mathfrak{Y}_1$ are 
\begin{align}
&\mathfrak{Y}:=\{ \check{g} \  |\  \|\check{g}\|_{\mathfrak{Y}}=\sqrt{(\check{g}|\check{g})_{\mathfrak{Y}}}<\infty,
 \int_0^R \check{g}(r)r^2dr=0 \}, \\
&(\check{g}_1|\check{g}_2)_{\mathfrak{Y}}=\int_0^R\check{g}_1(\check{g}_2)^*\frac{1}{\rho}\frac{dP}{d\rho}r^2dr, \\
&\mathfrak{Y}_1=\{ \check{g}\in \mathfrak{Y} \  |\   
\|\check{g}\|_{\mathfrak{Y}_1}= \sqrt{(\check{g}|\check{g})_{\mathfrak{Y}_1}}<\infty \}, \\
&(\check{g}_1|\check{g}_2)_{\mathfrak{Y}_1}=
(\check{g}_1|\check{g}_2)_{\mathfrak{Y}}+
\int_0^R
\Big(\frac{d\check{g}_1}{dr}\Big)\Big(\frac{d\check{g}_2}{dr}\Big)^*\frac{1}{\rho}\Big(\frac{dP}{d\rho}\Big)^2r^2dr.
\end{align}
\end{Definition}

Then the unit ball of $\mathfrak{Y}_1$ in
$\mathfrak{Y}$ is precompact, the resolvents of the self-adjoint operator
$\mathcal{N}_l$ is compact, and the spectrum of $\mathcal{N}_l$ is of the Sturm-Liouville type. The proof 
can be done by the same way as in Section 3. \\

Let us take an eigenvalue $\lambda_j^{(\mathcal{N}_l)}$ of the operator
$\mathcal{N}_l$ and an associated eigenfunction $\varphi_j^{(\mathcal{N}_l)}(r)$.

Then, for any $m \in \mathbb{N}, |m|\leq l, $ we see that
$$g_{[lmj]}^{\mathcal{N}}:=\varphi_j^{(\mathcal{N}_l)}(r)
Y_{lm}(\vartheta, \phi) $$
is an eigenfunction of $\mathcal{N}$ associated with the eigenvalue 
$\lambda_j^{(\mathcal{N}_l)}$. 

Moreover, suppose that $\lambda_j^{(\mathcal{N}_l)}\not=0$. Then we can put
$$\vec{U}_{[lj]}:=\frac{1}{\lambda_j^{(\mathcal{N}_l)}}
\begin{bmatrix}
\frac{1}{r}\frac{d}{dr} \\
 \\
\frac{\sqrt{l(l+1)}}{r^2}
\end{bmatrix}
\mathcal{G}_l(\varphi_j^{(\mathcal{N}_l)}),
$$
where
$$\mathcal{G}_l(f)=-\frac{1}{\rho}\frac{dP}{d\rho}f+
4\pi\mathsf{G}\mathcal{H}_lf.
$$
Then we have
$$\vec{L}_{l}\vec{U}_{[lj]}=\lambda_j^{(\mathcal{N}_l)}\vec{U}_{[lj]}
$$
and $\vec{U}_{[lj]}\not=\vec{0}$, since
$$\varphi_j^{(\mathcal{N}_l)}=\frac{1}{r^2}\frac{d}{dr}
(r^3\rho U_{[lj]}^1)-\sqrt{l(l+1)}\rho U_{[lj]}^2\not=0,
$$ 
where $\vec{U}_{[lj]}=(U_{[lj]}^1, U_{[lj]}^2)^T$. In other words,
$\lambda_j^{(\mathcal{N}_l)}$ is an eigenvalue of
$\vec{L}_{l}$ and $\vec{U}_{[lj]}$ is an associated eigenfunction. 
Here we note that $m \in \mathbb{Z}, |m|\leq l$ is arbitrary. 

Finally put
$$
\mbox{\boldmath$\xi$}_{[lmj]}:=rU^1_{[lj]}(r)Y_{lm}(\vartheta,\phi)\mathbf{e}_r+
\frac{rU_{[lj]}^2(r)}{\sqrt{l(l+1)}}\nabla_s
Y_{lm}.
$$
Of course $\mbox{\boldmath$\xi$}_{[lmj]}\not=\mathbf{0}$. 

Then we see that $\lambda_j^{(\mathcal{N}_l)}$ is an eigenvalue of
$\mathbf{L}$ and $\mbox{\boldmath$\xi$}_{[lmj]}$ is an associated eigenfunction.

Note that 
$\lambda_j^{(\mathcal{N}_l)}$ is independent of $m$. 
Therefore the multiplicity of the eigenvalue is no smaller than $2l+1$, both for $\mathcal{N}$ and
$\mathbf{L}$. Summing up, we claim the following

\begin{Theorem}\label{Th.4}
For any $l \in \mathbb{N}$, the spectrum of the self-adjoint operator $\mathcal{N}_l$ is
of the Sturm-Liouville type. If $\lambda$ is an eigenvalue of some
$\mathcal{N}_l$, and if $\lambda \not=0$, then $\lambda$ is an eigenvalue of
$\mathbf{L}$ with multiplicity $\geq 2l+1$.
\end{Theorem}

Moreover, as for the multiplicity of the eigenvalues of the operator
$\mathcal{N}_l$, we have the following observation.\\

\begin{Theorem}\label{Th.5}
Let $l \geq 1$. 
Then eigenvalues of the operator $\mathcal{N}_l$ are simple.
\end{Theorem}

Proof. 
 Let us decompose the operator $\mathcal{N}_l$ as

\begin{align}
&\mathcal{N}_l=
\mathcal{N}_{l00}+\mathcal{N}_{l01}, \\
&\mathcal{N}_{l00}\check{g}:=
-\frac{\rho}{r^2}\frac{d\rho}{dP}\frac{d}{dr}\Big(\frac{r^2}{\rho}\Big(\frac{dP}{d\rho}\Big)^2
\frac{d \check{g}}{d r}\Big)+q_l\check{g}, \\
&\mathcal{N}_{l01}\check{g}:=
4\pi\mathsf{G}\frac{d\rho}{dr}\frac{d}{d r} \mathcal{H}_l\check{g}, 
\end{align}
with
\begin{align*}
&q_l:=-\frac{\rho}{r^2}\frac{d}{dr}\Big(\frac{r^2}{\rho}\frac{d}{dr}\Big(\frac{dP}{d\rho}\Big)\Big)+\Big[\frac{1}{r^2}\frac{d}{dr}\Big(\frac{r^2}{\rho}\frac{d\rho}{dr}\Big)
+\frac{l(l+1)}{r^2}\Big]\frac{dP}{d\rho}-4\pi\mathsf{G}\rho, \\
&\frac{d}{dr}\mathcal{H}_l\check{g}=\frac{1}{2l+1}\Big[l\int_r^R\check{g}(r')\Big(\frac{r}{r'}\Big)^{l-1}dr'
-(l+1)\int_0^r\check{g}(r')\Big(\frac{r}{r'}\Big)^{-l-2}dr\Big].
\end{align*}

Here
\begin{align*} 
\mathsf{D}(\mathcal{N}_{l00})&=\{ \check{g} \in \mathfrak{Y}_1 |
\mathcal{N}_{l00}\check{g} \in \mathfrak{Y}\quad\mbox{in distribution sense}\quad\} \\
&=\mathsf{D}(\mathcal{N}_l)
\end{align*}
and, of course,  $\mathcal{N}_{l00}$ is a self-adjoint operator with spectrum of Sturm-Liouville type, too. 

We are going to consider the multiplicities of eigenvalues of the operator
$\mathcal{N}_{l00}$.

In order to do it, we prepare the following
\begin{Lemma}\label{Lemma_LP}
For any $\lambda$, the equation
\begin{equation}
\mathcal{N}_{l00}\check{g}=\lambda \check{g} \label{LP.1}
\end{equation}
admits linearly independent set of solutions $\check{g}=\varphi_1(r), \varphi_2(r)$ such that
\begin{subequations}
\begin{align}
&\varphi_1(r)=r^l(1+O(r^2))\quad\mbox{as}\quad r \rightarrow +0, \label{LP.2a} \\
&\varphi_2(r)=r^{-l-1}(1+O(r^2))\quad\mbox{as}\quad r \rightarrow +0. \label{LP.2b}
\end{align}
\end{subequations}
\end{Lemma}

Proof. Let us use the Liouville transformation
$$x=\int_0^r\sqrt{\frac{d\rho}{dP}}dr,
\qquad
y=\Big(\frac{1}{\rho^2}\Big(\frac{dP}{d\rho}\Big)^3\Big)^{\frac{1}{4}}r\check{g}.
$$
Then the equation \eqref{LP.1} 
is transformed to the Liouville normal form
\begin{equation}
-\frac{d^2y}{dx^2}+\hat{q}y=\lambda y, \label{LP.3}
\end{equation}
where
\begin{align*}
\hat{q}=&q_l+\frac{1}{4}\frac{dP}{d\rho}\Big(
\frac{dA}{dr}-\frac{1}{4}A^2+AB\Big), \\
A:=&\frac{4}{r}+3\frac{d\rho}{dP}\frac{d}{dr}\Big(\frac{dP}{d\rho}\Big)-
\frac{2}{\rho}\frac{d\rho}{dr}, \qquad
B:=\frac{2}{r}+\frac{d\rho}{dP}\frac{d}{dr}\Big(\frac{dP}{d\rho}\Big)-
\frac{1}{\rho}\frac{d\rho}{dr}.
\end{align*}

See \cite[p.276, Theorem 6]{BirkhoffR}. Putting
$$x_+:=\int_0^R\sqrt{\frac{d\rho}{d P}}dr,
$$
we see that $x$ runs on the interval $[0, x_+]$ while
$r$ runs on $[0,R]$, and
$$x 
=\sqrt{\frac{d\rho}{dP}}\Big|_{\rho=\rho_{\mathsf{O}}}(1+O(r^2))r
\quad\mbox{as}\quad r\rightarrow +0,
$$
since $\rho=\rho_{\mathsf{O}}+O(r^2)$. 

We see
\begin{align*}
q_l&=\frac{l(l+1)}{r^2}\frac{dP}{d\rho}\Big|_{\rho=\rho_{\mathsf{O}}}(1+O(r^2)) \\
&=\frac{l(l+1)}{x^2}(1+O(x^2))
\end{align*}
as $r\rightarrow +0, x\rightarrow +0$, since
$$\rho=\rho_{\mathsf{O}}+O(r^2),\quad
\frac{d\rho}{dr}=O(r),\quad
\frac{d^2\rho}{dr^2}=O(1)$$
as $r\rightarrow +0$. 

Then it is easy to prove that \eqref{LP.3},
which is of the form
$$\frac{d^2y}{dx^2}=\frac{l(l+1)}{x^2}(1+\omega(x))y,\quad\mbox{where}\quad \omega(x)=O(x^2)
$$
 admits a set of independent solutions 
$y=\psi_1(x), \psi_2(x)$ such that
\begin{align*}
\psi_1(x)&=x^{l+1}(1+O(x^2))\quad\mbox{as}\quad x\rightarrow +0, \\
\psi_2(x)&=x^{-l}(1+O(x^2))\quad\mbox{as}\quad x\rightarrow +0.
\end{align*}

 In fact, putting
$$y=x^{\mu}Y,\quad \mu=\mu_{\pm}=\frac{1}{2}(1\pm(2l+1))=l+1, -l,
$$
the equation \eqref{LP.3} can be  transformed to the integral equation
\begin{align}
Y(x)&=1+ \nonumber \\
&\frac{1}{-2\mu+1}\Big[x^{-2\mu+1}\int_0^x
\Omega(x')Y(x')(x')^2dx'
-\int_0^x\Omega(x')Y(x')x'dx'\Big],\label{LP.Int}
\end{align}
where
$$\Omega(x):=l(l+1)\frac{\omega(x)}{x^2}=O(1),
$$
and
$$-2\mu+1=-2\mu_{\pm}+1=\pm(2l+1).
$$
It is easy to prove that \eqref{LP.Int} admits a solution $Y(x)=1+O(x^2)$
given by the iteration on a sufficiently short interval $[0,\delta]$ of $x$ with
$0<\delta \ll 1$. ---

Then the set of solutions of \eqref{LP.1}:
\begin{align*}
\varphi_1(r)=C
\Big(\frac{1}{\rho^2}\Big(\frac{dP}{d\rho}\Big)^3\Big)^{-\frac{1}{4}}\frac{1}{r}\psi_1(x), \\
\varphi_2(r)=C
\Big(\frac{1}{\rho^2}\Big(\frac{dP}{d\rho}\Big)^3\Big)^{-\frac{1}{4}}\frac{1}{r}\psi_2(x)
\end{align*}
with a suitable normalizing constant $C(>0)$ satisfies \eqref{LP.2a}, \eqref{LP.2b}. 
$\square$ \\

Let us suppose that an eigenvalue $\lambda_0$ of the eigenvalue problem \eqref{LP.1}
has two linearly independent associated eigenfunctions in
$\mathfrak{Y}$. Then it follows that any solution of \eqref{LP.1}
with $\lambda=\lambda_0$ would belong to $\mathfrak{Y}$. By Lemma \ref{Lemma_LP} it follows that $\varphi_2(r) \sim r^{-l-1} $ belongs to 
$\mathfrak{Y}$. But it is impossible, since, for $0<\delta\ll 1$, 
$$\int_0^{\delta}(r^{-l-1})^2\frac{P}{\rho^2}r^2dr
\gtrsim \int_0^{\delta}r^{-2l}dr =+\infty,
$$
since $-2l+1 \leq -1$ for $l\geq 1$. Therefore any eigenvalue of the operator $\mathcal{N}_{l00}$ is simple.

\begin{Remark}
The type of the boundary $x=x_+$ varies according to the value of $\gamma$.
Now we see 
$$
x_+-x \sim \frac{2}{\sqrt{(\gamma-1)K}}\sqrt{R-r}
\quad\mbox{as}\quad r\rightarrow R-0,
$$
where
$K$ is the positive constant such that $u \sim K(R-r)$ as $r\rightarrow R-0$, that is, 
$$K=-\frac{du}{dr}\Big|_{r=R}.$$

It is easy to see that 
$$\hat{q}\sim \frac{\kappa}{(x_+-x)^2}\quad\mbox{as}\quad
x\rightarrow x_+-0,
$$
provided that $\gamma \not=5/3$,
where
$$
\kappa:=\frac{(5-3\gamma)(3-\gamma)}{4(\gamma-1)^2} 
=\frac{3}{4}+\frac{3-2\gamma}{(\gamma-1)^2} 
=-\frac{1}{4}+\frac{(\gamma-2)^2}{(\gamma-1)^2}
$$
while  $\hat{q}=O(1)$ when $\gamma=5/3$. So, the boundary $x=x_+$ is of limit point type when $\gamma < \frac{3}{2}$, while it is of limit circle type when
$\frac{3}{2} < \gamma <2$. See \cite[p. 159, Theorem X.10]{ReedS}.

Since $\kappa >-\frac{1}{4}$, it is guaranteed that the operator
$-\frac{d^2}{dx^2}+\hat{q}(x)$ with domain $C_0^{\infty}(]0,x_+[)$ is bounded from below in $\L^2(0,x_+)$. See \cite[p. 340, Problem X.15]{ReedS}, \cite[p. 345, Footnote 1]{Kato}.
\end{Remark}

Let us go back to
$$\mathcal{N}_l=\mathcal{N}_{l00}+\mathcal{N}_{l01}.$$

Note that the perturbation $\mathcal{N}_{l01}$ is a bounded linear operator in 
$\mathfrak{Y}$. In fact we see
$$\Big\|\frac{d}{dr}\mathcal{H}_l\check{g}\Big\|_{L^{\infty}}
\leq \|\check{g}\|_{\mathfrak{Y}}\cdot
\Big[\int_0^R\rho\frac{d\rho}{dP}r^2dr\Big]^{1/2}
$$
and
$$\|\mathcal{N}_{l01}\check{g}\|_{\mathfrak{Y}}
\leq 
\Big\|\frac{d}{dr}\mathcal{H}_l\check{g}\Big\|_{L^{\infty}}
\cdot
4\pi\mathsf{G}
\Big[\int_0^R\Big(\frac{d\rho}{dr}\Big)^2\frac{1}{\rho}\frac{dP}{d\rho}r^2dr\Big]^{1/2}
$$
so that
$$
\|\mathcal{N}_{l01}\|_{\mathcal{B}(\mathfrak{Y})}
\leq4\pi\mathsf{G}
\cdot
\Big[\int_0^R\rho\frac{d\rho}{dP}r^2dr\Big]^{1/2}
\cdot
\Big[\int_0^R\Big(\frac{d\rho}{dr}\Big)^2\frac{1}{\rho}\frac{dP}{d\rho}r^2dr\Big]^{1/2}.
$$
Here note that
$$\rho\frac{d\rho}{dP}, \quad \Big(\frac{d\rho}{dr}\Big)^2\frac{1}{\rho}\frac{dP}{d\rho}
\quad \lesssim \quad (R-r)^{\frac{2-\gamma}{\gamma-1}}.
$$

Let us consider the one-parameter family of operators $T(t), 0\leq t\leq 1$, defined by
$$ T(t)=\mathcal{N}_{l00}+t\mathcal{N}_{l01} $$
in $\mathfrak{Y}$. Of course, $T(0)=\mathcal{N}_{l00}$ and
$T(1)=\mathcal{N}_l$. Moreover we know that, for $0\leq \forall t \leq 1$, $T(t)$ is a self-adjoint operator uniformly bounded from below with spectrum of the Sturm-Liouville type,
and we have shown that the eigenvalues of $T(0)$ are simple. We are going to prove that eigenvalues of $T(1)$ are simple.

Let us label eigenvalues of $T(t)$ as
$$\lambda_1(t)\leq \lambda_2(t)\leq \cdots \leq \lambda_j(t)\leq \lambda_{j+1}(t)\leq \cdots
\rightarrow +\infty.$$
Here we count the numbering taking into account the multiplicities. Recall that we know that the multiplicities are finite.

Let us fix an arbitrarily large integer $J$, and let us consider the proposition \\

$\mathcal{P}(t)$: {\it Eigenvalues $\lambda_1(t), \cdots, \lambda_J(t)$ are simple, that is, 
$$\lambda_1(t) <\lambda_2(t) <\cdots <
\lambda_{J-1}(t)<\lambda_J(t)<\lambda_{J+1}(t). $$}\\

Put
$$A:=\{ t\in [0,1] \quad |\quad  \mathcal{P}(t)\quad\mbox{holds valid.} \}.$$

We know that $0\in A$. 

If $t\in A$, then 
$$T(t+\Delta t)=T(t)+\Delta t\cdot\mathcal{N}_{l01} $$
satisfies
$$\hat{\delta}(T(t+\Delta t), T(t))\leq |\Delta t|\cdot
|\|\mathcal{N}_{l01}\||_{\mathcal{B}(\mathfrak{Y})}.$$
Here $\hat{\delta}(T(t+\Delta t), T(t))$ means the gap between the graphs of the closed operators 
$T(t+\Delta t)$ and $T(t)$ in $\mathfrak{Y}\times\mathfrak{Y}$.
See \cite[p.203, Theorem 2.14]{Kato}. We are supposing
$$\lambda_1(t)<\lambda_2(t)<\cdots <\lambda_{J-1}(t)<\lambda_J(t)<\lambda_{J+1}(t). $$
Applying the stability theorem \cite[p.212, Theorem 3.16]{Kato}
to these finite number of eigenvalues, we can claim that if $|\Delta t|$ is sufficiently small, then for $1\leq \forall j \leq J$, $\lambda_j(t+\Delta t)$ is near to
$\lambda_j(t)$ and simple. That is, $t+\Delta t \in A$, provided that $|\Delta t|$ is sufficiently small. Thus, $A$ is open in $[0,1]$. Inversely, let $t\in ]0,1]$  belong to the compliment of $A$, that is $t\not\in A$. Then there would exit an integer $j \leq J$ and an integer $m \geq 1$ such that
$$j\geq 2 \quad\mbox{and}\quad
\lambda_{j-1}(t)<\lambda_j(t)=\lambda_{j+1}(t)=
\cdots =\lambda_{j+m}(t)<\lambda_{j+m+1}(t) $$
or
$$j=1\quad\mbox{and}\quad
\lambda_1(t)=\lambda_2(t)=\cdots =
\lambda_{m+1}(t)<\lambda_{m+2}(t)
, $$
that is,
$$\mbox{dim}\mbox{Ker}(\lambda_j(t)-T(t))=m.$$
Applying the stability theorem, we can claim that $T(t+\Delta t )$ has the same property:
$$\mbox{dim}\mbox{Ker}(\lambda_j(t+\Delta t)-T(t+\Delta t))=m \geq 1, $$
provided that $|\Delta t|$ is sufficiently small. Then $t +\Delta t \not\in A$.
Thus $[0,1]\setminus A$ is open. 

Therefore we see $A=[0,1]$. Since $J$ was arbitrary, this completes the proof.
$\square$

\begin{Remark}
It is difficult to deduce from Theorem \ref{Th.4} and Theorem \ref{Th.5} the exact information on the multiplicities of the non-zero eigenvalues of the operator $\mathbf{L}$, since it is not clear whether eigenvalues which comes from different $l$'s
can coincide or not. 
\end{Remark}

\section{Regularity of $g$ and its expansion with respect to spherical harmonics }

Let us consider a vector field
$\mbox{\boldmath$\xi$}$ which belongs to the domain
\begin{align*}
\mathsf{D}(\mathbf{L})=&
\{ \mbox{\boldmath$\xi$} \in \mathfrak{F} \  |\  \mathrm{div}(\rho\mbox{\boldmath$\xi$})\in\mathsf{D}(\mathcal{N}) \} \\
=&
\{ \mbox{\boldmath$\xi$} \in \mathfrak{H} \  |\  \mathrm{div}(\rho\mbox{\boldmath$\xi$})\in\mathsf{D}(\mathcal{N}) \}.
\end{align*}
Then we can put
$$\Psi=-\mathcal{K}g,\qquad g=\mathrm{div}(\rho\mbox{\boldmath$\xi$}). $$
Since $g \in \mathfrak{G}$,  we easily see that   $\Psi \in C(\mathbb{R}^3)$, using the Schwartz inequality. Moreover $\mathrm{grad}\Psi \in L^2(B_R)$. ( See 
the discussion on $\mathcal{N}_{01}$ in Section 3.) 
Note that $\triangle \Psi=g $ in distribution sense. \\

In this Section,  we consider $g, G=-\frac{1}{\rho}\frac{dP}{d\rho}g +4\pi\mathsf{G}\mathcal{K}g, \Psi$ using their series expansion by the spherical harmonics $\{ Y_{lm}  |  l, m \in \mathbb{Z}, l \geq 0, |m| \leq l \}$.


We shall use the notations introduced before the  proof of Proposition \ref{Prop3}. And we denote 
$\nabla\cdot\mathbf{A}, \nabla f$
for 
$\mathrm{div}\mathbf{A},  \mathrm{grad}f $
for any vector field $\mathbf{A}$ and scalar field $f$.


\subsection{Degenerate elliptic problem}

Considering the functional space 
$\displaystyle W^{1,2}(B_R, \rho\frac{d\rho}{dP},\rho)=
W^{1,2}(B_R, \frac{d\rho}{du},\rho)=
W^{1,2}(\mathsf{d}^{\nu-1},\mathsf{d}^{\nu})
$, we introduce
the  following :

\begin{Definition}

1) The set of all $\displaystyle  U \in W^{1,2}(B_R, \rho\frac{d\rho}{dP}, \rho)$ such that
\begin{equation}
(U|1)_{L^2(\rho d\rho/dP)}=\int U(\vec{x})\frac{d\rho}{du}d\vec{x}=0
\end{equation}
will be denoted by $\mathfrak{E}$,
which is a closed subspace of $\displaystyle W^{1,2}(B_R, \rho\frac{d\rho}{dP},\rho)$
endowed with the same norm.

2) For any $\displaystyle U\in W^{1,2}(B_R, \rho\frac{d\rho}{dP},\rho)$ we put
\begin{equation}
U^{[Z]}:=U- \|1\|_{L^2(\rho d\rho/dP)}^{-2}(U|1)_{L^2(\rho d\rho/dP)},
\end{equation}
which clearly belongs to $\mathfrak{E}$.
\end{Definition}

Note that
$$\|1\|_{L^2(\rho d\rho/dP)}=\int_{B_R}\rho\frac{d\rho}{dP}r^2dr=\int_{B_R}\frac{d\rho}{du}r^2dr,$$
and
$$(U|1)_{L^2(\rho d\rho/dP)}=\int_{B_R}U\rho\frac{d\rho}{dP}r^2dr=\int_{B_R}U\frac{d\rho}{du}r^2dr,$$
so that
$$\|1\|_{L^2(\rho d\rho/dP)}^{-2}(U|1)_{L^2(\rho d\rho/dP)}=
\frac{\int_{B_R}U\frac{d\rho}{du}r^2dr}{\int_{B_R}\frac{d\rho}{du}r^2dr}.$$

\

\begin{Lemma}\label{LemmaPo}(Poincar\'{e} inequality)
If $ U\in \mathfrak{E}$,
then 
\begin{equation}
\|U\|_{L^2(\mathsf{d}^{\nu-1})}\leq C\|\nabla U\|_{L^2(\mathsf{d}^{\nu})}. \label{PO}
\end{equation}
\end{Lemma}

Proof. Suppose \eqref{PO} is false for every $C$. Then there exists a sequence of functions $U_j \in \mathfrak{E}$ such that $\|U_j\|_{L^2(\mathsf{d}^{\nu-1})}=1$ and 
$\|\nabla U_j\|_{L^2(\mathsf{d}^{\nu})} \rightarrow 0$ as $j \rightarrow \infty$. But the imbedding 
of $W^{1,2}(B_R, \mathsf{d}^{\nu-1}, \mathsf{d}^{\nu})$ into $L^2(B_R,\mathsf{d}^{\nu-1})$
is compact, according to \cite[8.8 Theorem]{GurkaO.2}. (Recall $\nu >1$ and 
$W^{1,2}(B_R, \mathsf{d}^{\nu-1},\mathsf{d}^{\nu}) \hookrightarrow 
W^{1,2}(B_R, \mathsf{d}^{\nu}, \mathsf{d}^{\nu})$.) 

Therefore a subsequence
$U_{j(k)}$ converges to $U \in W^{1,2}(B_R,\mathsf{d}^{\nu-1},\mathsf{d}^{\nu})$ weakly and strongly in $L^2(B_R, \mathsf{d}^{\nu-1})$. Since 
$\|\nabla U_j\|_{L^2(\mathsf{d}^{\nu})} \rightarrow 0$,
we have $\|\nabla U\|_{L^2(\mathsf{d}^{\nu})}=0$. Since $B_R$ is connected, $U =\mbox{Const.}$ a. e. 
But, since $(U_j|1)_{L^2(\rho d\rho/dP)}=0$ for $\forall j$,
we have $(U|1)_{L^2(\rho d\rho/dP)}=0$. Hence $U=0$. On the other hand, since
$\|U_j\|_{L^2(\mathsf{d}^{\nu-1})}=1$ for $\forall j$, we have $\|U\|_{L^2(\mathsf{d}^{\nu-1})}=1$
because of the strong convergence. A contradiction.
$\square$\\

Thanks to this  Lemma \ref{LemmaPo}, $\|U\|_{\mathfrak{E}}:=\|\nabla U\|_{L^2(\rho)} $ turns out to be an equivalent norm of 
$\mathfrak{E}$. 

As a result of Lemma \ref{LemmaPo}, the Lax-Milgram theorem reads

\begin{Lemma}\label{LemmaLaxM}  1) For any $f \in L^2(B_R, \frac{1}{\rho}\frac{dP}{d\rho})=L^2(B_R,\mathsf{d}^{-\nu+1})$, there is a unique 
 $U$ in
$\mathfrak{E}$ which satisfies
\begin{equation}
(\nabla U|\nabla V)_{L^2(\rho)}=-\int_{B_R}f V^*d\vec{x}
\end{equation}
for any $V \in \mathfrak{E}$. Moreover
\begin{equation}
\|U\|_{\mathfrak{E}}\leq C\|f\|_{L^2(\mathsf{d}^{-\nu+1})}.
\end{equation}

2) If $f \in \mathfrak{G}$, then the solution $U$ is a weak solution of  the equation
\begin{equation}
\nabla\cdot(\rho\nabla U)=f. \label{7.6}
\end{equation}
\end{Lemma}

Proof.  1) In fact
$$V \mapsto \int_{B_R}fV^* $$
is a continuous conjugate-linear functional on
 $\mathfrak{E}$, provided that $f \in L^2(B_R, \mathsf{d}^{-\nu+1})$,
and
$a(U,V)=(\nabla U|\nabla V)_{L^2(\rho)}$ is the inner product of 
$\mathfrak{E}$ thanks to Lemma \ref{LemmaPo}.

2) In general the field $U$ is a weak solution of the equation
$$ \nabla\cdot(\rho\nabla U)=f -\|1\|_{L^2(d\rho/du)}^{-2}\frac{d\rho}{du}\int f. $$
This is verified, since $V=\varphi^{[Z]} \in \mathfrak{E}$ for any $\varphi \in
C_0^{\infty}(B_R)$. If $f \in \mathfrak{G}$, then $\int f=0$ and the statement \eqref{7.6} follows.
 $\square$\\

We shall use the following

\begin{Lemma}\label{LemmaReg}(Elliptic regularity) Let
$U \in \mathfrak{E}$ with
$f=\nabla\cdot(\rho\nabla U) \in \mathfrak{G}\subset
L^2(B_R, \mathsf{d}^{-\nu+1})$. Then $U \in W_{\text{loc}}^{2,2}(B_R)$ and
\begin{equation}
\Big\|\frac{\partial^2U}{\partial x^j\partial x^k}\Big\|_{L^2(\mathsf{d}^{\nu+1})}
\leq C\|f\|_{L^2(\mathsf{d}^{-\nu+1})}. \label{7.11}
\end{equation}
\end{Lemma}

Proof. By the theory of interior elliptic regularity, we know
$U \in W_{\text{loc}}^{2,2}(B_R)$ and
\begin{equation}
\rho\triangle U+
(\nabla\rho|\nabla U)=f \label{R1}
\end{equation}
holds a. e. in $B_R$. A proof can be found in \cite[Theorem 8.8]{GilbargT}.

1) First let us show $\triangle U \in L^2(B_R,\mathsf{d}^{\nu+1})$.

In order to do it, we introduce the function
$w=\rho^{1/\nu} \in C^{\infty}(B_R)$. 

Of course
$$\frac{1}{C} u \leq w \leq C u, $$
and
$$\frac{1}{C}w \leq \mathsf{d} \leq C w.$$
Moreover, since we have 
$$w=\frac{\gamma-1}{\mathsf{A}\gamma}u(1+\Lambda_w(u))\quad\mbox{for}\quad u >0
$$
with some $\Lambda_w\in C^{\infty}(\mathbb{R})$, we can consider $w$ is a smooth function of
$u \in \mathbb{R}$ such that
$$w \lessgtr 0 \quad
\Leftrightarrow \quad u \lessgtr 0.
$$

Supposing that $U$ is smooth, \eqref{R1} implies 
\begin{align*}
f^2w^{-\nu+1}&=(\triangle U)^2w^{\nu+1}+
2\nu(\triangle U)(\nabla U|\nabla w)w^{\nu} \\
&+\nu^2(\nabla U|\nabla w)^2w^{\nu-1}.
\end{align*}
Consider
\begin{align*} 
Q:=&(\triangle U)(\nabla U|\nabla w)w^{\nu} \\
&=\sum_{i,j}(\partial_j\partial_jU)(\partial_iU)(\partial_iw)w^{\nu},
\end{align*}
where $\partial_k=\partial/\partial x^k, k=1,2,3$. 

We have
\begin{align*}
Q&=\sum\partial_j\Big[(\partial_jU)(\partial_iU)(\partial_iw)w^{\nu}\Big] \\
&-\sum(\partial_jU)(\partial_iU)(\partial_j\partial_i w)w^{\nu}
-\nu\sum(\partial_jU)(\partial_iU)(\partial_jw)(\partial_iw)w^{\nu-1} \\
&-\sum(\partial_jU)(\partial_j\partial_iU)(\partial_iw)w^{\nu} \\
&=\sum\partial_j\Big[(\partial_jU)(\partial_iU)(\partial_iw)w^{\nu}\Big] \\
&-\sum(\partial_jU)(\partial_iU)(\partial_j\partial_i w)w^{\nu}
-\nu\sum(\partial_jU)(\partial_iU)(\partial_jw)(\partial_iw)w^{\nu-1} \\
&-\frac{1}{2}\sum\partial_i((\partial_jU)^2)(\partial_iw)w^{\nu} \\
&=\sum\partial_j\Big[(\partial_jU)(\partial_iU)(\partial_iw)w^{\nu}\Big] 
-\frac{1}{2}\sum \partial_i\Big[(\partial_jU)^2(\partial_iw)w^{\nu}\Big]\\
&-\sum(\partial_jU)(\partial_iU)(\partial_j\partial_i w)w^{\nu}
-\nu\sum(\partial_jU)(\partial_iU)(\partial_jw)(\partial_iw)w^{\nu-1} \\
&+\frac{1}{2}\sum(\partial_jU)^2(\partial_i\partial_iw)w^{\nu}
+\frac{1}{2}\sum(\partial_jU)^2(\partial_iw)^2w^{\nu-1} \\
&=\sum\partial_j\Big[\cdots (\partial_jw)w^{\nu}\Big] 
-\frac{1}{2}\sum\partial_i\Big[\cdots (\partial_iw)w^{\nu}\Big]\\
&-(\nabla U|D^2w.\nabla U)w^{\nu} -\nu(\nabla U|\nabla w)^2w^{\nu-1} \\
&+\frac{1}{2}\|\nabla U\|^2(\triangle w)w^{\nu}
+\frac{\nu}{2}\|\nabla U\|^2\|\nabla w\|^2w^{\nu-1}.
\end{align*}

Therefore we have
\begin{align*}
f^2w^{-\nu+1}&=(\triangle U)^2w^{\nu+1}+
\nu^2\Big(\|\nabla U\|^2\|\nabla w\|^2-
(\nabla U|\nabla w)^2\Big)w^{\nu-1} + \\
&+\nu\|\nabla U\|^2(\triangle w)w^{\nu}-
2\nu (\nabla U|D^2w.\nabla U) w^{\nu} +\\
&+2\nu\sum \partial_j\Big[\cdots (\partial_jw)w^{\nu}\Big].
\end{align*}
Note that 
$$\|\nabla U\|^2\|\nabla w\|^2 -(\nabla U|\nabla w)^2 \geq 0.$$

Supposing that $U \in C^{\infty}(\overline{B_R})$, the integration of the above relation implies
$$\|f\|_{L^2(\mathsf{d}^{-\nu+1})}^2
\gtrsim \|\triangle U\|_{L^2(\mathsf{d}^{\nu+1})}^2 - C
\|\nabla U\|_{L^2(\mathsf{d}^{\nu})}^2,
$$
since $w \in C^{2,\alpha}(\overline{B_R})$ so that $D^2w =O(1)$. Thus
\begin{align*}
\|\triangle U\|_{L^2(\mathsf{d}^{\nu+1})}^2&\lesssim
\|f\|_{L^2(\mathsf{d}^{-\nu+1})}^2+
\|\nabla U\|_{
L^2(\mathsf{d}^{\nu})}^2 \\
&\lesssim \|f\|_{
L^2(\mathsf{d}^{-\nu+1})}^2,
\end{align*}
or
\begin{equation}
\|\triangle U\|_{L^2(\mathsf{d}^{\nu+1})} \leq C \|f\|_{L^2(\mathsf{d}^{-\nu+1})}.\label{7.13}
\end{equation}

Using this estimate of $\triangle U$, we show
\begin{equation}
\sum_{j,k}\int_{B_R}(\partial_j\partial_kU)^2w^{\nu+1} \leq C \|f\|_{L^2(\mathsf{d}^{-\nu+1})}^2,
\end{equation}
which says \eqref{7.11}.

Supposing that $U$ is smooth, we have
\begin{align*}
\sum (\partial_j\partial_kU)^2w^{\nu+1}&=
\sum (\partial_j\partial_kU)(\partial_j\partial_kU)w^{\nu+1} \\
&=\sum\partial_j\Big[(\partial_kU)(\partial_j\partial_kU)w^{\nu+1}\Big] \\
&-\sum(\partial_kU)(\partial_j^2\partial_kU)w^{\nu+1}
-(\nu+1)\sum
(\partial_kU)(\partial_j\partial_kU)(\partial_jw)w^{\nu} \\
&=\sum\partial_j\Big[(\partial_kU)(\partial_j\partial_kU)w^{\nu+1}\Big]
-\sum\partial_k\Big[(\partial_kU)(\partial_j^2U)w^{\nu+1}\Big] \\
&+\sum (\partial_k^2U)(\partial_j^2U)w^{\nu+1}+
(\nu+1)\sum(\partial_kU)(\partial_j^2U)(\partial_kw)w^{\nu} \\
&-(\nu+1)\sum(\partial_kU)(\partial_j\partial_kU)(\partial_jw)w^{\nu} \\
&=\sum\partial_j\Big[(\partial_kU)(\partial_j\partial_kU)w^{\nu+1}\Big]
-\sum\partial_k\Big[(\partial_kU)(\partial_j^2U)w^{\nu+1}\Big] \\
&+(\triangle U)^2w^{\nu+1}+(\nu+1)(\triangle U)(\nabla U|\nabla w)w^{\nu} \\
&-(\nu+1)\sum(\partial_kU)(\partial_j\partial_kU)(\partial_jw)w^{\nu} \\
&=\sum\partial_j\Big[(\partial_kU)(\partial_j\partial_kU)w^{\nu+1}\Big]
-\sum\partial_k\Big[(\partial_kU)(\partial_j^2U)w^{\nu+1}\Big] \\
&-(\nu+1)\sum\partial_k\Big[(\partial_kU)(\nabla U|\nabla w)w^{\nu}\Big] \\
&+(\triangle U)^2w^{\nu+1}+(\nu+1)(\triangle U)(\nabla U|\nabla w)w^{\nu} \\
&+(\nu+1)\nabla\cdot(w^{\nu}\nabla U)(\nabla U|\nabla w)
+(\nu+1)(\nabla U|D^2w.\nabla U)w^{\nu} \\
&=\sum\partial_j\Big[(\partial_kU)(\partial_j\partial_kU)w^{\nu+1}\Big]
-\sum\partial_k\Big[(\partial_kU)(\partial_j^2U)w^{\nu+1}\Big] \\
&-(\nu+1)\sum\partial_k\Big[(\partial_kU)(\nabla U|\nabla w)w^{\nu}\Big] \\
&+(\triangle U)^2w^{\nu+1}+(\nu+1)(\triangle U)(\nabla U|\nabla w)w^{\nu} \\
&+(\nu+1)f(\nabla U|\nabla w)
+(\nu+1)(\nabla U|D^2w.\nabla U)w^{\nu} .
\end{align*}

Supposing that $U \in C^{\infty}(\overline{B_R})$, the integration of the above relation implies
\begin{align*}
\int \sum (\partial_j\partial_kU)^2w^{\nu}
&=
\int (\triangle U)^2w^{\nu+1}+(\nu+1)\int(\triangle U)(\nabla U|\nabla w)w^{\nu} \\
&+(\nu+1)\int f(\nabla U|\nabla w)
+(\nu+1)\int(\nabla U|D^2w.\nabla U)w^{\nu}.
\end{align*}
But we have $(\nabla U|\nabla w) \in L^2(\mathsf{d}^{\nu-1})$ and
\begin{align*}
\|(\nabla U|\nabla w)\|_{L^2(\mathsf{d}^{\nu-1})}&\lesssim \|f\|_{L^2(\mathsf{d}^{-\nu+1})}
+\|\triangle U\|_{L^2(\mathsf{d}^{\nu+1})} \\
&\lesssim \|f\|_{L^2(\mathsf{d}^{-\nu+1})},
\end{align*}
since
$$(\nabla U|\nabla w)=\frac{1}{\nu}\Big[w^{1-\nu}f -w\triangle U \Big].
$$
Therefore
\begin{align*}
\int \sum (\partial_j\partial_kU)^2w^{\nu}&
\lesssim \|\triangle U\|_{L^2(\mathsf{d}^{\nu+1})}^2+
\|f\|_{L^2(\mathsf{d}^{-\nu+1})}^2+
\|\nabla U\|_{L^2(\mathsf{d}^{\nu})}^2 \\
&\lesssim \|f\|_{L^2(\mathsf{d}^{-\nu+1})}^2.
\end{align*}

This completes the proof by using the following Proposition. $\square$

\

\begin{Proposition}
Let $U\in \mathfrak{E}$ and $\nabla\cdot(\rho\nabla U)=f \in \mathfrak{G}$.
Then there exists a sequence $U_n\in C^{\infty}(\overline{B_R})$ such that 
$(U_n|1)_{L^2(\rho d\rho/dP)}=0 ,  \quad U_n\rightarrow U \quad\mbox{in}\quad
(\mathfrak{E}, \|\cdot\|_{\mathfrak{E}})$ and
$f_n:=\nabla\cdot(\rho\nabla U_n) \rightarrow f $ in $\mathfrak{G}\hookrightarrow L^2(\mathsf{d}^{-\nu+1})$.
\end{Proposition}

Proof. $C^{\infty}(\overline{B_R})$ is dense in
 $W^{1,2}(B_R, \mathsf{d}^{\nu-1},\mathsf{d}^{\nu-1})$ according to \cite[7.2 Theorem]{Kufner} applied for $\varepsilon=\nu-1>0$ with $p=2$. Hence $C^{\infty}(\overline{B_R})$ is dense in $W^{1,2}(B_R,\mathsf{d}^{\nu-1},\mathsf{d}^{\nu})$. Let $f_n \in C^{\infty}(\overline{B_R}) \rightarrow
U$ in $W^{1,2}(B_R,\mathsf{d}^{\nu-1}, \mathsf{d}^{\nu})$. Put
$U_n:=f_n^{[Z]}=f_n-(\|1\|_{L^2(\rho d\rho/dP)})^{-2}(f_n|1)_{L^2(\rho d\rho/dP)}$.
Then $U_n \in C^{\infty}(\overline{B_R})$ and $(U_n|1)_{L^2(\rho d\rho/dP)}=0$ and $U_n \rightarrow U$ in $(\mathfrak{E}, \|\cdot\|_{\mathfrak{E}})$. On the other hand, $f_n:=
\nabla\cdot(\rho\nabla U_n)$ satisfies 
$$ (\nabla U_n|\nabla V)=-\int f_n V^*
$$
for $\forall V \in C^{\infty}(\overline{B_R})$. Therefore
\begin{align*}
\|f_n\|_{L^2(\mathsf{d}^{-\nu+1})}&=
\sup\{ \frac{|\int f_nV^*|}{ \|V\|_{L^2(\mathsf{d}^{\nu-1})} }\quad|\quad V \in
L^2(B_R, \mathsf{d}^{\nu-1}) \} \\
&=
\sup\{ \frac{|\int f_nV^*|}{ \|V\|_{L^2(\mathsf{d}^{\nu-1})} }\quad|\quad V \in
\mathfrak{E} \} \\
&=
\sup\{ \frac{ |(\nabla U_n|\nabla V)_{L^2(\rho)}| }{ \|V\|_{L^2(\mathsf{d}^{\nu-1})} }\quad|
\quad V \in C^{\infty}(\overline{B_R}) \cap \mathfrak{E} \} \\
&\rightarrow \\
&
\sup\{ \frac{ |(\nabla U|\nabla V)_{L^2(\rho)}| }{ \|V\|_{L^2(\mathsf{d}^{\nu-1})} }\quad|
\quad V \in C^{\infty}(\overline{B_R}) \cap \mathfrak{E} \} \\
&=\|f\|_{L^2(\mathsf{d}^{-\nu+1})}.
\end{align*}
$\square$\\
 
\subsection{Spherical harmonics expansion of $g$}

We next claim the following:

\begin{Lemma}\label{Lemma.glm}
Let $\mbox{\boldmath$\xi$} \in \mathsf{D}(\mathbf{L})$ so that $g=\mathrm{div}(\rho\mbox{\boldmath$\xi$})\in \mathsf{D}(\mathcal{N})$. Then 

i)  $g \in C(B_R), \Psi:=-\mathcal{K}g \in C(\mathbb{R}^3), 
\displaystyle G=-\frac{1}{\rho}\frac{dP}{d\rho} g+4\pi\mathsf{G}\mathcal{K}g \in C(B_R)$ and
\begin{align}
g_{lm}(r)&=\int_0^{2\pi}\int_0^{\pi}
g(r\mbox{\boldmath$\omega$})Y_{lm}(\vartheta,\phi)^*\sin\vartheta d\vartheta d\phi, \\
G_{lm}(r)&=\int_0^{2\pi}\int_0^{\pi}
G(r\mbox{\boldmath$\omega$})Y_{lm}(\vartheta,\phi)^*\sin\vartheta d\vartheta d\phi, \\
\Psi_{lm}(r)&=\int_0^{2\pi}\int_0^{\pi}
\Psi(r\mbox{\boldmath$\omega$})Y_{lm}(\vartheta,\phi)^*\sin\vartheta d\vartheta d\phi \\
&=-\int_0^{2\pi}\int_0^{\pi}
(\mathcal{K}g)(r\mbox{\boldmath$\omega$})Y_{lm}(\vartheta,\phi)^*\sin\vartheta d\vartheta d\phi
\end{align}
are well-defined for $0<\forall r < R$, where
$$\mbox{\boldmath$\omega$}=\sin\vartheta\cos\phi\frac{\partial}{\partial x^1}
+\sin\vartheta\cos\phi\frac{\partial}{\partial x^2}+
\cos\phi\frac{\partial}{\partial x^3},
$$
while $g_{lm}, G_{lm} \in C([0,R[)\cap L^2([0,R]), \Psi_{lm}\in C([0,R])$;

ii) As $L \rightarrow \infty$,
$$g^{(L)}:=\sum_{|l|\leq L}g_{lm}Y_{lm}$$
converges to $g$ in $L^2(B_R)$, and converges uniformly on $\overline{B_{(1-\delta)R}}$
for each fixed $\delta \in ]0,1[$;

iii)
It holds that
\begin{equation}
-\frac{1}{\rho}\frac{dP}{d\rho}g_{lm}+4\pi\mathsf{G}\mathcal{H}_l(g_{lm})=G_{lm},
\end{equation}
where 
$$\Psi_{lm}=-\mathcal{H}_l(g_{lm})$$
and
$$
\mathcal{H}_l\check{g}(r)=\frac{1}{2l+1}
\Big[\int_r^{+\infty}\check{g}(r')
\Big(\frac{r}{r'}\Big)^lr'dr'+
\int_0^r\check{g}(r')
\Big(\frac{r}{r'}\Big)^{-l-1}r'dr'\Big], \eqno\eqref{LetHlm}
$$
which satisfies
$$
-\frac{1}{r^2}\frac{d}{dr}\Big(r^2\frac{d}{dr}\mathcal{H}_l\check{g}\Big)+\frac{l(l+1)}{r^2}
\mathcal{H}_l\check{g}=
\check{g}
$$
in distribution sense for $\check{g} \in L^2([0,R])$.
\end{Lemma}

Proof. Since $g \in \mathfrak{G} \hookrightarrow L^2(B_R)$, we see $\Psi=-\mathcal{K}g \in C(\mathbb{R}^3)$ and
$\|\nabla \Psi\|_{L^2(B_R)}\lesssim \|g\|_{L^2(B_R)}$. See the discussion on 
$\mathcal{N}_{01}$ in Section 3. Therefore we have
$$G=-\frac{1}{\rho}\frac{dP}{d\rho}g+4\pi\mathsf{G}\mathcal{K}g \in
 W^{1,2}(B_R, \mathsf{d}^{\nu-1},\mathsf{d}^{\nu}), $$
since
$\displaystyle \frac{1}{\rho}\frac{dP}{d\rho} \in C^{\infty}(B_R)$ and $\lesssim 
\mathsf{d}^{-\nu+1}$, while $g \in \mathfrak{G}_1 \hookrightarrow
W^{1,2}(B_R, \mathsf{d}^{-\nu+1}, \mathsf{d}^{-\nu+2})$. 
By the Lax-Milgram, Lemma \ref{LemmaLaxM}, we have a unique weak solution $V \in \mathfrak{E}$ of the equation
$$\mathrm{div}(\rho\mathrm{grad}V)=\mathcal{N}g \in \mathfrak{G}.$$
On the other hand we are considering
$$\mathrm{div}(\rho\mathrm{grad}G)=\mathcal{N}g. $$
Since $G\in W^{1,2}(B_R, \mathsf{d}^{\nu-1}, \mathsf{d}^{\nu})$, we can claim
$V=G^{[Z]}$.  Then by the elliptic regularity theorem, Lemma \ref{LemmaReg}, we have $G \in W^{2,2}(B_{(1-\delta)R})$ for $0<\forall \delta <1$. By the Sobolev imbedding theorem, we can claim $G \in C(\overline{B_{(1-\delta)R}})$ so that 
$G \in C(B_R)$ and
$$g=-\rho\frac{d\rho}{dP}\Big[G+4\pi\mathsf{G}\Psi\Big] \in C(B_R).
$$
This completes the proof of  i). 

Passing to the limit from
$$\|g^{(L)}\|_{L^2(B_R)}^2=
\|\sum_{|l|\leq L} g_{lm}Y_{lm}\|_{L^2(B_R)}^2=
\sum_{|l|\leq L}\|g_{lm}\|_{L^2([0,R], r^2dr)}^2,$$
we have
$$\|g\|_{L^2(B_R)}^2=\sum \|g_{lm}\|_{L^2([0,R], r^2dr)}^2.$$
Thus $g^{(L)} \rightarrow g$ in $L^2(B_R)$ and the convergence is uniformly on each compact subset of $B_R$. Thus we have proved ii). 

Since $\mathcal{K}g^{(L)} \rightarrow\mathcal{K}g$ in $C(\overline{B_R})$, we see
$$\Psi_{lm}^{(L)}(r)=-\int_0^{2\pi}\int_0^{\pi}
(\mathcal{K}g^{(L)})(r\mbox{\boldmath$\omega$})Y_{lm}(\vartheta,\phi)^*\sin\vartheta d\vartheta d\phi $$
tends to
$$\Psi_{lm}(r)=-\int_0^{2\pi}\int_0^{\pi}
(\mathcal{K}g)(r\mbox{\boldmath$\omega$})Y_{lm}(\vartheta,\phi)^*\sin\vartheta d\vartheta d\phi $$
as $L \rightarrow \infty$. A direct calculation gives
$$\Psi_{lm}^{(L)}(r)=
\begin{cases}
-\mathcal{H}_l(g_{lm})(r)\quad\mbox{if}\quad |l|\geq L \\
0 \quad\mbox{if}\quad |l| <L,
\end{cases}
$$
since
$$\mathcal{K}(g_{lm}Y_{lm})=-\mathcal{H}_l(g_{lm})Y_{lm}.$$
This implies
$$\Psi_{lm}(r)=-\mathcal{H}_l(g_{lm})(r),$$
which completes the proof of iii). $\square$\\

\begin{Remark} The regularity of $g$ given in i) of Lemma \ref{Lemma.glm} can be improved to $g \in C^{\frac12}(B_R)$  as the Sobolev imbedding theorem gives $W^{2,2}  \hookrightarrow C^{\frac12} $ and hence $G\in C^{\frac12}(B_R)$. For our purpose, the continuity of $g$ is sufficient to justify the spherical harmonics expansion of $g$. 
\end{Remark}

\begin{Lemma}\label{Lemma.Ulm}
Let $\mbox{\boldmath$\xi$} \in \mathsf{D}(\mathbf{L})$ so that $g=\mathrm{div}(\rho\mbox{\boldmath$\xi$}) \in
\mathsf{D}(\mathcal{N})$. Then we have:

i)  There exists a unique $U \in \mathfrak{E}$ and a solenoidal vector field $\hat{\mathbf{B}}
\in \hat{\mathfrak{H}}, \mathrm{div}\hat{\mathbf{B}}=0$, such that
\begin{equation}
\mbox{\boldmath$\xi$}=\mathrm{grad}U+\frac{1}{\rho}\hat{\mathbf{B}}.
\end{equation}

ii) Moreover $ U \in C(B_R)$ and
\begin{equation}
U_{lm}(r)=
\int_0^{2\pi}\int_0^{\pi}
U(r\mbox{\boldmath$\omega$})Y_{lm}(\vartheta,\phi)^*\sin\vartheta d\vartheta d\phi
\end{equation}
where
$$\mbox{\boldmath$\omega$}=\sin\vartheta\cos\phi\frac{\partial}{\partial x^1}
+\sin\vartheta\cos\phi\frac{\partial}{\partial x^2}+
\cos\phi\frac{\partial}{\partial x^3},$$
are well-defined for $0<r<R$. As $L \rightarrow \infty$,
$$U^{[L]}=\sum_{|l|\leq L}U_{lm}Y_{lm} $$
converges to $U$ in $L^2(B_R, \mathsf{d}^{\nu-1})$ 
and converges uniformly on $\overline{B_{(1-\delta)R}}$ for any fixed $\delta \in ]0,1[$.
\end{Lemma}

Proof. By the Lax-Milgram theorem, Lemma \ref{LemmaLaxM}, there exists a unique weak solution $U \in \mathfrak{E}$ of
$$\mathrm{div}(\rho\mathrm{grad}U)=g. $$
Putting $$\hat{\mathbf{B}}:=\rho\mbox{\boldmath$\xi$} -\rho\mathrm{grad}U,
$$
we have i). In fact, note $\rho\mathrm{grad}U \in \rho
L^2(B_R, \mathsf{d}^{\nu})=\hat{\mathfrak{H}}$.

By the elliptic regularity theorem, Lemma \ref{LemmaReg}, 
we have $U \in W_{\mathrm{loc}}^{2,2}(B_R)$ so that, by the Sobolev imbedding theorem, we have $U \in C(B_R)$. 
Therefore $U_{lm}$ and $\partial_r U_{lm}$ are well-defined for $0<\forall r <R$, and
passing to the limit from
$$\|U^{[L]}\|_{L^2(\mathsf{d}^{\nu-1})}^2=\sum_{|l|\leq L}
\|U_{lm}\|_{\mathfrak{X}_{\nu-1}}^2, $$
we have ii). $\square$\\

\begin{Remark}
 $U$ enjoys a higher interior regularity  $U \in C^{1,\frac12}(B_R)$. To see this, we recall that $g \in \mathfrak{G}_1 \hookrightarrow W^{1,2}(B_R, \mathsf{d}^{-\nu+1}, \mathsf{d}^{-\nu+2})$ and hence $U \in W_{\mathrm{loc}}^{3,2}(B_R)$ by the standard interior regularity. The claim follows from the Sobolev imbedding theorem. However, we note that such higher regularity of $U$ does not lead to the regularity of $\hat{\mathbf{B}}$. 
\end{Remark}

Now we have a decomposition of  $\mbox{\boldmath$\xi$}$ of the form
$$\mbox{\boldmath$\xi$}=\mathrm{grad}U+\mathbf{B},
$$
where $U$ is a scalar field and $\hat{\mathbf{B}}=\rho\mathbf{B}$ is a solenoidal vector field, that is,
$\mathrm{div}\hat{\mathbf{B}}=0$. According to \cite[p.225]{Chandra1961}, the
solenoidal field $\hat{\mathbf{B}}$ can be represented as the sum of a toroidal vector filed $\mathbf{T}$ and a poloidal vector field
$\mathbf{S}$ of the form
$$ \mathbf{T}=\mbox{curl}(T\mathbf{e}_r),
\qquad
\mathbf{S}=\mbox{curl}(\mbox{curl} (S\mathbf{e}_r)),
$$
that is,
$$\hat{\mathrm{B}}=\mbox{curl}(T\mathbf{e}_r)+
\mbox{curl}(\mbox{curl}(S\mathbf{e}_r)) ,$$
where $T$ and $S$ are scalar fields. For a mathematically rigorous proof of this decomposition, we refer to \cite{Schmidt}. Thus we have the decomposition of $\mbox{\boldmath$\xi$}$
as follows:
\begin{align}
\mbox{\boldmath$\xi$}&=\mathrm{grad}U+\frac{1}{\rho}\Big[\mbox{curl}(T\mathbf{e}_r)+
\mbox{curl}(\mbox{curl}(S\mathbf{e}_r))\Big] \nonumber \\
&=\Big(\frac{\partial U}{\partial r}-\frac{1}{r^2}\triangle_s\frac{S}{\rho}\Big)\mathbf{e}_r+
\frac{1}{r}\nabla_s\Big(U+\frac{1}{\rho}\frac{\partial S}{\partial r}\Big)+
\frac{1}{r}\nabla_s^{\perp}\frac{T}{\rho}.
\end{align}
See \cite[Appendix III]{Chandra1961}.

Let us take the expansion of the scalar fields $U, T, S$ by the spherical harmonics formally as
\begin{subequations}
\begin{align}
&U=\sum_{l,m}U_{lm}(t,r)Y_{lm}(\vartheta, \phi), \\
&T=\sum_{l,m}T_{lm}(t,r)Y_{lm}(\vartheta,\phi), \\
&S=\sum_{l,m}S_{lm}(t,r)Y_{lm}(\vartheta,\phi).
\end{align}
\end{subequations}

Let us introduce the fields
$\psi_{lm}, \chi_{lm}, \kappa_{lm}$ by
\begin{subequations}
\begin{align}
&\psi_{lm}=\frac{1}{r}\Big[\frac{\partial U_{lm}}{\partial r}+\frac{1}{r^2}l(l+1)\frac{S_{lm}}{\rho}\Big], \\
&\chi_{lm}=\frac{1}{r^2}\Big[U_{lm}+\frac{1}{\rho}\frac{\partial S_{lm}}{\partial r}\Big], \\
&\kappa_{lm}=\frac{1}{r^2}\frac{T_{lm}}{\rho}.
\end{align}
\end{subequations}

Then we have
\begin{equation}
\mbox{\boldmath$\xi$}=\sum_{l,m}\mbox{\boldmath$\xi$}_{lm}
\end{equation}
with
\begin{equation}
\mbox{\boldmath$\xi$}_{lm}=
r\psi_{lm}Y_{lm}\mathbf{e}_r+
r\chi_{lm}\nabla_sY_{lm}+
r\kappa_{lm}\nabla_s^{\perp}Y_{lm}. \label{8}
\end{equation}

\

This series expansion is merely formal at the moment. We shall not use it with effective convergence. Nonetheless, we mention that the form of an eigenvalue \eqref{xilm} of the preceding Section heuristically comes from \eqref{8}.

\section{Estimate of eigenvalues from below}

Let us consider an eigenvalue $\lambda $ of the operator $\mathbf{L}$ and let $\mbox{\boldmath$\xi$} \in \mathsf{D}(\mathbf{L})$ be an associated eigenfunction such that
$$\mathbf{L}\mbox{\boldmath$\xi$}=\lambda\mbox{\boldmath$\xi$},\quad \|\mbox{\boldmath$\xi$}\|_{\mathfrak{H}}=1.$$

Then we can put
$$\Psi=-\mathcal{K}g,\qquad g=\mathrm{div}(\rho\mbox{\boldmath$\xi$}). $$
Since $g \in \mathfrak{G}$,  we easily see that   $\Psi \in C(\mathbb{R}^3)$, using the Schwartz inequality. Moreover $\mathrm{grad}\Psi \in L^2(B_R)$. ( See 
the discussion on $\mathcal{N}_{01}$ in Section 3.) 
Note that $\triangle \Psi=g $ in distribution sense. \\

Putting $$\mathbf{C}=\rho\mbox{\boldmath$\xi$}-\mathrm{grad}\Psi,$$
we have
\begin{align*}
\mathrm{div}\mathbf{C}&=g-\mathrm{div}(\mathrm{grad}\Psi) \\
&=g-\triangle \Psi =g-g=0\quad\mbox{in distribution sense},
\end{align*}
in other words, we have the Helmholtz decomposition of $\rho\mbox{\boldmath$\xi$}$ as
$$\rho\mbox{\boldmath$\xi$}=\mathrm{grad}\Psi+\mathbf{C},\qquad \mathrm{div}\mathbf {C}=0
\quad\mbox{(in distribution sense)}. $$

Since
$$\lambda\mbox{\boldmath$\xi$}=\mathrm{grad}G=\mathrm{grad}\Big(
-\frac{1}{\rho}\frac{dP}{d\rho}g-4\pi\mathsf{G}\Psi\Big),$$
we see that $\lambda=(\lambda\mbox{\boldmath$\xi$}|\mbox{\boldmath$\xi$})_{\mathfrak{H}}$ satisfies 
\begin{equation}
\lambda =\int
\Big[\frac{1}{\rho}\frac{dP}{d\rho}|\triangle \Psi|^2-4\pi\mathsf{G}
\|\mathrm{grad}\Psi\|^2\Big]d\vec{x},
\end{equation}
where $\triangle \Psi=g$.\\

Thanks to Lemma \ref{Lemma.glm}, we may take the expansion of $g$ with respect to the spherical harmonics
\begin{equation}
g=\sum_{l,m}g_{lm}(r)Y_{lm}(\vartheta,\phi), \label{6.2}
\end{equation}
and put
\begin{equation}
\Psi_{lm}:=-\mathcal{H}_l({g}_{lm}) \label{6.3}
\end{equation}
with
$$
\mathcal{H}_l\check{g}(r)=\frac{1}{2l+1}
\Big[\int_r^{+\infty}\check{g}(r')
\Big(\frac{r}{r'}\Big)^lr'dr'+
\int_0^r\check{g}(r')
\Big(\frac{r}{r'}\Big)^{-l-1}r'dr'\Big]. \eqno\eqref{LetHlm}
$$

Then we get the expansion of $\Psi$ as 
\begin{equation}
\Psi=\sum_{l,m}\Psi_{lm}(r)Y_{lm}(\vartheta,\phi),
\end{equation}
and 
\begin{align}
&g=\triangle \Psi=\sum_{l,m}\triangle^{\langle l \rangle}\Psi_{lm}(r)Y_{lm}(\vartheta,\phi), \\
& \mbox{where,} \nonumber \\
&\triangle^{\langle l \rangle}=\frac{1}{r^2}\frac{d}{dr}r^2\frac{d}{dr}-\frac{l(l+1)}{r^2}, \\
&\mathrm{grad} \Psi=\sum_{l,m}
\frac{d\Psi_{lm}}{dr}Y_{lm}\mathbf{e}_r+\frac{\Psi_{lm}}{r}\nabla_sY_{lm}. \label{7.7}
\end{align}\\

 Of course, here and hereafter, we should read
$ \triangle^{\langle l \rangle}\Psi_{lm}={g}_{lm} \in \mathfrak{Y} $.
We point out that the well-definedness of $g_{lm}$ and the justification of the expansion 
\eqref{6.2} and the equation \eqref{6.3} are guaranteed by  Lemma \ref{Lemma.glm}.
Note that $g_{lmn} \in \mathfrak{Y}$ implies that $\Psi_{lm} \in C^1([0,+\infty[)$
by the definition \eqref{6.3} with \eqref{LetHlm}, since the derivative  $d\Psi_{lm}/dr$ is given by \eqref{Let18}.  Therefore the series expansion \eqref{7.7} is convergent in $L^2(B_R)$, while $ \triangle^{\langle l \rangle}\Psi_{lm}={g}_{lm} \in \mathfrak{Y} $
is in distribution sense. \\

Therefore we can put
\begin{equation}
\lambda=\sum_{l,m}\Lambda_{lm},
\end{equation}
with
\begin{equation}
\Lambda_{lm}=\int_0^R
\Big[\frac{1}{\rho}\frac{dP}{d\rho}|\triangle^{\langle l \rangle}\Psi_{lm}|^2-
4\pi\mathsf{G}\|\nabla^{\langle l \rangle}\Psi_{lm}\|^2\Big]r^2dr,
\end{equation}
where
\begin{equation}
\|\nabla^{\langle l \rangle}V\|^2:=\Big|\frac{dV}{dr}\Big|^2+\frac{l(l+1)}{r^2}|V|^2.
\end{equation}

Here and hereafter, by diverting the symbols,  $\displaystyle \rho, u, \frac{1}{\rho}\frac{dP}{d\rho},$ etc stands for
$\displaystyle  r \mapsto \rho(\vec{x}), u(\vec{x}), \frac{1}{\rho}\frac{dP}{d\rho}(\vec{x}),$ etc with $\|\vec{x}\|=r$.\\

Now let us consider $\Lambda_{00}$. We see
\begin{align}
\Lambda_{00}&=\int_0^R
\Big[\frac{1}{\rho}\frac{dP}{d\rho}|g_{00}|^2+4\pi\mathsf{G}\Psi_{00}g_{00}^*\Big]r^2dr \nonumber \\
&=\int_0^R(\mathcal{L}^{\mathsf{ss}}y)y^*\rho r^4dr,
\end{align}
where, taking $\displaystyle y :=\frac{1}{\rho r}\frac{d\Psi_{00}}{dr}$, 
$$g_{00}=\frac{1}{r^2}\frac{d}{d}(r^3\rho y), $$
and
\begin{equation}
\mathcal{L}^{\mathsf{ss}}y=-\frac{1}{\rho r^4}\frac{d}{dr}\Big(r^4\Gamma P\frac{dy}{dr}\Big)
-\Big(3\frac{\rho}{\Gamma}\frac{d\Gamma}{d\rho}+3\Gamma-4\Big)
\frac{1}{r}\frac{du}{dr}y
\end{equation}
with $\displaystyle \Gamma=\frac{\rho}{P}\frac{dP}{d\rho}$. 
Therefore, if we denote by $\lambda_1^{\mathsf{ss}}$ the least eigenvalue of the problem of spherically symmetric perturbations, then we have
\begin{equation}
\Lambda_{00}\geq \lambda_1^{\mathsf{ss}},
\end{equation}
provided that $y \not=0$, that is, $\Psi_{00}\not=\mbox{Const.}$.

Of course $\Lambda_{00}=0$ if $\Psi_{00}=\mbox{Const.}$
\\

 Note that, even if $\Psi=g=0$, there is an eigenfunction $\rho\mbox{\boldmath$\xi$}=\mathbf{C}\not=0, \mathrm{div}\mathbf{C}=0$, associated with the eigenvalue $0$.\\

Let us consider $\Lambda_{lm}$ with $l \geq1$. 

We see
\begin{align*}
\Lambda_{lm}&=\int\Big|
\sqrt{\frac{1}{\rho}\frac{dP}{d\rho}}\triangle^{\langle l \rangle}\Psi_{lm}+
4\pi\mathsf{G}\sqrt{\rho\frac{d\rho}{dP}}\Psi_{lm}\Big|^2r^2dr \\
&-4\pi\mathsf{G}\int\Big(
2\mathfrak{Re}[(\triangle^{\langle l \rangle}\Psi_{lm})\Psi_{lm}^*]+
\|\nabla^{\langle l \rangle}\Psi_{lm}\|^2+
4\pi\mathsf{G}\rho\frac{d\rho}{dP}|\Psi_{lm}|^2\Big)r^2dr \\
&\geq 4\pi\mathsf{G}\int
\Big(\|\nabla^{\langle l \rangle}\Psi_{lm}\|^2-
4\pi\mathsf{G}\rho\frac{d\rho}{dP}|\Psi_{lm}|^2\Big)r^2dr,
\end{align*}
since
$$\int(\triangle^{\langle l \rangle}\Psi_{lm})\Psi_{lm}^*r^2dr=-\int
\|\nabla^{\langle l \rangle}\Psi_{lm}\|^2r^2dr.
$$
This means
\begin{equation}
\Lambda_{lm}\geq4\pi\mathsf{G}\int
\Big[-\triangle^{\langle l \rangle}-4\pi\mathsf{G}\rho\frac{d\rho}{dP}\Big]\Psi_{lm}\cdot \Psi_{lm}^*r^2dr.
\end{equation}

But, differentiating
$$-\triangle u=4\pi\mathsf{G}\rho$$
with respect to $r$, we see that $\displaystyle u':=\frac{du}{dr}
\in C^2([0,R])$ is a solution of the equation
$$ Au'=0,$$
where
\begin{equation}\label{EqA}
A=-\triangle^{\langle 1 \rangle}-4\pi\mathsf{G}\frac{d\rho}{du}
=-\frac{1}{r^2}\frac{d}{dr}
r^2\frac{d}{dr}+\frac{2}{r^2}-4\pi\mathsf{G}\frac{d\rho}{du}.
\end{equation}
Note that 
$$ 0< {\rho}\frac{d\rho}{dP}=\frac{d\rho}{du}
\leq C\rho^{\gamma-2}\leq C\rho_{\mathsf{O}}^{2-\gamma}<\infty, $$
and we know
\begin{align*}
& u=u_{\mathsf{O}}-a_1r^2+O(r^4),\qquad u'=-2a_1 r+O(r^3), \\
&\frac{du'}{dr}=-2a_1 +O(r^2),\qquad \frac{d^2u'}{dr^2}=O(r)
\end{align*}
as $r \rightarrow +0$, where $a_1$ is a positive coefficient, so that
$$-\frac{2}{r}\frac{du'}{dr}+\frac{2}{r^2}u'= O(r). $$

Since $\mathcal{A}_0$ defined on $\mathsf{D}(\mathcal{A}_0)=
C_0^{\infty}(]0,R[)$ with $\mathcal{A}_0y=Ay$ is symmetric and bounded from below in the Hilbert space $\mathfrak{X}_0:=L^2([0,R];  r^2dr)$, its Friedrichs extension $\mathcal{T}_0$ is self-adjoint operator in $\mathfrak{X}_0$. The associated quadratic form is
$$
Q[y]=\int_0^R
\Big[\Big|\frac{dy}{dr}\Big|^2+
\Big(\frac{2}{r^2}-4\pi\mathsf{G}\frac{d\rho}{du}\Big)|y|^2\Big]r^2dr +\kappa\int_0^R|y|^2r^2dr,
$$
$\kappa$ being a constant such that 
$\displaystyle -4\pi\mathsf{G}\sup \frac{d\rho}{du} + \kappa \geq 1$. Since $\{ y | Q[y] \leq 1\}$ is precompact in $\mathfrak{X}_0$, we see that the spectrum of $\mathcal{T}_0$ is of the Sturm-Liouville type. 

Hereafter, diverting the symbol, we shall denote by $A$ the Friedrichs extension $\mathcal{T}_0$.

 Note that the Liouville normal form of the eigenvalue problem $Ay=\lambda y$ is
$$-\frac{d^2\hat{y}}{dr^2}+\hat{q}(r)\hat{y}=\lambda\hat{y}, $$
where
$$\hat{y}=ry,\qquad \hat{q}(r)=\frac{2}{r^2}-4\pi\mathsf{G}\frac{d\rho}{du}. $$
Since $\displaystyle \hat{q} \sim \frac{2}{r^2}, 2>\frac{3}{4}$, as $r\rightarrow +0$, the boundary point $r=0$ is of the limit-point type, (see e.g., 
\cite[p.159, Theorem X 10]{ReedS}), and, since
$$\hat{q}=\frac{2}{R^2}+O((R-r)^{\frac{2-\gamma}{\gamma-1}\wedge 1})$$ as $r\rightarrow R-0$, the boundary point $r=R$ is regular, that is, the boundary condition is the Dirichlet: $\hat{y}=0$ at $r=R$. Hence the eigenvalues are simple, and the usual Sturm-Liouville theory can be applicable.

Let $\mu_1$ be the least eigenvalue of $A$ and $\varphi_1$ be an associated eigenfunction. We may assume that
$\varphi_1 \in C^1([0,R]), \varphi_1(r) >0$ for $0<r<R$ so that $d\varphi_1/dr \leq 0$ at
$r=R-0$.  Here note that $\varphi_1(r)$ should vanish at $r=R-0$, but $u'$ does not so,
$r=R$ being a physical vacuum boundary so that $u'$ has a nonzero negative boundary value at $r=R-0$. In other words, $u' \not\in \mathsf{D}(A)$.

Anyway
$$\int_0^R(AV)V^*r^2dr \geq \mu_1\int_0^R|V|^2r^2dr \qquad \forall V \in \mathsf{D}(A).
$$
By integration by parts, we have
\begin{align*}
0=&(A u'|\varphi_1)_{\mathfrak{X}_0}= (u'|A\varphi_1)_{\mathfrak{X}_0}+
r^2u'\frac{d\varphi_1}{dr}\Big|_{r=R-0} \\
&=\mu_1(u'|\varphi_1)_{\mathfrak{X}_0}+
r^2u'\frac{d\varphi_1}{dr}\Big|_{r=R-0}.
\end{align*}
(Here $A$ denotes the differential operator in the usual meaning, but not the self-adjoint operator in the Hilbert space $\mathfrak{X}_0$. Recall that $u' \in C^2([0,R])$.) 
Since $u'(r)<0$ for $0,r\leq R$, we see $(u'|\varphi_1)_{\mathfrak{X}_0} <0$, and
we see $\displaystyle r^2u'\frac{d\varphi_1}{dr}\Big|_{r=R-0}\geq 0$.
Therefore we can claim that $\mu_1\geq 0$. As a conclusion we see that
$$\int_0^R\Big[-\triangle^{\langle 1 \rangle}-4\i\mathsf{G}\frac{d\rho}{du}\Big]V\cdot V^* r^2dr \geq 0$$
for $\forall V \in \mathsf{D}(A)$ and the equality can hold only if $\mu_1=0$ and $V \propto \varphi_1$.

 In other words,
$\Lambda_{1m}\geq 0$ and $\Lambda_{1m}=0$ only if 
$\displaystyle \Psi_{1m}=K_m\varphi_1$ with some constants $K_m, m=-1,0,1$. If
$l \geq 2$ then we have
$$-\triangle^{\langle l \rangle}=-\triangle^{\langle 1 \rangle}+\frac{l(l+1)-2}{r^2},\qquad
l(l+1)-2\geq 4, $$
so that
$$\Lambda_{lm} \geq 16\pi\mathsf{G}\int |\Psi_{lm}|^2dr.$$
Therefore $\Lambda_{lm}\geq 0$ and $\Lambda_{lm}=0$ only if $\Psi_{lm}=0$.

Summing up, we can claim that 1) if $\Psi_{00}\not=\mbox{Const.}$, then
$\lambda\geq \lambda_1^{\mathsf{ss}}$, and 2) if $\Psi_{00}=\mbox{Const.}$, then $\lambda \geq 0$ and 
$\lambda=0$ only if
$$g=-4\pi\mathsf{G}\frac{d\rho}{du}
\varphi_1\sqrt{\frac{3}{4\pi}}
\Big(K_0\zeta-\sqrt{\frac{1-\zeta^2}{2}}(K_1e^{\sqrt{-1}\phi}
-K_{-1}e^{-\sqrt{-1}\phi})\Big),
$$
since
$$P_1^0(\zeta)=\zeta, \qquad P_1^{-1}(\zeta)=-P_1^1(\zeta)=\sqrt{1-\zeta^2}.$$
Here $\zeta=x^3/r$.
Note that $g$ is axisymmetric only if $K_{-1}=K_0=K_1=0$ and $g=0$.

Summing up, we can claim the following

\begin{Theorem}\label{Th.6}
Any eigenvalue  $\lambda$  of $\mathbf{L}$ satisfies $\lambda \geq \lambda_1^{\mathsf{ss}}$ or $\lambda \geq 0$, where $\lambda_1^{\mathsf{ss}}$ stands for  the least eigenvalue of the spherically symmetric problem. Moreover, if, for the associated  eigenvector
 $\mbox{\boldmath$\xi$}$, $g=\mathrm{div}
(\rho\mbox{\boldmath$\xi$})$ is axisymmetric and $g\not=0$, then 
$\lambda \geq \lambda_1^{\mathsf{ss}}$ or $\lambda >0$.
\end{Theorem}

\begin{Corollary} \label{Cor.6.1}
If $P=\mathsf{A}\rho^{\gamma}$ and if $\gamma \geq 4/3$, then
$$\inf_{\mbox{\boldmath$\xi$} \in \mathsf{D}(\mathbf{L})}
 \frac{ (\mathbf{L}\mbox{\boldmath$\xi$}|\mbox{\boldmath$\xi$})_{\mathfrak{H}}}{\|\mbox{\boldmath$\xi$}\|_{\mathfrak{H}}} 
=0.
$$
\end{Corollary} 

Proof. See \cite[Theorem 11.7, 'Max-Min principle']{Helffer}. In fact Theorem \ref{Th.6}
implies that for the operator $\mathbf{L}$ there are no negative eigenvalues and the bottom of the essential spectrum is 0, provided that $P=\mathsf{A}\rho^{\gamma}, \gamma \geq 4/3$.

\section{ Irrotational fields} 

In this section we show that the operator $\mathbf{L}$ restricted to irrotational perturbations has the spectrum of the Sturm-Liouville type.

\begin{Definition}\label{DefIrr} 
A vector field $\mbox{\boldmath$\xi$} \in \mathfrak{H}$ is said to be irrotational if there exists a scalar field $\displaystyle U \in W^{1,2}(B_R, \rho\frac{d\rho}{dP},\rho)$ such that
$\mbox{\boldmath$\xi$}=\nabla U$ in distribution sense, that is,
$$(\mbox{\boldmath$\xi$}|\mbox{\boldmath$\varphi$})_{\mathfrak{H}} =-\int U
(\nabla\cdot(\rho\mbox{\boldmath$\varphi$}))^*dx $$
for any $\mbox{\boldmath$\varphi$} \in C_0^{\infty}(B_R; \mathbb{C}^3)$. 
\end{Definition}

Note that, if such a $U$ exists, then $U^{[Z]} \in \mathfrak{E}$ satisfies $\mbox{\boldmath$\xi$}=\nabla U^{[Z]}$. There fore we can assume $U \in \mathfrak{E}$ in the above definition without loss of generality. Moreover, such $U \in \mathfrak{E}$, if exists, is uniquely determined by 
$\mbox{\boldmath$\xi$}$. This fact follows from the Poincar\'{e} inequality, Lemma \ref{LemmaPo}.

Thanks to this  Lemma \ref{LemmaPo}, $\|U\|_{\mathfrak{E}}:=\|\nabla U\|_{L^2(\rho)} $ turns out to be an equivalent norm of 
$\mathfrak{E}$. Thus 
we put

\begin{Definition}
The set $\mathfrak{H}^{\mathsf{irr}}$ of all irrotational $\mbox{\boldmath$\xi$}=\nabla U \in \mathfrak{H}$
is the Hilbert space endowed with the norm
$\|\cdot\|_{\mathfrak{H}^{\mathsf{irr}}}$ defined by
$$
\|\mbox{\boldmath$\xi$}\|_{\mathfrak{H}^{\mathsf{irr}}}=\|\mbox{\boldmath$\xi$}\|_{\mathfrak{H}}=\|\nabla U\|_{L^2(\rho)} \simeq \|\nabla U\|_{L^2(\mathsf{d}^{\nu})}. 
$$
\end{Definition}



Now we consider 

\begin{Definition}
The set $\mathfrak{F}^{\mathsf{irr}}$ of all $\mbox{\boldmath$\xi$} \in \mathfrak{H}^{\mathsf{irr}}$
such that $g=\nabla\cdot(\rho\mbox{\boldmath$\xi$}) \in \mathfrak{G}=L^2(B_R, 
\frac{1}{\rho}\frac{dP}{d\rho})$ is a Hilbert space endowed with the norm
$\|\cdot\|_{\mathfrak{F}^{\mathsf{irr}}}$ defined by
\begin{align}
\|\mbox{\boldmath$\xi$}\|_{\mathfrak{F}^{\mathsf{irr}}}^2&=
\|\mbox{\boldmath$\xi$}\|_{\mathfrak{H}^{\mathsf{irr}}}^2+\|g\|_{\mathfrak{G}}^2 \\
&\simeq \|\nabla U\|_{L^2(\mathsf{d}^{\nu})}^2+
\|g\|_{L^2(\mathsf{d}^{-\nu+1})}^2 \nonumber
\end{align}
In other words, $\mathfrak{F}^{\mathsf{irr}}=\mathfrak{H}^{\mathsf{irr}}\cap \mathfrak{F}$. (Recall Definition \ref{Def_F}).
\end{Definition}




We claim the following

\begin{Proposition}\label{PropIrrComp}
The imbedding of $\mathfrak{F}^{\mathsf{irr}}$ into $\mathfrak{H}^{\mathsf{irr}}$ is compact.
\end{Proposition}

Proof. Let $\{ \mbox{\boldmath$\xi$} = \nabla U \} \subset \mathfrak{H}$ 
with $\{ g=\nabla\cdot \rho\nabla U \} \subset \mathfrak{G} $ be bounded.  Then
by Lemma \ref{LemmaReg} we have
$$\|\partial_j\partial_k U\|_{L^2(\mathsf{d}^{\nu+1})}\leq C.$$
On the other hand, the imbedding of $W^{1,2}(B_R, \mathsf{d}^{\nu}, \mathsf{d}^{\nu+1})$
into
$L^2(B_R, \mathsf{d}^{\nu})$ is compact. See \cite[8.8 Theorem]{GurkaO.2}, which can be applied since $\nu+1 >1$ and $$ \frac{\nu}{2} - \frac{\nu+1}{2}+1>0.$$ Therefore $\{
\mbox{\boldmath$\xi$}=\nabla U\}$ contains a convergent sequence in $L^2(B_R, \mathsf{d}^{\nu})$. $\square$\\

Now the operator $\mathfrak{A}^{\mathsf{irr}}$ defined on $\{\mbox{\boldmath$\xi$} | \mbox{\boldmath$\xi$} =\nabla U\quad\mbox{with}\quad
U \in C_0^{\infty}(B_R)\}$ with $\mathfrak{A}^{\mathsf{irr}}\mbox{\boldmath$\xi$}=\mathbf{L}\mbox{\boldmath$\xi$}=\nabla G$
 is valued in $\mathfrak{H}^{\mathsf{irr}}$. In fact, if $ U \in C_0^{\infty}(B_R)$, then $g
=\nabla\cdot(\rho \nabla U) \in C_0^{\infty}(B_R)$  and
$$G=-\frac{1}{\rho}\frac{dP}{d\rho}g+4\pi\mathsf{G}\mathcal{K}g $$
belongs to $C^2(\overline{B_R}) \subset
W^{1,2}(B_R,\mathsf{d}^{\nu-1}, \mathsf{d}^{\nu})$ so that
$\mbox{\boldmath$\xi$}=\nabla G$ is an irrotational in the sense of Definition \ref{DefIrr}.

As in Section 2.1, we see that $\mathfrak{A}^{\mathsf{irr}}$ is symmetric and bounded from below, so that it has the Friedrichs extension $\mathbf{L}^{\mathsf{irr}}$ which is a self-adjoint operator in
$\mathfrak{H}^{\mathsf{irr}}$. Proposition \ref{PropIrrComp} guarantees the following

\begin{Theorem}
The self-adjoint operator $\mathbf{L}^{\mathsf{irr}}$ in $\mathfrak{H}^{\mathsf{irr}}$ is a compact resolvent operator,  therefore its spectrum is  of the Sturm-Liouville type.
\end{Theorem}

In fact, for $\mbox{\boldmath$\xi$} =\nabla U, U\in C_0^{\infty}(B_R)$, we see
$$(\mathfrak{A}^{\mathsf{irr}}\mbox{\boldmath$\xi$}|\mbox{\boldmath$\xi$})_{\mathfrak{H}^{\mathsf{irr}}} =
\int|g|^2\frac{1}{\rho}\frac{DP}{d\rho}
-4\pi\mathsf{G}\int (\mathcal{K}g)g^*, $$
with $g=\nabla\cdot(\rho\nabla U)$. As for the quadratic form
$$Q[g]=\int|g|^2\frac{1}{\rho}\frac{dP}{d\rho}-4\pi\mathsf{G}\int
(\mathcal{K}g)g^*,
$$
we see
$$Q[g]\geq \|g\|_{\mathfrak{G}}^2-C\|\mbox{\boldmath$\xi$}
\|_{\mathfrak{H}^{\mathsf{irr}}}^2.$$
(See the proof of Proposition \ref{Prop2}.) This means
$$((\mathfrak{A}^{\mathsf{irr}}+\kappa)\mbox{\boldmath$\xi$}|\mbox{\boldmath$\xi$})_{\mathfrak{H}^{\mathsf{irr}}}
\geq \|\mbox{\boldmath$\xi$}\|_{\mathfrak{F}^{\mathsf{irr}}}^2, $$
provided that $\kappa$ is a sufficiently large constant. Then the resolvent
$(\mathbf{L}^{\mathsf{irr}} +\kappa)^{-1} \in \mathcal{B}(\mathfrak{H}^{\mathsf{irr}})$ is compact thanks to Proposition \ref{PropIrrComp}. \\

\begin{Remark}
Suppose $P=\mathsf{A}\rho^{\gamma}, \frac{6}{5}<\gamma < 2$, and consider the particular perturbation
$\mbox{\boldmath$\xi$}=r\mathbf{e}_r$.
Then $\mbox{\boldmath$\xi$}$ is irrotational, since $\mbox{\boldmath$\xi$}=\nabla U$ with
 $\displaystyle U =\frac{r^2}{2} \in W^{1,2}(B_R, \rho\frac{d\rho}{dP}, \rho)$
 and it belongs to the domain of $\mathbf{L}^{\mathsf{irr}}$ as remarked in Remark \ref{R3}, where we consider $y=1$. 
We have
$$g=\nabla\cdot(\rho\mbox{\boldmath$\xi$})=\frac{1}{r^2}\frac{d}{dr}(r^3\rho) $$
and
$$ G=(-3\gamma+3)u-r\frac{du}{dr}-4\pi\mathsf{G}
\Big[\int_0^r\rho(r')dr'-\int_0^R\rho(r')dr'\Big]. $$
Using $-\triangle u=4\pi\mathsf{G}\rho$, we can show that
$$G=(-3\gamma+4)u+C\Big(\frac{1}{\gamma-1}\Big)u(O), $$
where the constant $C(\nu)$ is given by
$$C(\nu)=-1+\int_0^{\xi_1(\nu)}\Theta(\xi;\nu)^{\nu}\xi d\xi 
= -\Big(\xi\frac{d}{d\xi}\Theta(\xi;\nu)\Big)_{\xi=\xi_1(\nu)} >0$$
with the Lane-Emden function $\Theta(\cdot; \nu)$ of index $\nu=\frac{1}{\gamma-1}$
and its zero $\xi_1(\nu)$.

Therefore, when $\gamma=4/3$, then $G$ turns out to be a positive constant, and $0$ is an eigenvalue of
$\mathbf{L}^{\mathsf{irr}}$,  $\mbox{\boldmath$\xi$}=r\mathbf{e}_r$ being  an associated eigenvector.
\end{Remark}


Moreover we claim the following 

\begin{Theorem}\label{Theorem.Kerss}
We have 
\begin{equation}
\mathrm{dim.Ker}\mathbf{L}^{\mathsf{irr}} =\mathrm{dim.Ker}\mathcal{L}^{\mathsf{ss}} \quad (\leq 1),
\end{equation}
where $\mathcal{L}^{\mathsf{ss}}$ is the self-adjoint ODE operator 
$$\mathcal{L}^{\mathsf{ss}}y=-\frac{1}{\rho r^4}\frac{d}{dr}\Big(r^4\Gamma P\frac{dy}{dr}\Big)
-\Big(3\frac{\rho}{\Gamma}\frac{d\Gamma}{d\rho}+3\Gamma-4\Big)
\frac{1}{r}\frac{du}{dr}y $$
with $\displaystyle \Gamma=\frac{\rho}{P}\frac{dP}{d\rho}$. 
\end{Theorem}

Proof.  
Let $\mbox{\boldmath$\xi$} =\nabla U \in 
\mathrm{Ker}\mathbf{L}^{\mathsf{irr}}$ with 
$g=\nabla\cdot (\rho\mbox{\boldmath$\xi$}) \in \mathsf{D}(\mathcal{N})$. Since $\nabla G=0$, $G$ is a constant. Look at
\begin{equation}
g=-\rho\frac{d\rho}{dP}\Big[G-4\pi\mathsf{G}\mathcal{K}g\Big],\quad
G=\mathrm{Const.} \label{gEX}
\end{equation}

Since $\Psi=-\mathcal{K}g \in C(\overline{B_R})$ and $\displaystyle\rho\frac{d\rho}{dP} \in C
(\overline{B_R)})$, we have
$g \in C(\overline{B_R})$ and $|g(\vec{x})| 
\lesssim \mathsf{d}(\vec{x})^{\nu-1}$. 

Let $l \geq 1$. Then $G_{lm}=0$ for $G$ is a constant. That is,
$$-\frac{1}{\rho}\frac{dP}{d\rho}g_{lm}+4\pi\mathsf{G}\mathcal{H}_l(g_{lm})=0.
$$

Consider
$$ h(r,\zeta):=\frac{1}{\rho}\frac{dP}{d\rho}g_{lm}(r)P_l(\zeta)=\frac{du}{d\rho}g_{lm}(r)P_l(\zeta), $$
where we note
$$Y_{l0}(\vartheta,\phi)=\sqrt{\frac{2l+1}{4\pi}}P_l(\cos\vartheta).$$
Then $ h \in C([0,R]\times [-1,1])$ and we have
\begin{equation}
-h+4\pi\mathsf{G}\mathcal{H}_l(g_{lm})P_l =0. \label{H1}
\end{equation}
But 
\begin{align*}
-\triangle (\mathcal{H}_l(g_{lm})P_l)&=
\Big[-\frac{1}{r^2}\frac{d}{dr}r^2\frac{d}{dr}+\frac{l(l+1)}{r^2}\Big]\mathcal{H}_l(g_{lm})\cdot P_l \\
&=g_{lm}\cdot P_l =\frac{d\rho}{du}h.
\end{align*}
Therefore \eqref{H1} means
$$-h+4\pi\mathsf{G}\mathcal{K}\Big(\frac{d\rho}{du}h\Big)=0.
$$
That is, $h$ 
 satisfies
$$-h+D\mathcal{G}(u)h=0,$$
where
$$\mathcal{G}(\varphi)=\mathsf{
G}\int_{-1}^1\int_0^{\infty}
K(r,\zeta, r',\zeta')\mathsf{f}(\varphi(r',\zeta'))r'^2dr'd\zeta'. $$
See \cite{JJTMODE}. This nothing but the antecedent of the {\bf (HL)} condition defined by 
\cite[Definition 2]{JJTMODE}. Therefore the proof of \cite[Theorem 2]{JJTMODE} says
that $h=0$ and $g_{lm}=0$ for $l\geq 2$.

We next consider $l=1$. In this case, since $G_{1m}=0$, with $4\pi \mathsf{G}\frac{d\rho}{d u} \mathcal H_{1}(g_{1m}) =g_{1m} $, we have 
\begin{equation}
A\mathcal H = 0 
\end{equation}
where
\begin{equation}
A\mathcal H= - \frac{1}{r^2} \frac{d}{dr} \left( r^2 \frac{d\mathcal H}{dr}  \right) + \frac{2}{r^2} \mathcal H - 4\pi \mathsf{G} \frac{d\rho}{d u} \mathcal H 
\end{equation}
with $\mathcal H= \mathcal H_{1}(g_{1m})$. By recalling \eqref{EqA} and discussion thereafter, we know that the operator $A$ admits its Friedrichs extension, denoted by $A$ again and the usual Sturm-Liouville theory is applicable; in particular, the least eigenvalue $\mu_1$ of $A$ is non-negative: $\mu_1\geq 0$ and the corresponding eigenfunction $\varphi_1$ can be taken such that $\varphi_1\geq 0$ for $0<r<R$, $\varphi_1(R-0)=0$, and $\varphi_1'(R-0)\leq 0$. Moreover, it holds that $ \varphi_1'(R-0) < 0$ iff $\mu_1>0$ and $\varphi_1'(R-0)= 0$ iff $\mu_1=0$. 

We claim that $\mu_1 >0$. Suppose not, namely $\mu_1=0$. Then the eigenfunction $\varphi_1$ satisfies 
\begin{equation}
- ( \varphi_1' + \frac{2\varphi_1}{r})' =4\pi \mathsf{G} \frac{d\rho}{d u} \varphi_1
\end{equation}
Integrating over $(r,R)$ and using $\varphi_1(R-0)=\varphi_1'(R-0)=  0$, we obtain 
\begin{equation}
 \varphi_1'(r) + \frac{2\varphi_1(r)}{r}  =  \int_r^R 4\pi \mathsf{G} \frac{d\rho}{d u} \varphi_1
\end{equation}
which leads to 
\begin{equation}
(r^2\varphi_1)' = r^2 \int_r^R 4\pi \mathsf{G} \frac{d\rho}{d u} \varphi_1 
\end{equation}
Since the right-hand side is positive for $r<R$, $r^2\varphi_1$ is increasing. But this implies $r^2\varphi_1(r) < 0$ because $\varphi_1(R-0)=0$, which contradicts $\varphi_1>0$. 
 
This implies that if $A\mathcal H = 0 $, then $\mathcal H = 0 $
so that $\displaystyle g_{1m}=\rho\frac{d\rho}{dP}\cdot 4\pi\mathsf{G}\mathcal{H}=0$.

Therefore, 
$$g(\vec{x})=g_{00}(r)Y_{00}(\vartheta,\phi)=\frac{1}{\sqrt{4\pi}}g_{00}(r).$$

We can solve the equation
$$\frac{1}{r^2}\frac{d}{dr}\Big(r^2\rho\frac{dV}{dr}\Big)=\frac{1}{\sqrt{4\pi}}g_{00}(r) $$
as
\begin{align*}
\frac{dV}{dr}&=\frac{1}{\sqrt{4\pi}r^2\rho(r)}\int_0^rg_{00}(r')(r')^2dr' \\
&=-\frac{1}{\sqrt{4\pi}r^2\rho(r)}\int_r^Rg_{00}(r')(r')^2dr',
\end{align*}
where we have used the fact $\int_0^R g_{00}r^2dr =0$ , which follows from $g \in \mathfrak{G}$. Then we see $|dV/dr| =O(r)$ as $r \rightarrow +0$ and
$|dV/dr|=O(1)$ as $r\rightarrow R-0$, since $|g_{00}(r)|\leq C
(R-r)^{\nu-1}$. Hence 
$$\int_0^R\Big|\frac{dV}{dr}\Big|^2(R-r)^{\nu}r^2dr < \infty. $$
Integrating $dV/dr$, we get a solution $V=O(1)$ such that
$$\int_0^R|V|^2(R-r)^{\nu-1}r^2dr <\infty. $$
Thus $V^{\flat}: \vec{x} \mapsto V(\|\vec{x}\|)$ belongs to
$W^{1,2}(B_R, \mathsf{d}^{\nu-1}, \mathsf{d}^{\nu})$ and satisfies 
$\nabla\cdot(\rho\nabla V^{\flat})=g$. By the uniqueness of the solution of Lax-Milgram theorem, we can claim $(V^{\flat})^{[Z]}=U$ for $\mbox{\boldmath$\xi$}=\nabla U
\in \mathfrak{H}^{\mathsf{irr}}, U \in \mathfrak{E}$. Since $(V^{\flat})^{[Z]}$ is spherically symmetric, $U$ is so, and $\displaystyle \mbox{\boldmath$\xi$}=\frac{dV}{dr}\mathbf{e}_r$.

This completes the proof. $\square$.\\

\hspace{15mm}

{\bf\large Acknowledgments}\\

The first author acknowledges the support from NSF Grant DMS-1608494. The second author acknowledges the support from 
JPS KAKENHI Grant Number JP18K03371. 
A part of this work was done  during the stay of the second author at Korea Institute for Advanced Study, during 23 - 27 July, 2018. The second author expresses his sincere thanks to KIAS for the hospitality and financial support for this stay.

\hspace{15mm}


\begin{thebibliography}{99}
\bibitem{Batchelor} G. K. Batchelor, An Introduction to Fluid Dynamics,
Cambridge UP, 1967.
\bibitem{BeyerSch} H. R. Beyer and B. G. Schmidt, Newtonian strellar oscillations, Astron. Astrophys., 296(1995), 722-726.
\bibitem{Beyer} H. R. Beyer, The spectrum of radial adiabatic stellar oscillations, J. Math. Phys., 36 (1995), 4815-4825.
\bibitem{BirkhoffR} G. Birkhoff and G.-C. Rota, Ordinary Differential Equations, 3rd Ed.,  John Wiley and Sons, NY, 1978.
\bibitem{Chandra} S. Chandrasekhar, An Introduction to the Study of Stellar Structure, Univ. Chicago Press, 1936.
\bibitem{Chandra1961} S. Chandrasekhar, Hydrodynamic and Hydromagnetic Stability,
Clarendon Press, Oxford, 1961.
\bibitem{Chandra64} S. Chandrasekhar, A general variational principle governing the radial and the non-radial oscillations of gaseous masses, Astrophys. J., 139(1964), 664-674
\bibitem{CH} R. Courant und D. Hilbert,
Methoden der Mathematischen Physik, Zweiter Band, Springer, 1937.
\bibitem{Davies} E. B. Davies, Spectral Theory and Differential Operators, Cambridge UP., 1995.
\bibitem{GilbargT} D. Gilbarg and N. S. Trudinger, Elliptic Partial Differential
Eqiuations of Second Order, Sprnger, 1998.
\bibitem{Gough} D. O. Gough, Liner adiabatic stellar pulsation, in J. -P. Zahn et al eds., Les Houches, Session XVII, 1987, Dynamique des fluides astrophysiques, North-Holland, Amsterdam,London, New York,Tokyo, 1993, pp. 399-560.
\bibitem{GurkaO} P. Gurka and B. Opic, Continuous and compact imbeddings of weighted Sobolev spaces. I, Czechoslovak Math. J., 38(1988), 730-744.
\bibitem{GurkaO.2} P. Gurka and B. Opic, Continuous and compact imbeddings of weighted
Sobolev spaces. II, Czechoslovak Math. J., 39(1989), 78-94.
\bibitem{Helffer} B. Helffer, Spectral Theory and Its Applications, Cambridge UP., 2013.
\bibitem{Jackson} J. D. Jackson, Classical Electrodynamics, Wiley, NY, 1962.
\bibitem{JosephL} D. D. Joseph and T. S. Lindgren, Quasilinear Dirichlet problems
driven by positive sources, Arch. Rational Mech. Anal., 49(1972/73), 241-269.
\bibitem{JJ1} Juhi Jang, Nonlinear instability theory of Lane-Emden stars,
Comm. Pure Appl. Math., LXVII(2014), 1418-1465.
\bibitem{JJTMODE} Juhi Jang and T. Makino, On rotating axisymmetric solutions of the Euler-Poisson equations, J. Differential Equations,
266(2019), 3942-3972. 
\bibitem{Kato} T. Kato, Perturbation Theory for Linear
Operators, Springer, 1980.
\bibitem{Kufner} A. Kufner, Weighted Sobolev spaces, A Wiley-Interscience Publication, John Wiley \& Sons Inc., New York, 1985.
\bibitem{Lebovitz} N. R. Lebovitz, The virial tensor and its application to self-gravitating fluids, Astrophys. J., 134(1961), 500-536.
\bibitem{LedouxW} P. Ledoux and Th. Walraven, Variable stars, in Handbuch der Physik, Band LI : Sternaufbau,  Springer, Berlin,G\"{o}ttingen, Heidelberg, 1958.
\bibitem{Lin} S.-S. Lin, Stability of gaseous stars in spherically symmetric motions, SIAM J. Math., 28 (1997), 539-569.
\bibitem{LyndenBO} D. Lynden-Bell and J. P. Ostriker, On the stability of
differentially rotating bodies, Mon. Not. Astr. Soc., 136(1967), 293-310.
\bibitem{TM1984} T. Makino, On the existence of positive solutions at infinity for ordinary differential equations of Emden type, Funkcialaj Ekvacioj, 27(1984), 319-329.
\bibitem{OJM} T. Makino, On spherically symmetric motions of a gaseous star
governed by the Euler-Poisson equations, Osaka J. Math., 52(2015), 545-580.
\bibitem{ReedS} M. Reed and B. Simon, Methods of Modern Mathematical Physics, II. Fourier Analysis, Self-adjointness, Academic Press, NY, 1975.
\bibitem{Schmidt} B. J. Schmitt, The poloidal-toroidal representation of
solenoidal fields in spherical domains, Analysis, 15(1995), 257-277.
\bibitem{Tassoul} J. -L. Tassoul, Stellar Rotation, Cambridge UP, Cambridge, 2000.



\end{thebibliography}
\end{document}